\documentclass[12pt]{amsart}

\usepackage{amssymb} 

\usepackage{enumerate, amsfonts, latexsym,epsfig, color, wasysym}
\usepackage{epstopdf}

\input xy
\xyoption{all}

\def\widecheck{\check}

\newtheorem {theorem}{Theorem} [section]
\newtheorem {lemma} [theorem] {Lemma}

\newtheorem {proposition} [theorem] {Proposition}
\newtheorem {example} [theorem] {Example}
\newtheorem {corollary} [theorem] {Corollary}
\newtheorem {coroll} [theorem] {Corollary}
\newtheorem {definition} [theorem] {Definition}

\newtheorem {remark} [theorem] {Remark}

\newtheorem {notation} [theorem] {Notation}

\newtheorem {convention} [theorem] {Convention}

\newtheorem{thmA}{Theorem}

\newtheorem{corA}[thmA]{Corollary}

\newfont{\msb}{msbm10 scaled 1200}
\newfont{\euf}{eufm10 scaled 1200}

 \def\nf {nibbled future}

\def\gep{{\rm GEP}}
\def\pep{{\rm $\Psi$EP}}
\def\rpep{{\rm PEP}}
\def\atom{atom}
\def\lbite{\LEFTcircle}
\def\rbite{\RIGHTcircle}
\def\lhnp{\Leftcircle}
\def\rhnp{\Rightcircle}

\def\n  {|\partial\Delta|}

\def\G{\Gamma}
\def\AR{\langle\Cal A\mid\Cal R\rangle}
\def\Cal {\mathcal}
\def\e{\varepsilon}

\def\At{(\mathcal A\cup\{t\})^{\pm 1}}
\def\A{\mathcal A}
\def\Apm{\mathcal A^{\pm}}
\def\FAt{(\mathcal A^{\pm}\cup\{t^{\pm 1}\})^*}
\def\FAtp{(\mathcal A^{\pm}\cup\{t^{\pm p}\})^*}
\def\FA{F(A)}

\def\Z {\mathbb Z}
\def\wt{\widetilde}

\def\time{\text{\rm time}}
\def\bonus{\text{\rm bonus}}
\def\F{\mathcal F} 
\def\Pin{\Pi} 
\def\life{\text{\rm Life}}

\def\T{\mathcal T}

\def\pTmm{\hat\T(\mu,\mu')} 
\def\pT{\hat\T} 
\def\vin        {\in_v}  

\def\vecZ{\mathcal Z} 

\def\QT{Q(\T)} 
\def\Gthree{\mathcal G_3}
\def\down{\text{\rm{down}}}

\begin{document}

\title[Mapping tori of free group automorphisms, II]
{The quadratic isoperimetric inequality
for mapping tori of free group automorphisms II:
The general case}

\author[Martin R. Bridson]{Martin R.~Bridson}
\address{Martin R.~Bridson\\
Mathematics Department\\
180 Queen's Gate\\
London, SW7 2BZ\\
U.K. }
\email{m.bridson@ic.ac.uk} 

\author[Daniel Groves]{Daniel Groves}
\address{Daniel Groves\\
Department of Mathematics \\
California Institute of Technology \\
Pasadena, CA, 91125, USA} 
\email{groves@caltech.edu} 

\date{10 October, 2006}



\thanks{The first author's work was supported in part
by Fellowships from the EPSRC and by a Royal Society-Wolfson
Research Merit Award. The second author's 
work was supported in part by a Junior Research Fellowship
at Merton College, Oxford, and by NSF Grant DMS-0504251.
We thank these organisations for their support.}

\begin{abstract} {If $F$ is a finitely generated free group and $\phi$ is 
an automorphism of $F$ then $F \rtimes_{\phi} \mathbb Z$ satisfies a
quadratic isoperimetric inequality.}
\end{abstract}

\maketitle

\section{Introduction}

This is the third and final paper in a series whose
purpose is to prove the following
theorem.

\begin{thmA} \label{MainThm} If $F$ is a finitely generated 
free group and $\phi$ is 
an automorphism of $F$ then $F \rtimes_{\phi} \mathbb Z$ satisfies a
quadratic isoperimetric inequality.
\end{thmA}

For an account of the history and context of Theorem \ref{MainThm},
we refer the reader to the introduction of \cite{BG1}.  We note here
just one additional consequence.
In \cite[Theorem 2.5]{OS}, Ol'shanskii and Sapir proved that if a 
multiple HNN extension of a free group has Dehn function less
than $n^2 \log n$ (with a somewhat
technical definition of `less than') then
it has a solvable conjugacy problem.  Theorem \ref{MainThm} shows
that free-by-cyclic groups fall into this class, and so we have
the following result.

\begin{corA} \label{Conj} 
If $F$ is a finitely generated 
free group and $\phi$ is 
an automorphism of $F$ then the
conjugacy problem for $F \rtimes_{\phi} \mathbb Z$ is solvable.
\end{corA}

Corollary \ref{Conj} was first proved in \cite{BMMV} using
different methods.

In \cite{BG1}, we proved Theorem \ref{MainThm} in the case of
{\em positive} automorphisms. That proof proceeded via an analysis
of van Kampen diagrams in the universal cover of the mapping
torus $R\times[0,1]/\langle(x,0)\sim (f(x),1)\rangle$, 
where $R$ is a 1-vertex graph with fundamental
group $F$ and  $f$ is the obvious
homotopy equivalence with $f_*=\phi$.
 
Such  $f$ are the prototypes for the 
{\em improved relative train track maps} of
Bestvina, Feighn and Handel
\cite{BFH}. Our strategy for proving Theorem \ref{MainThm} in the general
case is to refine and study these maps so as to tease-out features that
allow us to adapt the crucial arguments from \cite{BG1}. A
vital ingredient in this approach is the identification
of basic units that will play the
role in the general case that single edges (letters)
played in the positive case. We achieved this in \cite{BG2}
 with the development of {\em{beads}}, whose claim to the 
 role was clinched by the
{\em{Beaded Decomposition Theorem}}.

With this technical innovation in hand, we now set about the task of
adapting the arguments of \cite{BG1} to the general case, following
the proof  from \cite{BG1} as closely as possible and providing the
(often fierce) technical details needed to translate each step into
the more general context provided by
\cite{BG2}.  We shall not repeat the proofs of technical lemmas 
from \cite{BG1} when the adaptation is obvious. Nor shall
we repeat our account of the intuition underlying our overall strategy
of proof and intermediate strategies at  key stages. 

Unfortunately, at times we are obliged to break from the narrative
that parallels \cite{BG1} in order to deal with phenomena that do not
arise in the case of positive automorphisms ---
Section \ref{TrappedHNP}, for example. But we as far as
possible we have organised matters so that, having taken
account of the new phenomena, we can return to the main
narrative with the new phenomena controlled and
packaged into concise terminology. Thus, with
considerable technical exertions in our wake, we are
able to arrange  
matters so that  the
final stages of the proof of our main theorem
consist only of references to the corresponding
sections of \cite{BG1} with a brief explanation of what changes, if any,
must be made in the general setting.

We have already noted that,
from the analysis of improved relative train tracks in \cite{BG2}, it
emerged that beads are the correct analogue for the
role played by `letters' in the positive case. An important
manifestation of this is that
Theorem \ref{MainThm} can be reduced to a statement concerning the 
existence of a linear bound (in terms of $|\partial\Delta|$) on the
number of beads along the bottom of any corridor in a van Kampen diagram
$\Delta$ in the universal cover of the mapping tori that we consider.
In contrast to
the positive case, however, the existence of such a bound does not immediately
imply  Theorem \ref{MainThm},  because there 
is no global bound
on the length of a bead.  

Nevertheless,
proving a bound on the number of beads is
by far the bulk of our work, occupying Sections 
\ref{FastBeads}--\ref{BonusSection}, which closely 
follow \cite[Sections 6--10]{BG1} (with
different numbering and modified structure).  In Section 
\ref{LongGepsandPepsSection} we explain how the bound on the
number of beads, together with the ideas from the {\em{Bonus Scheme}}
in Section \ref{BonusSection}, finally gives Theorem \ref{MainThm}.
In Section 
\ref{BracketingSection} we explain how to deduce estimates on the
geometry of van Kampen diagrams for {\em all}
mapping tori of free group automorphism
from the specially-crafted  ones that we work
with during our main proof. The key estimate -- the linear
bound on the length
of $t$-corridors -- admits the following algebraic formulation. This
clarifies the manner in which our results concerning the
geometry of van Kampen diagrams give
rise to a non-deterministic quadratic time algorithm for the word 
problem in free-by-cyclic groups (for an alternative approach see
\cite{schleimer}).

Fix a set of generators $\mathcal B$ for $F$ and let
$d_{F}$ be the corresponding word metric.
We
consider words over the alphabet $(\mathcal B\cup\{t\})^{\pm 1}$,
where $t$ is a generator of the righthand factor of
$F\rtimes_\phi\Z$.
A {\em{bracket}} $\beta$ in a word $w$ is a decomposition
$w 
\equiv
w_1(w_2)w_3$; the subword $w_2$ is the {\em{content}} of 
$\beta$, and the 
initial and terminal letters of $w_2$ are its {\em{sentinels}}.
A  second bracket $\beta'$, giving $w \equiv w_1'
(w_2') w_3'$
is {\em{compatible}}
with $\beta$ if $w_2' \subset w_i$ for some 
$i \in\{1, 2, 3\}$ or $w_2\subset w_i'$.
A 
{\em{$t$-complete}} bracketing is a set of pairwise compatible 
brackets $\beta_1,\dots,\beta_m$ such that the sentinels of
 each $\beta_i$
are $\{t, t^{-1}\}$ and every $t^{\pm 1}$ in $w$ is a sentinel of a 
unique bracket.  In such 
a bracketing, the content of each bracket is equal in 
$F\rtimes_\phi \Z$ to an element of $F$. 

\begin{thmA}\label{Bracket}  There exists a constant $K = K(\phi,\mathcal B)$
such that any word $w \equiv e_1 \dots e_n$ that represents
the identity in $F\rtimes_\phi\Z$ 
admits a $t$-complete bracketing $\beta_1,\dots , \beta_m$
such that the content $c_i$ of each $\beta_i$
satisfies $d_F(1,c_i)\le Kn$.
\end{thmA}

In an appendix to this paper we explain how
our proof of Theorem \ref{MainThm} allows one to reprove
the main result of 
\cite{BrinkDyn}.

We suggest that readers approach this paper as follows.
First, they must be familiar with the structure of the argument
in \cite{BG1} and the vocabulary of beads in \cite{BG2}.
This will enable them to skim smoothly through Sections
\ref{DiagSection}--\ref{Section:Iterate} of the current paper.
Next, they can gain an accurate overview of the proof of Theorem 
\ref{MainThm} by reading the introduction to each of Sections 
\ref{DiagSection}--\ref{LongGepsandPepsSection} together
with the titles of their subsections (and the introductions to
subsections when they exist).
There is then no alternative but to delve into the details
of the proof.

Section \ref{BracketingSection} can be read independently.
The argument in Appendix
\ref{BrinkSection} is easy to understand in outline, but the
proof appeals to detailed results from Sections \ref{FastBeads},
\ref{TeamsSection} and \ref{BonusSection}.

\setcounter{tocdepth}{1}
\tableofcontents

\section{The Structure of Diagrams} \label{DiagSection}

Associated to any finite group-presentation $\G=\AR$ one has the
standard
combinatorial 2-complex $K(\Cal A:\Cal R)$
with fundamental group $\G$ and directed edges labelled by the
$a\in\A$.  
There is a 1-1 correspondence between words in the letters $\A^{\pm 1}$ and
combinatorial loops in the 1-skeleton of  $K(\Cal A:\Cal R)$. Words such that
$w=1$ in $\G$ correspond to loops that are null-homotopic.
Van Kampen's Lemma explains the connection\footnote{For a complete account of
the equivalences in this subsection, see \cite{steer}.} between free equalities
demonstrating the membership $w\in\langle\!\langle \Cal R\rangle\!\rangle$
and combinatorial null-homotopies for the corresponding loops. 

Such a null-homotopy is given by a van Kampen diagram over  $\AR$, which is
a 1-connected, combinatorial planar 2-complex $\Delta$ in $\mathbb R^2$ with a 
basepoint; each oriented edge is labelled by a generator $a_i^{\pm 1}$ with $a_i\in\A$
and the boundary label on each face is some $r_j^{\pm 1}$ with $r_j\in\Cal R$
(read from a suitable basepoint).  There is a  unique label-preserving map
from the 1-skeleton of $\Delta$ to the 1-skeleton of the standard
2-complex $K(\Cal A: \Cal R)$, and this extends to a combinatorial map
$\Delta\to K(\Cal A: \Cal R)$.  

Van Kampen's Lemma implies that 
the number of faces in a least-area van Kampen diagram 
with boundary label $w$ is the least number $N$ of factors among free
equalities $w=\prod_{j=1}^Nu_jr_ju_j^{-1}$.  Thus the {\em Dehn
  function} of $\AR$ can be defined to be the minimal function
$\delta(n)$ such that every null-homotopic edge-loop of length at most $n$
in $K(\Cal A: \Cal R)$ is the  restriction to $\partial\Delta$ of a
combinatorial map $\Delta\to K(\Cal A: \Cal R)$ where $\Delta$ is a
$1$-connected, planar combinatorial 2-complex.  When described in this
manner, it is natural to call  the Dehn function  the {\em
combinatorial isoperimetric function} of $K(\Cal A: \Cal R)$; the
combinatorial isoperimetric function of an arbitrary compact
combinatorial 2-complex is defined in the same way.

There is a standard  diagrammatic argument for showing that the Dehn
functions of quasi-isometric groups are $\simeq$ equivalent --- see
\cite{Alonso}. In that argument, it is unimportant that the
complexes considered have only one vertex. Thus if $K$ is any compact
combinatorial 2-complex with  fundamental group $\G$, then the
combinatorial  isoperimetric function of $K$ is $\simeq$ equivalent to
the Dehn function of $\G$. We shall exploit the freedom stemming from
this equivalence.  Specifically, we shall prove Theorem \ref{MainThm}
by establishing a quadratic upper bound on the combinatorial
isoperimetric function of a carefully-crafted 2-complex $M$ with
fundamental group $F\rtimes_{\phi^r}\mathbb Z$, where $r>0$. In other
words, we identify a constant $C>0$ such that  every null-homotopic
combinatorial loop of length at most $n$ in $M^{(1)}$ is the boundary of a
combinatorial map to $M$ from a 1-connected planar 2-complex  with  at
most $Cn^2$ 2-cells.  In fact, we prove something more refined than this (see 
Section \ref{Linear} below).

\begin{remark}\label{f.index}
Note that we are free to pass from $F\rtimes_\phi\mathbb Z$ to the
finite-index subgroup  $F \rtimes_{\phi^r}\mathbb Z$ because the
$\simeq$ class of the Dehn function of a group is an  invariant  of
commensurability.
\end{remark}

Henceforth we shall use the term {\em van Kampen diagram} to refer to the
domain of a combinatorial map to $M$ from a 1-connected planar
2-complex, with
oriented edges {\em labelled} by letters representing the oriented edges of
the target. (Note that this agrees with the standard terminology
in the special case $M=K(\Cal A:
\Cal R)$.) Such a diagram is said to be {\em least-area} if it has the
least number of 2-cells among all diagrams with the same boundary
label.

\subsection{The Mapping Torus}

Let $G$ be a compact graph and let $f:G \to G$ be a continuous
map that sends each
edge $e_i$ of $G$ to an immersed edge-path $u_i=\e_1\dots\e_m$ in $G$.  We
attach to each vertex $v\in G$ a new edge $t_v$ joining $v$ to
$f(v)$. We then attach one 2-cell to this augmented graph for each
edge $e_i$; the 2-cell is attached  along the edge path
$t_v^{-1}e_it_{v'}u_i^{-1}$, where $v$ and $v'$ are the initial
and terminal vertices of $e_i$ and where the inverse is taken in the path groupoid
(i.e. $u_i^{-1}$ is $u_i$ traversed backwards). The resulting 2-complex
is the {\em mapping torus} of $f$, which we shall denote $M(f)$.

In this paper we are primarily concerned with van Kampen diagrams over
$M(f)$, where $f$ is a homotopy equivalence representing a given
free-group automorphism $\phi$.  In this case $\pi_1(M(f)) \cong
\pi_1(G) \rtimes_{\phi} \Z$.  The $1$-cells in such a diagram
$\Delta_0$ are either labelled by some $t_u$ or by an edge $e \in G$.
We will refer to all of the edges $t_u$ as {\em $t$-edges} and, when
it does not cause confusion, denote them simply by $t$.  For the other
edges in $\Delta_0$, it is necessary to distinguish between the edge
and its {\em label} in $G$.

\begin{notation}[Labels $\widecheck\rho$]
If an edge $\e$ in a van Kampen diagram  over
$M(f)$ is labelled by an edge in $G$, then we write  $\widecheck{\e}$
to denote that label. More generally, if an edge-path $\rho$ in such a
diagram contains no $t$-edges,  we write $\widecheck{\rho}$
to denote the path in $G$ that labels $\rho$.
\end{notation}

\subsection{Time, folded $t$-corridors, singularities and bounded
  cancellation} \label{Folding}

Assume we are in the setting of the previous paragraph.
A {\em $t$-corridor} (more simply, {\em corridor}) is then defined exactly
as in \cite[Section 1.4]{BG1}, and we have the corresponding
notion of {\em{time}} (which may be thought of as a map to $\mathbb R$
that is constant on non-$t$ edges, integer-valued on vertices, and
sends the endpoints of each $t$-edge to integers that differ by $1$.
As in  \cite[Subsections 1.5, 1.6]{BG1},
 we see that each least-area
diagram is the union of its corridors, and we may assume that the tops
of all corridors are {\em folded}. (In Subsection \ref{ss:foldBead}
we shall specify how this folding is to be done, but for the
results in this subsection it is not necessary to prescribe it.)

We write $\bot(S)$ and $\top(S)$ to denote
the {\em top} and {\em bottom} of a (folded) corridor,
respectively. {\em Singularities} are defined exactly as in
\cite{BG1}.

We restrict our attention to least-area disc diagrams. The argument
used to prove \cite[Lemma 2.1]{BG1} applies {\em verbatim}  in the
present setting  to prove:

\begin{lemma}If $S$ and $S'$ are distinct corridors in a least-area
diagram, then $\bot(S)\cap\bot(S')$ consists of at most one point.
\end{lemma}

Let $L$ be the maximum length of $f(E)$ for $E$ an edge in $G$.
As in \cite[Proposition 2.3]{BG1} we have

\begin{proposition}[Bounded singularities] \label{SingularityProp}$\ $

\begin{enumerate} 
\item[1.] If the tops of two corridors in a  least-area 
diagram meet, then their intersection is a singularity. 
\item[2.] 
There exists a constant $B$ depending only 
on $\phi$  such that less than $B$ 2-cells 
hit each singularity in any  least-area diagram over $M(f)$. 
\item[3.]  
If $\Delta$ is a least-area diagram over $M(f)$, 
then there are less than $2|\partial \Delta|$ non-degenerate singularities 
in $\Delta$, and each has length at most $LB$.  
\end{enumerate} 
\end{proposition} 
\begin{proof} Except for one minor difficulty,  the proof from
\cite{BG1} translates directly to the current setting. The minor
difficulty is that in the current context the map
$f$ is a homotopy equivalence rather than a group automorphism,  
 and $f^{-1}$ is not defined as a topological map. 
Thus, given a path $\rho$,
we need a canonical path $\sigma$ in $G$ such that $f_\#(\sigma)=\rho$,
where $ f_\#$ is tightening rel endpoints.  

Consider
$\wt{M(f)}$, the universal cover of $M(f)$.  Its 1-skeleton consists of
a collection of trees (copies of the universal cover of  $G$)   joined by
$t$-edges.  Consider a lift to $\wt{M(f)}$ of the unique edge-path
$\tau_0\rho \tau_1^{-1}$ such that the $\tau_i$ are $t$-edges.
Both endpoints of this lift lie in one
of the
 trees
$T\cong \tilde G$; define
$\wt{\sigma}$ to be the unique injective path 
which joins them in $T$, and define
$\sigma$ to be the image of $\wt \sigma$ in $M(f)$. 
\end{proof}

As in \cite[Lemma 2.4]{BG1}, the above result yields as a
special case (cf. \cite{Cooper} and \cite[Lemma 2.3.1,
  pp.527--528]{BFH}):

\begin{lemma}[Bounded Cancellation Lemma] \label{BCL} There is a
  constant $B$, depending only on $f$, so that if 
$I$ is an interval consisting of $|I|$ edges 
on the bottom of a (folded) corridor $S$ in a least-area diagram over
$M(f)$, and every edge of $I$ dies in $S$, then $|I| < B$. 
\end{lemma}

\subsection{Past, Future and Colour in Diagrams}\label{s:colour}

These concepts, for edges and 2-cells in van Kampen diagrams $\Delta$,
are defined exactly as in \cite[Section 3]{BG1}.
The {\em immediate past} (or {\em ancestor}) of an edge at the top of a
corridor in any diagram
 is the unique edge at the bottom of the corridor 
that lies in
the same 2-cell; the {\em entire past} of an edge is defined by taking the
transitive closure of the relation ``is the immediate past of''. The
past of a 2-cell is defined similarly. The {\em future} of an edge $e_0$ is
the set of edges that have $e_0$ in their past. The future of 2-cells
is defined similarly. The evolution of edges is
described by a graph $\Cal F$  whose vertices are the 1-cells $e$ of $\Delta$,
which has an edge connecting each $e$ to its  immediate ancestor.
Note that $\mathcal F$ is a forest.
Its connected components  define  {\em colours} in $\Delta$;
each edge not labelled $t$ is
assigned a unique colour, as is each 2-cell. Note that colours are in bijection with
  a subset of the edges of the boundary of the
diagram. 
The union of the 2-cells in a corridor $S$ that have colour
$\mu$ will be denoted $\mu(S)$.

As in \cite{BG1},  simple separation arguments yield the following
observations.

\begin{lemma} Each $\mu(S)$ is
 connected and intersects each of $\top(S)$ and $\bot(S)$ in an interval.
\end{lemma}

\begin{lemma} [cf. Lemma 5.9, \cite{BG1}]
Let $\e_1, \e_2$ and $\e_3$ be three (not necessarily adjacent) edges
that appear in order of increasing subscript as one reads from left to
right along the bottom of a corridor.  If the future of $\e_2$
contains an edge of $\partial \Delta$ or of a singularity, then no 
edge in the future of $\e_1$ can cancel with any edge in the future of
$\e_3$. 
\end{lemma}

Again following \cite{BG1}, given a diagram $\Delta$ we define $\Cal Z$ to
be the set of pairs $(\mu,\mu')$ such that the coloured regions
$\mu(S)$ and $\mu'(S)$ are adjacent in some corridor $S$. The proof of
\cite[Lemma 6.3]{BG1} establishes: 

\begin{lemma} \label{NoOfAdj}
$$|\Cal Z| \le 2\,|\partial\Delta| - 3.$$
\end{lemma}

\section{Adapting Diagrams to the Beaded Decomposition}

We refer the reader to \cite{BG2} for the definitions and results
which we require here about improved relative train track maps,
nibbled futures, monochromatic paths,
hard splittings and the language of {\em{beads}}
 --- including
$(J,f)$-atoms, {\gep}s and {\pep}s and what it means for a path
to be $(J,f)$-beaded.
We shall proceed under the assumption that the reader is
familiar with each of these terms, and work axiomatically with
the following outputs from \cite{BG2}.

\begin{theorem}[Beaded Decomposition Theorem, \cite{BG2}] \label{BDT}
For every $\phi\in{\rm{Out}}(F_r)$,
there exist positive integers $k$, $r$ and $J$ such that $\phi^k$ has an improved
relative train-track representative $f_0:G\to G$ 
with the property that every $(f_0)_{\#}^r$-monochromatic
path in $G$ is $(J,f_0)$-beaded.
\end{theorem}

Beads are either monochromatic paths (in case they are \atom s) or else \gep s
or \pep s (which may be monochromatic, but do not have to be). 
Thus, by the above theorem and \cite[Proposition 6.10]{BG2}, any nibbled
future of a $(J,f_0)$-bead is $(J,f_0)$-beaded. Any hard splitting of an edge-path is 
inherited by its (nibbled) futures, by definition. And if one refines a hard splitting
by decomposing the factors in a hard splitting, the result is again
a hard splitting (\cite[Lemma 2.6]{BG2}). Thus we have:

\begin{coroll} \cite[Theorem 8.4]{BG2} 
Let $f = (f_0)_{\#}^r$ be as in the Beaded Decomposition Theorem above. If an edge-path $\sigma$ in $G$
is $(J,f_0)$-beaded, then any $f$-nibbled future of $\sigma$ is $(J,f_0)$-beaded.
In particular, $f_{\#}(\sigma)$ is also $(J,f_0)$-beaded.
\end{coroll}

\begin{remark} An important point to recall from \cite{BG2} is that the
decomposition of an edge-path into $(J,f_0)$-beads is canonical.
\end{remark}

The value of the constant $J$ in the Beaded Decomposition Theorem will be
of no importance in what follows, so we drop it from the terminology.
Similarly, we will fix the map $f_0$.  Once we have passed to the power $f = (f_0)_{\#}^r$, the above results remain true when $f$ is replaced by an iterate.  Therefore, we
refer simply to  ``beads" and ``beaded paths".

\subsection{Refolding corridors according to the Beaded Decomposition}\label{ss:foldBead}

\ 
\medskip

{\em{Henceforth\footnote{There exceptions
to this in Theorem \ref{BoundS}, Section \ref{BracketingSection} and Appendix 
\ref{BrinkSection}}, we consider only diagrams over the mapping torus of $M(f)$, where
$f$ is an iterate of $(f_0)_{\#}^r$ as in the Beaded Decomposition Theorem.} In Section \ref{Section:Iterate}, we will fix the map $f$ once and for all.}

\medskip

We return to the matter of how best to fold the tops of corridors in
least area diagrams over $M(f)$. Given an arbitrary least-area
diagram, we refold the tops
of corridors in order of increasing time. The process begins with
edges at the minimal time on the boundary of the diagram, where
there is no folding to be done provided the boundary label is
reduced.

Focussing on a
particular corridor $S$, our folding up to  ${\rm{time}}(S)$
defines the histories of all edges up to this time and hence
assigns colours to the edges on $\bot(S)$, decomposing
it as a concatenation of monochromatic paths, one for
each of the colours $\mu(S)$. Theorem \ref{BDT} decomposes
each of these labels as a hard splitting of beads $\sigma_i$.
The hardness of the splitting means that after tightening the
$f(\sigma_i)$, their concatenation will be a tightening of 
$f_{\#}(\check\mu(S))$. We insist that the first step in 
the tightening of the naive top of $S$, is that determined by
the tightening of labels just described: i.e.~we first tighten
beads {\em{within}} colours, each according to a left-to-right
convention (which labels inherit
from the orientation of the corridors within
the diagram). Then, as a second step, we tighten (again with
a left-to-right convention) the concatenation of the tightened
images of the colours.  A diagram which is folded according to
these conventions will be called {\em well-folded}.

The key point of this convention is that the hard splitting 
of the label on each colour is carried into the future --- of course the futures of
the original beads may split into a concatenation of several
beads, and some beads at the ends of each colour may be cancelled by interaction 
with neighbouring colours, but {\em{each bead (more precisely\footnote{we shall
generally drop this cumbersome distinction in the sequel},
bead-labelled arc) in the beaded decomposition
of each coloured interval on $\top(S)$ is contained into the future of a 
unique bead-labelled arc of the same colour on $\bot(S)$.}} Thus 
$\top(S)$ is a concatenation of beads, each with a definite colour,
where neighbouring beads are separated by a hard splitting if
they are of the same colour but perhaps not if they are of a 
different colour.
 (It also becomes sensible to discuss
the future of a bead in a [well-folded] diagram.)

\medskip
{\em{
We henceforth suppose (usually without comment) that our diagram 
has been refolded according to this convention.}}
\smallskip
 
\begin{definition}\label{BeadLength} [cf.~Definition \ref{beadNorm}]
The {\em bead length} of $[S]_\beta$, of a corridor $S$ in a well-folded
diagram is the number of beads along $\bot(S)$. 
\end{definition}

\begin{remark} 
It is important to note that the decomposition of $\bot(S)$
and $\top(S)$ into coloured intervals is
{\em not} a hard splitting in general.  Indeed it
is the analysis of 
the cancellation between these intervals as one
flows $S$ forwards in time  that forms the meat of
this paper.
\end{remark}

\subsection{Abstract Futures of Beads}\label{stackDiags}
Given an edge-path $\rho$ in $G$, expressed as a concatenation
of monochromatic edge-paths $\rho=\rho_1\dots\rho_m$, consider
the van Kampen diagram $\Delta(l,\rho)$ with boundary
label equal to $t^{-l}\rho t^l \overline{f_{\#}^l(\rho)}$; 
this is a simple stack of corridors. The above convention dictates
how we should fold the corridors of $\Delta$ and
determines the future at each time up to $l$ for each  bead in the
beaded decompositions of the $\rho_i$. 
 
We define the (full) {\em{abstract future of a bead in $\rho$}}   to be
(the label on) its future in $\Delta(l,\rho)$.

\section{Linear Bounds on the Length of Corridors} \label{Linear}

In any least-area diagram, 
each corridor has at least two edges on the boundary, namely its $t$-edges.
The {\em length} of a  
corridor $S$ is defined to be the number of 2-cells that it 
contains. The area of a least-area diagram is the sum of the lengths of
its corridors, and therefore Theorem \ref{MainThm}
 is an immediate consequence of:

\begin{theorem} \label{BoundS}
Let $\phi$ be an automorphism of a finitely generated free group
and let $f$ be a topological representative for a positive power of $\phi$.  
There is a 
constant $K$, depending only on $f$, so that each corridor in a least-area 
diagram $\Delta$ over $M(f)$ has length at most $K\n$.
\end{theorem}

Note that Theorem \ref{MainThm} actually depends only on establishing 
Theorem \ref{BoundS}
 for a single  topological representative $f^k$ of a suitable
 power of our given free group automorphism
$\phi$; in the next section we shall articulate what that suitable
power is.
The bulk of this paper will then be devoted to 
 proving the existence of the constant $K$ for this
particular $f^k$. (In Section \ref{BracketingSection} we shall 
deduce Theorem \ref{BoundS} from this special case.)

Having restricted attention to a particular $f^k$, we may
further restrict our attention to diagrams that are well-folded
in the sense of Subsection \ref{ss:foldBead}, since
refolding the corridors of an arbitrary a diagram does
not change the configuration of corridors or their length.
In a well-folded diagram,
the top of each corridor $S$ is a concatenation of beads, and
the vast majority of our work (up to and including Section \ref{BonusSection}) goes into proving  the following result.

\begin{theorem}\label{noName} If $f$ and $k$ are as  above, then
there is a constant $K_1$ such that  all
corridors $S$ in  well-folded, least-area diagrams $\Delta$ over  $M(f_{\#}^k)$,
have bead length $[S]_\beta\le K_1\ |\partial\Delta|$. 
\end{theorem}

The linear bound on the length of $S$ that we require for
Theorem \ref{BoundS} does not follow directly from this estimate because
there is no uniform bound on the length of certain beads, namely
\gep s and \pep s. However, we shall see in Section \ref{LongGepsandPepsSection} 
that the
ideas developed in \cite{BG1} to implement the bonus scheme adapt
 to the current setting to provide the following estimate:

\begin{proposition}\label{longGEPS} There are constants $J$ and $K_2$,
depending only on $f$, such that the beads $\beta$ on $\bot(S)$ of length
greater than $J$ satisfy
$$
\sum_{\beta} |\beta| \le K_2\, |\partial\Delta|.
$$
\end{proposition} 
The constant $J$ in the above statement is the one from Theorem \ref{BDT}.

\section{Replacing $f$ by a Suitable Iterate} \label{Section:Iterate}

In order to establish the  bound on the 
length of corridors required to prove Theorem \ref{BoundS}, 
we must analyse how 
corridors grow as they flow into the future and 
assess what cancellation can take place to inhibit this 
growth.  This is much more difficult than in \cite{BG1} because
now we must cope with the cancellation that takes place
within colours. But in common with our approach in \cite{BG1},
we can appeal to Remark \ref{f.index}
 repeatedly in order to replace our
topological representative $f$ by some iterate of $f$ that
affords a more stable situation in which 
cancellation phenomena are more amenable to analysis.

In the present setting, we have to be a little careful about specifying
what we mean by ``an iterate", because we wish to consider only
topological  representatives  whose restriction to
each edge is an immersion, and this property is not inherited by
powers of the map. To avoid this problem, we deem the phrase\footnote{and obvious
variations on it}
{\em replacing $f$ by an iterate}, to mean that for fixed
$k\in \mathbb N$, we pass from consideration of $f : G \to G$
to consideration of  the map $f_{\#}^k : G \to G$ that sends
each edge $E$ in $G$ to  the tight edge-path $f_{\#}^k(E)$ that is 
homotopic rel endpoints to $f^k(E)$.

When we replace $f$ by $f_{\#}^k$, we leave behind the   
mapping torus $M(f)$ and consider instead $M(f_{\#}^k)$, which although 
homotopic to a $k$-sheeted covering of  $M(f)$ is distinct
from it. 

A corridor in a van Kampen diagram over $M(f_{\#}^k)$ can be divided into
a stack of $k$ corridors in order to yield a van Kampen diagram over 
$M(f)$. This observation will play little role in our arguments, but it highlights
one reason for hoping to simplify diagrams by passing to an iterate of $f$: the van Kampen
diagrams over $M(f_{\#}^k)$ are a proper subset (after subdivision\footnote{the obvious subdivision of a diagram $\Delta$ is called the {\em{$k$-refinement}}}
of $\Delta$) of 
the diagrams over $M(f)$; in the diagrams of this subset, corridors  
flow unhindered for at least $k$ steps in time.

\subsection{Finding the desired iterate} \label{Iterate}

We have already passed to a large iterate in order to obtain
the Beaded Decomposition Theorem. In the present subsection we
pass to further iterates in order to control the behaviour of
the images of beads.

Before settling on a specific $f$ for the remainder of the
paper, we must remove an irritating ambiguity concerning the ordering
of strata in the filtration associated to the train track structure.
This is required in order to render the choices in Section 
\ref{s:PrefFuture} coherent.

\begin{definition} \label{interchangeable}
Suppose that $f : G \to G$ is an improved relative train track map, and
that $H_i, H_j$ are strata for $f$.  We say that $H_i$ and $H_j$ are
{\em interchangeable} if one can reorder the strata, so that one
still has an improved relative train track structure, but
the order of $H_i$ and $H_j$ is reversed.
\end{definition}

If $H_i$ and $H_j$ are interchangeable, and $i > j$, then no iterate
of any edge in $H_i$ crosses an edge in $H_j$ (and neither do the
iterates of any edges occurring in the iterated images of edges in $H_i$).

\begin{convention} \label{convention:interchangeable}
We suppose that for any improved relative train track map that
we consider, if $H_i$ and $H_j$ are interchangeable strata so
that $H_i$ is an exponential stratum and $H_j$ is a parabolic
stratum then $i > j$.  

We further assume that if $H_i
= \{ E_i \}$ and $H_j = \{ E_j \}$ are interchangeable
 parabolic strata  and $n\mapsto |f^n(E_i)|$
grows exponentially while $n\mapsto |f^n(E_j)|$
grows polynomially, then $i > j$.  And if both these functions
grow polynomially, then the degree of polynomial growth of
the former is at least as great as the latter. 
\end{convention}

In the following lemma, $\omega$ is the number of strata
in the train track structure for $f$. Also recall that an
edge $\e$ in a path $\sigma$ is said to be {\em displayed} if there is a hard
splitting $\sigma = \sigma_1\odot \e \odot \sigma_2$.  The definition of
a {\em displayed sub edge-path} is entirely analogous, and will be used later.

\begin{lemma} \label{fofatom}
One can replace $f$ by an iterate to ensure  that
 if $\rho$ is any atom then
either the beads of  $f_{\#}^\omega (\rho)$ are
Nielsen paths and \gep s only, or else there is a displayed edge
$\epsilon $ in $ f_{\#}^\omega (\rho)$ so that 
\begin{enumerate}
\item $\e$ is 
 of highest weight amongst all displayed
edges in all $f_{\#}^k(\rho)$, for $k \ge 1$, and 
\item the growth of $n\mapsto |f^n_{\#}(\epsilon)|$ is at least 
as large as that of
any displayed edge in any $f_{\#}^k(\rho)$.
\end{enumerate}
\end{lemma}

\begin{proof} Lemma \ref{fofatom} from \cite{BG2} contains
all but statement (2), whose validity is assured by
 Convention 
\ref{convention:interchangeable}.
\end{proof}

Our next two results capture the {\em{end
stability}} that \cite[Proposition 4.5]{BG1} provided in the
case of positive automorphisms. This is the first stage in our
analysis at which we encounter an awkward point that does not
arise in \cite{BG1}, namely   there may exist beads (more
specifically atoms) $\rho$ such that $f_{\#}(\rho)$ is a single vertex.

\begin{definition} \label{vanishing} A {\em vanishing bead (atom)}
$\rho$ is one with $f_{\#}(\rho)$  a single vertex. 
\end{definition}

\begin{lemma} \label{EndStab1}
There exists a constant $k_0$, depending only on $f$ so that the
map $f_0 = f_{\#}^{k_0}$ satisfies the following properties.  Let
$\rho$ be a non-vanishing bead,
let $i \in \{ 1 ,\ldots , \omega \}$, and let  $\sigma_i$ be the 
leftmost bead in $(f_0)_{\#}(\rho)$ of weight at least $i$.
\begin{enumerate}
\item If $\sigma_i$ is not a \gep\ or a \pep\ then the leftmost bead
of weight at least $i$ in $(f_0)_{\#}^j(\rho)$ is the same for all $j \ge 1$.
Furthermore, in this case $\sigma_i$ is a single (displayed) edge or
a Nielsen bead.
\item If $\sigma_i$ is a \gep\ or a \pep\ then the leftmost bead of weight
at least $i$ in $(f_0)_{\#}^j(\rho)$ is contained in the (abstract)
 future of $\sigma_i$
for all $j \ge 1$.
\end{enumerate}
\end{lemma}
\begin{proof}
If $\sigma$ is a bead then
all iterated images of $\sigma$ are beaded paths, 
and a simple finiteness
argument shows that there is a
bound on the number of beads which are not
\gep s or \pep s.
\end{proof}

An entirely similar argument applies to rightmost beads, of course.
In order to deal with the different types of beads, we also
need the following variant.

\begin{lemma} \label{EndStab2}
There exists a constant $k_1$, depending only on $f$, so that
the map $f_1 = f_{\#}^{k_1}$ satisfies the following properties.   Let
$\rho$ be a non-vanishing bead 
and let $\sigma$ be the leftmost bead in 
$(f_1)^j_{\#}(\rho)$ which is not a Nielsen bead.
\begin{enumerate}
\item  If $\sigma$ is not a \gep \ or a \pep\ then for all $j \ge 1$
the leftmost bead in $(f_1)^j_{\#}(\rho)$ which is not a Nielsen bead
is $\sigma$.  Furthermore, in this case $\sigma$ is a (displayed)
edge.
\item  If $\sigma$ is a \gep\ or a \pep\ then for all $j \ge 1$
the leftmost bead in $(f_1)_{\#}(\rho)$ which is not a Nielsen bead
is in the future of $\sigma$.
\end{enumerate}
\end{lemma}

We are finally in a position to articulate all of the properties
that we want to arrange for $f$ by replacing it with an iterate.

\begin{proposition} \label{FinalPower} There is a constant $D_2$
that depends only on $f$, so that if we replace $f$ by 
 $f_{\#}^{D_2}$ then,
\begin{enumerate}
\item\label{ExpLong}  the conclusion of \cite[Lemma 5.1]{BG2} holds with $k_1 = 1$: in particular, if $\e$ is an exponential edge of weight $i$, then $f(\e)$
is longer than the unique indivisible Nielsen path of weight $i$
(if it exists);
\item the conclusion of \cite[Theorem 8.1]{BG2} holds with $D_1 = 1$;
\item the conclusion of Lemma \ref{fofatom} holds;
\item the conclusions of Lemmas \ref{EndStab1} and \ref{EndStab2} hold; and
\item\label{DisplayExp} if $\rho$ is a bead  then
$f_{\#}(\rho)$ contains at least three displayed copies of any
exponential edge that is displayed 
in any  $f_{\#}^j(\rho),\ j\ge 1$. Moreover, the leftmost (and rightmost)
such displayed edge $\e$ is contained in a displayed path
of the form $f(\e)$.
\end{enumerate}
\end{proposition}

\medskip
\noindent{\bf{Power Decree:}}{\em{
For the remainder of the paper, we will assume
that $f : G \to G$ is an improved
relative train track map that satisfies the properties 
in Proposition \ref{FinalPower}. We shall also
operate under Convention \ref{convention:interchangeable}.
}}

{\em Let $L$ be the maximal length of $f(E)$, for
edges $E \in G$.}\footnote{In \cite{BG1}, the symbol
`$M$' was used for the analogous quantity.  We use $L$ here
(and in \cite{BG2}) in 
order to avoid
confusion with the mapping torus $M(f)$.}

\section{Preferred Futures of Beads} \label{s:PrefFuture}

The reader who is comparing our progress to \cite{BG1}
will find that we are now in the position that we were at the
start of Section 5 of that paper. Thus we now want to define
the preferred future of a bead $\rho$ (in three senses\footnote{in $f_{\#}(\rho)$,
in a diagram, and in a concatenation of beaded paths})
 and then begin a study of  fast beads.

Unfortunately, the definition of the preferred future of a bead 
in a diagram is much more 
cumbersome than the analogue in \cite{BG1}.

\subsection{Abstract Preferred Futures and Growth}
 
First we note that if beads (or more generally edge paths
in $G$) are ever going to vanish in the sense
of Definition \ref{vanishing}, then they do so immediately.

\begin{lemma}
If $\sigma$ is an edge path in $G$ and
$f_{\#}^k(\sigma)$ is a vertex for some $k\ge 1$,
then $f_{\#}(\sigma)$ is already a vertex.
\end{lemma}
\begin{proof}
For all vertices $v \in G$, $f(v)$ is a fixed point of $f$.  Therefore,
the endpoints of $f_{\#}^j(\sigma)$ are the same for all $j \ge 1$.
If $f^k_{\#}(\sigma)$ is a point, then the endpoints of $f_{\#}^k(\sigma)$
are equal, hence the tight path $f_{\#}(\sigma)$ is  a loop. Since $f$ is a 
homotopy equivalence, this loop must be trivial.
\end{proof}

\begin{definition} [Abstract preferred futures]\label{PrefFuture}
The (immediate) {\em preferred future} of a 
non-vanishing bead $\sigma$ is a particular
 bead in the beaded decomposition of $f_{\#}(\sigma)$, as
defined below. The {\em $k$-step preferred future} is then
defined by an obvious recursion.
\begin{enumerate}
\item If $\sigma$ is a \gep\  then $f_{\#}(\sigma)$ is
also a \gep, and we define the preferred future of $\sigma$
to be $f_{\#}(\sigma)$.
\item If $\sigma$ is a \pep\ then either $\sigma$ or $\overline{\sigma}$
has the form $\sigma = E \overline{\tau}^k\nu\gamma$.
If it is $\sigma$, then
by \cite[Corollary 6.11]{BG2}, $f_{\#}(\sigma)$ is
either of the form $\sigma' \odot \xi$, where $\sigma'$
is a \pep\ (which has the same weight as $\sigma$), 
or else of the form $E \odot \xi$, where $E$
has the same weight as $\sigma$ and is the unique highest
weight edge in $f_{\#}(\sigma)$.  In the first case, the preferred
future of $\sigma$ is $\sigma'$.  In the second case, the 
preferred future of $\sigma$ is $E$. The preferred future of a \pep \ $\sigma$
where $\overline{\sigma}$ has the above form is
defined in an entirely analogous way.
\item If $\sigma$ is a Nielsen path then 
the preferred future of
$\sigma$ is $f_{\#}(\sigma) = \sigma$. 
\item Finally, we consider  a non-vanishing
\atom\ $\sigma$.

\noindent
(a) If the beaded decomposition of
 $f_{\#}(\sigma)$  consists entirely of Nielsen paths
and \gep s, 
 then we fix a highest weight \gep\ to
be the preferred future of $\sigma$;  
otherwise, we fix a highest weight
Nielsen path.

\noindent
(b) If not, then let $\e$ be the edge described in Lemma \ref{fofatom},
 fix a displayed occurrence of $\e$ in $f_{\#}(\sigma)$ (in case $\e$
 is exponential, choose a displayed occurrence that is neither
leftmost nor rightmost\footnote{this exists by Proposition \ref{FinalPower}}) 
and define this to be the preferred future of
$\e$.
\end{enumerate}
\end{definition}

\begin{remark}
Suppose that $\e$ is an edge in $G$, considered as a bead, and suppose
that $\e$ is not contained in a zero-stratum.  Then $\e$ has a preferred
future, which is an edge contained in the same stratum as $\e$.  We always
assume that the preferred future of $\e$ is a (fixed) occurrence of $\e$ in 
$f_{\#}(\e)$ which satisfies the requirements of the above definition.  This
situation is very close in spirit to the definition of preferred future in \cite{BG1}.
\end{remark}

We now divide the beads into classes according to the growth
of the paths $f_{\#}^k(\sigma)$, $k = 1,2,\ldots$.  Specifically,
we define left-fast and left-slow beads in accordance with
\cite[Subsection 5.1]{BG1}.

\begin{definition} [Left-fast beads] \label{LeftFast}
\gep s and Nielsen paths are {\em left-slow}.

Suppose that $\alpha$ is an \atom\ or a \pep.  Then $\alpha$ 
is {\em left-fast} if the distance between the left end of 
$f_{\#}^k(\alpha)$ and the left end of the 
preferred future of $\alpha$ in 
$f_{\#}^k(\alpha)$ grows at least quadratically with $k$, and {\em left-slow}
otherwise.
\end{definition}
Note that if a \pep\ $\sigma$ is left-fast then it is $\overline{\sigma}$
which it is of the form $E \overline{\tau}^k \nu\gamma$.

\begin{remark}
We only care that fast growth be super-linear, but it happens that
this is the same as being at least quadratic  (cf.~\cite{BGgrowth}).
\end{remark}

The concepts of {\em right-fast} and {\em right-slow} beads are
entirely analogous.

\subsection{Preferred future in diagrams} \label{PrefFutureDiag}

In this subsection we define the notion of `preferred futures' within
van Kampen diagrams.  We also
define `biting' and `consumption', which are the analogues in this
paper of `consumption' from \cite[Section 5]{BG1}.

The folding convention of Subsection \ref{Folding} expresses
$\bot(S)$ as the concatenation of coloured paths $\mu(S)$, each labelled
by a monochromatic path in $G$. The Beaded Decomposition Theorem
gives us a hard splitting into beads
\[      \widecheck{\mu(S)} = \widecheck{\beta}_1 \odot 
\widecheck{\beta}_2 
\odot \cdots \odot \widecheck{\beta}_{m_\mu}    ,      \]
and it is convenient to refer to the sub-paths 
$\beta_i \subseteq \bot(S)$ carrying the labels $\widecheck\beta_i$
as beads, as we did in Subsection \ref{Folding}.

If $\mu_1, \ldots , \mu_k$ are the colours appearing in $S$, in order,
then the label on $\top(S)$ is obtained by tightening
\[ f_{\#}(\widecheck{\mu_1(S)}) \cdots f_{\#}(\widecheck{\mu_k(S)}).    \]
The path $f_{\#}(\widecheck{\mu_1(S)}) \cdots 
f_{\#}(\widecheck{\mu_k(S)})$
is called the {\em semi-naive future} of $S$.

We have adopted a left-to-right convention to remove any
ambiguity in how one tightens the semi-naive
future to obtain the label of $\top(S)$.

We previously defined the (immediate) future of a 
bead $\beta \subset \bot(S)$ to consist of those edges of $\top(S)$ 
whose immediate past lies  in $\beta$. Since it is integral to what
we shall do now, we re-emphasize:

\begin{lemma}
The immediate future of a bead $\beta\subset\bot(S)$ is a (possibly
empty) interval equipped with a hard-splitting into beads.
\end{lemma}

If $\rho$ is the immediate
future of $\beta$, then $\rho$ is also an interval in the semi-naive
future of $S$, and hence its label
$\widecheck{\rho}$ is a specific sub-path of 
$f_{\#}(\widecheck{\beta})$.   [Note that one has more than the
path $\widecheck{\rho}$ here, one also has its position within
$f_{\#}(\widecheck{\beta})$; thus, for example, 
we would distinguish between the two visible copies
of $\widecheck{\rho}$ in  $f_{\#}(\widecheck{\beta})=
\widecheck{\rho}\sigma\widecheck{\rho}$.]

\begin{definition}[Preferred and tenuous futures in $\Delta$] \label{PreferredFutureDef} 
Consider a bead $\beta\subset\mu(S)\subset\bot(S)$ in $\Delta$
whose immediate future $\rho\subset\top(S)$ determines the
subpath $\widecheck{\rho_0}$ of $\widecheck{\beta}$ in $G$.

If the (abstract) preferred future $\widecheck\beta_+$ of  $\widecheck\beta$, as
defined in Definition \ref{PrefFuture}, is entirely contained in  
$\widecheck{\rho_0}$, then the corresponding sub-path
$\beta_+$ of $\rho$ 
is the {\em preferred future} of $\beta$. 

If $\widecheck{\rho_0}$ does not contain $\widecheck\beta_+$,
 then $\beta$ does not have a preferred future. In this situation
we say that the future of $\beta$ is {\em{tenuous}}.
\end{definition}

\begin{remark} \label{rem:Pref}
Note that, if it exists, the preferred future of a bead $\beta\subset\mu(S)$ is 
a bead in the beaded decomposition of both $\rho$ and the $\mu$-coloured
interval of $\top(S)$.

Also, if a bead happens to be a single edge  $\e$ whose label is not contained in 
a zero stratum, the preferred future is a single (displayed) edge, with 
the same label as $\e$.
\end{remark}

\begin{definition} [Biting and consumption] If the future of a bead
$\beta\subset\bot(S)$ is tenuous,  
we say that $\beta$ is {\em bitten in $S$}. If, in the notation of
(\ref{PreferredFutureDef}),  {\em no edge} of the preferred future of 
$\widecheck{\beta}$ appears in $\widecheck{\rho}$, then we say that
$\beta$ is {\em consumed} in $S$.
\end{definition}

\begin{remark}
The above definition says in particular that any bead whose label
is a vanishing atom is consumed.
\end{remark}

Let $\beta'\subset\bot(S)$ be a bead whose label is non-vanishing.
If $\beta'$ is bitten in $S$, 
there is a specific edge $\e$ in the semi-naive future of $S$ 
that, during the tightening process, is the first to  cancel with
an edge $\e'$ in the interval labelled by the preferred future of 
$\widecheck{\beta'}$.  The edge $\e$ is in the immediate future
of a bead $\beta$, necessarily of a different colour than $\beta'$.

\begin{definition} \label{d:bite}
In the above situation,
we say that $\beta$ {\em bites $\beta'$ from the left} if $\beta$
lies to the left of $\beta'$ in $S$, and that $\beta$ {\em bites
$\beta'$ from the right} if $\beta$ lies to the right of $\beta'$ in $S$.
We say that the edges $\e$ and $\e'$ discussed above {\em exhibit} the biting.
\end{definition}

The above concepts
of biting and consumption replace the single, simpler, notion
of consumption from \cite[Section 5]{BG1}: there, since
the preferred future was a single edge,  if it was bitten 
it was consumed.   In \cite{BG1}, a frequently used concept was
for an edge to be `eventually consumed'.  In this paper, we need 
the following replacement:

\begin{definition} \label{d:eventbite}
Suppose that $\rho_1 \subset \mu_1(S)$ and $\rho_2 \subset
\mu_2(S)$ are beads in $\bot(S)$.  We say that $\rho_1$ is
{\em eventually bitten} by $\rho_2$ if there is a corridor $S'$
which contains a preferred future $\beta_1$ of $\rho_1$ and a bead $\beta_2$
in the future of $\rho_2$ so that $\beta_2$ bites $\beta_1$ in $S'$.
\end{definition}

With these definitions in hand, we have the following, which is an appropriate
replacement for \cite[Lemma 5.3]{BG1}.

\begin{lemma} [cf.  Lemma 5.3, \cite{BG1}]  \label{C0Lemma}
There exists a constant $C_0$ with the following property: if $\rho$
is a bead such that $f_{\#}(\rho)$ contains a left-fast displayed
edge $E$ and if $UV\rho$ is a (tight) path with $V\rho = V \odot \rho$ 
 and $|V | \geq C_0$ then for all $j \geq 1$ the preferred
future of $E$ is not bitten when $f^j(UV\rho)$ is tightened.
Moreover, $|f_{\#}^j(UV\rho)| \to \infty$ as $j \to \infty$.
\end{lemma}
\begin{proof}
We first prove the result in the special case that $V\rho$ is a \nf \ of
a left-fast edge $E_1$, where $\rho$ is the preferred future of $E_1$.
In other words, we will prove the existence of a constant $C_0'$ so that
if $|V| \ge C_0'$ then the statement of the lemma holds for the particular
path $U V\rho$.
(We will later reduce to this special case.)

Note that $V$ and $V\rho$ are monochromatic paths, and thus admit a
beaded decomposition.  Suppose first that $V$ does not contain any
beads of length greater than $J$.  In this case, the proof is entirely
parallel to that of \cite[Lemma 5.3]{BG1}, where we count using the number 
of non-vanishing beads rather than the number of edges.

In case $V$ contains long \gep s or long \pep s, we note that the cancellation
by $U$ on the left, and possibly by one of the edges in the \gep \ or \pep \ on
the right can only decrease the length of a \gep \ or \pep \ by at most $2B$
at each iteration.  Thus it is straightforward to include long \gep s and \pep s
into the above calculation.  We now turn to the general case.

Suppose that $V$ is an arbitrary path so that $V\rho = V \odot \rho$.
Then $V$ can shrink of its own accord (it needn't be beaded), and can be cancelled
by the future of $U$.  However,  there is certainly a constant $C_0$ so that
if $|V| \ge C_0$ then by the time this shrinking of $V$ combined with cancelling
by the future of $U$ can have reduced $V$ to the empty path, the future of the edge $E$ has
at least $C_0'$ edges to the left of its preferred future.  
We are then in the
special case that we dealt with first.
\end{proof}

The following two lemmas are proved in an entirely similar manner to \cite[Lemma 5.5]{BG1}. Recall that displayed edges are particular types of beads,
and the (abstract) preferred futures of beads were defined in
Definition \ref{PrefFuture}.  Recall from Remark \ref{rem:Pref} that the 
preferred future of a displayed edge whose label is not contained in a 
zero stratum is a single displayed edge.

\begin{lemma} \label{OldLemma} Let $\chi_1\sigma\chi_2$
be a tight path in $G$.
Suppose that $\chi_1$ and $\chi_2$ are monochromatic 
and that, for $i = 1,2$, 
the edge $E_i$ is displayed in $\chi_i$ and that $E_i$ is not
in a zero stratum.  Suppose that $\sigma$ is a concatenation of beaded paths.
Then the preferred futures of $E_1$ and $E_2$ cannot cancel
each other in any tightening of $f_{\#}(\chi_1)f_{\#}(\sigma)f_{\#}(\chi_2)$.

Suppose that $S$ is a corridor in a well-folded diagram, and that $\mu_1(S)$
and $\mu_2(S)$ are non-empty paths in $\bot(S)$, where $\mu_1$ and $\mu_2$ are
colours.  Suppose further that for $i = 1,2$ there is a  displayed edge $\e_i$
such that $\widecheck{\e_i}$ is not contained in a zero stratum.  Then the edges
in the semi-naive future of $S$ corresponding to the preferred futures of
$\e_1$ and $\e_2$ do not cancel each other when folding the semi-naive future
of $\bot(S)$ to form $\top(S)$.
\end{lemma}

\begin{lemma} \label{ParabolicOldEdges}
Let $S$ be a corridor and suppose that $\e_1$ and $\e_2$ are edges in
$\bot(S)$ whose labels lie in parabolic strata. In the naive future of
each $\e_i$ (that is, before even the beads have been tightened),
there is a unique edge $\e_i'$ with the same label as
$\e_i$. At no stage during the tightening of $\top(S)$
can $\e_1'$ cancel with $\e_2'$.
\end{lemma}

\begin{corollary} \label{ParaCancelledbyHigher}
A displayed edge in any coloured interval $\mu(S)$ which is 
 labelled by a parabolic edge $\widecheck E_i \in H_i$ can
only be consumed by an edge whose label is in 
$\overline{G \smallsetminus G_i}$.
\end{corollary}

\subsection{Abstract paths, futures and biting}
\label{AbstractFolding}

In many of the arguments in later sections, we wish to work with
concatenations of beaded paths in $G$
 rather than sides of  corridors in diagrams.  This is done as in 
Subsection \ref{stackDiags} by associating to such a path
$\rho = \rho_1\dots\rho_m$, with the $\rho_i$ beaded, 
the van Kampen diagram $\Delta(l,\rho)$ with boundary
label $t^{-l}\rho t^l \overline{f_{\#}^l(\rho)}$. But we modify
the usual definition of colour by defining the colours on the bottom
of the first (earliest) corridor not to be single edges but
rather to be intervals labelled $\rho_i$.
We then use the 
definitions of the previous subsection (biting, preferred future {\em{etc.}})
to define the associated
concepts {\em{for beads in $\rho$}}.

We emphasize, $\rho$ itself need not
beaded; only the $\rho_i$ are. We also emphasize that edges do not
have preferred futures, only beads do.

However, some beads are single, displayed edges, and when considered
as beads they do have a preferred future.

\section{Counting Fast Beads} \label{FastBeads}

This section is the analogue of \cite[Section 6]{BG1}; it is here
that the proof of Theorem \ref{MainThm} begins in earnest.

Let $\Delta$ be a minimal area van Kampen diagram, folded
according to the convention of Section \ref{Folding},
and fix a corridor $S_0$ in $\Delta$. As explained in 
Section \ref{Linear}, the core of our task is to bound
the number of beads in the decomposition of $\bot(S_0)$.
In order to do so, we must undertake a detailed study of
the preferred futures of these beads.

First we dispense 
with the case that $\widecheck{\beta}$ is a
vanishing \atom.  

\begin{lemma} \label{VanishLemma}
Suppose that $\mathcal{S}$ is the collection of beads in $S_0$
which are not vanishing \atom s.  If $\sum_{\beta \in \mathcal{S}}
|\beta| = D$ then $|S_0| \le B(D+1)$.
\end{lemma}
\begin{proof} This follows immediately from the Bounded Cancellation
Lemma.
\end{proof}
  
Narrowing our focus in the light of this lemma, we define:

\begin{definition} [Bead norm]\label{beadNorm} 
Given a concatenation
 $\rho=\rho_1\dots\rho_m$ 
 of  beaded paths, we
define the {\em{bead norm}} of $\rho$, denoted $\| \rho \|_\beta$, 
to be the number of non-vanishing beads in the concatenation.
(This is poor notation, since the norm depends
on the decomposition into the $\rho_i$ and not just the edge-path $\rho$. But
in the contexts we shall use it, specifically $\bot(S_0)$, it will
always be clear which decomposition we are considering.)
\end{definition}

\begin{remark}
All beads have length at least $1$.  Thus bead norm is dominated by
length.  In particular, estimates concerning
Bounded Singularities and Bounded Cancellation
remain true when distance is replaced by bead norm; cf.~Lemma 
\ref{EasyBits}.
\end{remark}

\begin{remark} \label{BeadsDontVanish} An important advantage
of bead norm over edge-length is that when one takes the
repeated images $f_{\#}^k(\chi)$ of a monochromatic path,
its length can decrease,
due to cancellation within beads, whereas bead norm cannot.
\end{remark}

In Definition \ref{BeadLength} we defined the bead length $[S]_\beta$ of a
corridor $S$ in a well-folded diagram.  It is convenient for our future
arguments to concentrate on non-vanishing \atom s, and hence on
bead norm rather than bead length.  However,
an immediate consequence of the Bounded Cancellation Lemma
is the following bi-Lipschitz estimate:

\begin{lemma} \label{BLengthNorm}
Suppose $S$ is a corridor in a well-folded corridor.  Then
\[	\| S \|_{\beta} \le [S]_\beta \le B \| S \|_{\beta}	.	\]
\end{lemma}

\subsection{The first decomposition of $S_0$} [cf. \cite{BG1}, 
Subsection 6.1]

Let $\beta$ be a bead in $S_0$ that is not a vanishing \atom. 
As we follow the preferred future of $\beta$ forwards in time, one of the
following events must occur:
\begin{enumerate}
\item[1.] The last preferred future of $\beta$ intersects the boundary
of $\Delta$ nontrivially.
\item[2.] The last preferred future of $\beta$ intersects a singularity
nontrivially.
\item[3.] The last preferred future of $\beta$ is bitten in a corridor $S$.
\end{enumerate}
We remark that, unlike in \cite{BG1}, these events are not
mutually exclusive; this is because a bead can consist of
more than one edge.

We shall bound the bead norm of $S_0$ by finding a bound on the 
number of non-vanishing beads  in each of the
three cases.

We divide Case (3) into two sub-cases:
\begin{enumerate}
\item[3a.] The preferred future of $\beta$ is bitten by a bead
that is not in the future of $S_0$.
\item[3b.] The preferred future of $\beta$ is bitten by a bead
that is in the future of $S_0$.
\end{enumerate}

\subsection{Bounding the easy bits} [cf. \cite{BG1}, Subsection 6.2]

Label the non-vanishing beads which fall into the above classes
$S_0(1)$, $S_0(2)$, $S_0(3a)$ and $S_0(3b)$, respectively.  We shall
see, just as in \cite{BG1}, that $S_0(3b)$ is by far the most 
troublesome of these sets.

The following lemma is proved in an entirely similar way to
\cite[Lemmas 6.1 and 6.2]{BG1}, using the Bounded Cancellation
Lemma and simple counting arguments.

\begin{lemma} \label{EasyBits}
\ 

\begin{enumerate}
\item $\| S_0(1) \|_\beta \le \n$.
\item $\| S_0(2) \|_\beta \le 2B\n$.
\item $\| S_0(3a) \|_\beta \le B \n$.
\end{enumerate}
\end{lemma}

We have thus reduced our task of bounding $\| S_0 \|_\beta$
to bounding the numbers of beads in $S_0(3b)$, i.e. to
understanding cancellation {\em within} the future of $S_0$.
The bound on the number of beads in $S_0(3b)$ is proved
in an analogous way to \cite{BG1}, and takes up a large
part of the remainder of this paper (through Section 
\ref{BonusSection}).

\subsection{The chromatic decomposition} \label{Chromatic}
[cf. \cite{BG1}, Subsection 6.3]

Fix a colour $\mu$ and consider the interval $\mu(S_0)$ in
$\bot(S_0)$ consisting of beads coloured $\mu$.

We shall subdivide $\mu(S_0)$ into five (disjoint but possibly
empty) subintervals according to
the fates of the preferred futures of the beads.

Let $l_\mu(S_0)$ be the rightmost bead $\beta$ in $\mu(S_0)$
such that $f_{\#}(\widecheck{\beta})$ contains a left-fast displayed edge
$\epsilon$
so that the preferred future of $\epsilon$ is eventually
bitten from the left from within the future of $S_0$.
Let $A_1(S_0,\mu)$ be the set of beads in $\mu(S_0)$ from the left
end up to and including $l_\mu(S_0)$.

Let $A_2(S_0,\mu)$ consist of those beads which are not in 
$A_1(S_0,\mu)$ but whose preferred futures are bitten from 
the left from within the future of $S_0$.

Let $A_3(S_0,\mu)$ denote those beads which do not lie in
$A_1(S_0,\mu)$ or $A_2(S_0,\mu)$ and which fall into
the set $S_0(1) \cup S_0(2) \cup S_0(3a)$.

All of the beads which are not in $A_1(S_0,\mu)$, $A_2(S_0,\mu)$
or $A_3(S_0,\mu)$ must have their preferred future bitten from
the right from within the future of $S_0$.

Analogous to the definition of $l_\mu(S_0)$, we define
a bead $r_\mu(S_0)$: the bead $r_\mu(S_0)$ is the
leftmost bead $\beta'$ so that $f_{\#}(\widecheck{\beta'})$ contains
a right-fast displayed edge whose preferred future is eventually bitten
from the right from within the future of $S_0$.

Let $A_4(S_0,\mu)$ denote those beads which are not in
$A_1(S_0,\mu)$, $A_2(S_0,\mu)$ or $A_3(S_0,\mu)$ and which
lie strictly to the left of $r_\mu(S_0)$. 

Finally, let $A_5(S_0,\mu)$ denote those edges not in 
$A_1(S_0,\mu)$, $A_2(S_0,\mu)$, $A_3(S_0,\mu)$ or
$A_4(S_0,\mu)$ which lie to the right of $r_\mu(S_0)$
(include $r_\mu(S_0)$ in $A_5(S_0,\mu)$ if it has not 
already been included in one of the earlier sets).

Now Lemma \ref{EasyBits} immediately implies

\begin{lemma}
\[      \sum_{\mu} \| A_3(S_0,\mu) \|_\beta \le (3B+1) \n       .       \]
\end{lemma}

We also have

\begin{lemma} \label{A1Short}
Let $C_0$ be the constant from Lemma \ref{C0Lemma} above.  Then
\begin{enumerate}
\item $\| A_1(S_0,\mu) \|_\beta, \| A_5(S_0,\mu) \| \le C_0$; and
\item $| A_1(S_0,\mu) \smallsetminus l_\mu(S_0)|,
| A_5(S_0,\mu) \smallsetminus r_{\mu}(S_0)|  \le C_0$.
\end{enumerate}
\end{lemma}

\begin{proof}
We prove the bounds only for $A_1(S_0,\mu)$, the proofs for
$A_5(S_0,\mu)$ being entirely similar.

The entire future of beads in 
$A_1(S_0,\mu)$ other than $l_\mu(S_0)$ must be eventually
consumed from the left from within the future of $S_0$;
cf.~\cite[Lemma 5.9]{BG1}.

If $\| A_1(S_0,\mu) \|_{\beta}$ or $| A_1(S_0,\mu) \smallsetminus 
l_\mu(S_0)|$ were greater than $C_0$ then we would conclude 
from Lemma \ref{C0Lemma} that
no left-fast bead in the immediate future of $l_\mu(S_0)$
could be bitten at any stage from the left from within the future of $S_0$,
contrary to the definition of $l_\mu(S_0)$.
\end{proof}

As we continue to follow the proof from \cite{BG1}, our next goal is to reduce
the task of bounding the bead norm of $S_0$ to that of bounding
the number of Nielsen beads contained in $A_2(S_0,\mu)$
and $A_4(S_0,\mu)$.  We focus exclusively on $A_4(S_0,\mu)$, the arguments
for $A_2(S_0,\mu)$ being entirely similar.

In outline, our argument proceeds  in analogy with the subsections
beginning with  \cite[Subsection 6.4]{BG1}, commencing with the decomposition
of $A_4(S_0,\mu)$ into subintervals $C_{(\mu,\mu')}$. But we quickly
encounter  a new phenomenon that requires an additional section
of argument ---
HNP cancellation, which does not arise in the case of positive
automorphisms.

\subsection{The decomposition of $A_4(S_0,\mu)$ into the $C_{(\mu,\mu')}$.}

All beads in $A_4(S_0,\mu)$ are eventually bitten from the right
from within the future of $S_0$.  For a colour $\mu' \neq \mu$, 
define a subset $C_{(\mu,\mu')}$ of $A_4(S_0,\mu)$ as follows:
given a bead $\sigma \in A_4(S_0,\mu)$, there is a bead $\sigma'$
in $S_0$ so that $\sigma$ is eventually bitten by $\sigma'$.
If $\sigma'$ is coloured $\mu'$ then $\sigma \in C_{(\mu,\mu')}$. 

The sets $C_{(\mu,\mu')}$ form intervals in $S_0$.

\section{HNP-Cancellation and Reapers} 
\label{TrappedHNP}

The results of the previous section
 reduce the task of bounding $\| S_0 \|_\beta$ to
that of establishing a bound on the sum
of the bead norms of the monochromatic
intervals $C_{(\mu,\mu')}$.  In 
\cite{BG1}, the corresponding intervals (also labelled $C_{(\mu,\mu')}$)
contained no exponential edges.  In the current context, however, 
there
may be exponential edges {\em{trapped}} in 
Nielsen paths, which may themselves be contained in beads of any type. 
This raises the concern that our attempts to control the length of the
$C_{(\mu,\mu')}$ in the manner of \cite{BG1} will be undermined by
the {\em{release}} of these trapped edges when the Nielsen path
is bitten, leading to  rapid growth in subsequent
nibbled futures of the Nielsen path. Our purpose in this section
is to develop  tools to control  this situation, specifically Lemmas
\ref{DealWithExp} and \ref{WholeNielsen}.

We must also deal with a second threat 
that arises from the phenomenon described in Example \ref{Ex:HNP};
we call this {\em Half Nielsen Path (HNP-) cancellation}.

Recall that a {\em{\pep}} is an edge path $\rho$ in $G$;
it is associated to a \gep and
 either $\rho$ or $\bar{\rho}$ is of the
form $E \bar{\tau}^k \bar{\nu} \gamma$
where $E$ is an edge with
$f_\#(E) = E\odot\tau^m$, where $\tau$ and $\nu$
are  Nielsen paths, and $\bar{\gamma}\nu$ is a terminal segment of $\tau$
(and $m, k > 0$).
These are the prototypes of the following types of paths.

\begin{definition}
Suppose that $E$ is a linear edge with $f_{\#}(E) = E \odot \tau^m$, where $\tau$ is a 
Nielsen path and $m > 0$.  Suppose further that $\nu$ is a Nielsen path and $\gamma$
an edge-path so that $\bar{\gamma}\nu$ is a terminal segment of $\tau$.  

A {\em \rpep} is a path $\rho$ so that either $\rho$ or $\bar{\rho}$ has the form
$E \bar{\tau}^k \bar{\nu} \gamma$ where $k > 0$.
\end{definition}

\begin{remark} Every
\pep\ is  a \rpep, but an arbitrary
 \rpep has no \gep associated to it.
\end{remark}

It is important to note that in the following definition the \rpep \ 
being discussed is {\em{not}} assumed
to be a bead in the decomposition of $\bot(S)$. (Beads along
$\bot(S)$ are monochromatic whereas we want to discuss HNP
cancellation, as in Definition \ref{HNPBite}, in the context
of adjacent colours interacting.)

\begin{definition}[HNP cancellation] \label{HNPDef} 
Let $S$ be a corridor in
a well-folded diagram, let
$\e$ and $\e'$ be edges in the 
naive (unfolded) future of $\bot(S)$ that cancel in the passage to $\top(S)$
and assume that  $\e$ is to the
left of $\e'$.

Suppose further that the past of $\e$ is $e$ with
label $\check e = E$ a linear edge 
and that $\e'$ is in the future of an edge $e_\gamma$
whose label is an edge $\gamma$.

We call the cancellation of $\e$ and $\e'$ {\em left HNP-cancellation}
and write $\e \lhnp \e'$
if the interval from $e$ to $e_\gamma$ in $\bot(S)$ (inclusive) is
labelled by a \rpep\
of the form $E \bar{\tau}^k \bar{\nu}\text{\em\o}\gamma$, where
$\tau$ is a Nielsen path so that $\tau = \xi \nu$, where $\xi$ and $\nu$
are Nielsen paths, and 
$\overline{\text{\em\o}\gamma}$ is a terminal sub edge-path of $\xi$.

{\em Right HNP-cancellation} is defined by reversing the roles of $\e$ and 
$\e'$
and insisting upon a \rpep\ in $\bot(S)$ of the form
$\bar{\gamma}\bar{\text{\em\o}} \nu \tau^k \bar{E}$.  It is denoted $\e \rhnp \e'$.

When we are unconcerned about the distinction between left and right, we refer
simply to {\em HNP-cancellation}.

We extend this definition to concatenations of beaded paths in $G$ by using the
obvious stack-of-corridors diagram as in Subsection \ref{stackDiags}.
\end{definition}

\begin{remark}
HNP-cancellation occurs at the `moment of
death' of the \rpep; see \cite[Section 6]{BG2} for an explanation of the
significance of this moment  and an analysis of it (in
the language of \pep s).
\end{remark}

\begin{lemma} \label{l:oEmpty}
Suppose that $E \overline{\tau}^k \overline{\nu} \text{\em\o}\gamma$ is a \rpep \
which exhibits an HNP-cancellation, as in Definition \ref{HNPDef}.
Then $\text{\em\o}$ is empty, so $\gamma$ is the first edge of $\overline{\xi}$.
\end{lemma}
\begin{proof}
The assumption that HNP-cancellation occurs means that we can
restrict our attention to cancellation when tightening
\[	f(E\overline{\tau}^k\overline{\nu}\text{\em\o}\gamma).	\]
This can be written as
\[	E \tau^m f(\overline{\tau}^k \overline{\nu}) f(\text{\em\o}\gamma).	\]
The path $\overline{\tau}^k\overline{\nu}\text{\em\o}\gamma$ admits a
hard splitting $\overline{\tau} \odot \cdots \odot  \overline{\tau}
\odot \overline{\nu} \odot \text{\em\o}\gamma$.
Therefore, under any choice of tightening, the $m$ copies of 
$\tau$ cancel with the $k$ copies of $f(\overline{\tau})$ (partially
tightened), then with $f(\overline{\nu})$; they then begin to
interact with $f(\text{\em\o}\gamma)$.  Just as in the proof of 
\cite[Proposition 6.9]{BG2}, under the assumptions of 
\cite[Lemma 5.1]{BG2}, there is only a single edge in $\text{\em\o}\gamma$
whose future can interact with $f(E)$ when tightening. 
\end{proof}

We now present the deferred example that explains the
need to consider HNP-cancellation. This will also lead
us to a further definition --- {\em{HNP biting}} --- that
encodes a genuinely troublesome situation where HNP
cancellation must be accounted\footnote{We usually
account for it by excluding it from our definitions.  When it
cannot be excluded, we often sidestep it, using the notions of `robust future'
and `robust past' given in Definitions \ref{Robust} and \ref{RobustPast}
below.}  for. Fortunately, many
other instances of HNP-cancellation are swept-up 
by our general cancellation and finiteness arguments, allowing
us to avoid a detailed analysis of the possible outcomes.

The problem at the heart of the following example did
not arise in \cite{BG1} because the natural
realisation of a positive automorphism does not map
any linear edge   across other linear edges.

\begin{example} \label{Ex:HNP}
Suppose that $u$ is a Nielsen path, and that $E_1$ and 
$E_2$ are edges so that $f(E_i) = E_i u^k$ for $i = 1,2$
and some integer $k > 0$.  For any integer $j$, the path
$\tau_j = E_1 u^j \overline{E_2}$ is an indivisible Nielsen path.

Suppose that $E_3$ is an edge so that $f(E_3) = E_3 \tau_j^l$,
for some integers $j$ and $l$ (with $l > 0$).  For ease of 
notation, we will assume that $l =1$.

Consider the path $\rho = E_3 \overline{\tau}_j^r E_2$,
for some $r > 0$.  Then $\rho$ is a \rpep.

In the iterated images $f_{\#}(\rho)$, the visible copy of $E_2$
has a unique future labelled $E_2$, which we will call the `preferred
future' of $E_2$ for the purposes of this example.
After $r+1$ iterations of $\rho$ under $f_{\#}$ (and any choice
of tightening at each stage), the future of $E_3$ cancels the preferred future
of the visible copy of $E_2$.  If we encode the evolution of $\rho$
in a stack diagram as in Subsection \ref{stackDiags} then the cancellation
of $E_2$ is HNP-cancellation.
\end{example}

In the following discussion, we assume that the reader is
familiar with  \cite{BG1}, in particular the vocabulary of
teams and reapers.

The phenomenon described in
the above example causes problems when 
the sub-path $\rho_1 = \overline{\tau}_j^r E_2$ of $\rho$ is
monochromatic and $E_2$ is displayed in $\rho_1$.  In
this situation, it shows that the most obvious
adaptation of \cite[Lemma 6.7]{BG1} would be false.
It is for this reason that we must exclude
HNP-biting in Definition \ref{DefEi}.

Similarly, because  Example \ref{Ex:HNP}
renders a naive version of the results of \cite[Section 8]{BG1} false,
HNP-biting must be excluded from
the Two Colour Lemma and the associated results
in Section \ref{PincerSection}.

A situation in which we cannot exclude HNP-biting by decree
arises in the analysis of teams and in particular
the definition of a {\em{reaper}} (subsection 8.3).
Suppose that  $\rho$ labels some interval in the bottom
of a corridor, with many copies of $\overline{u}$ to its immediate right.
In this case, the edge $\e_2$ labelled
$E_2$ will consume copies of $\overline{u}$
in the first $r$ units of time, but its future
will then be cancelled (assuming no
other cancellation occurs from either side, and that there are no
singularities, etc.).  Since $\e_2$ was acting as the reaper
of a team, we must find a continuing manifestation of it
at subsequent times, for otherwise we
will lose control over the length of teams
($r$ being arbitrary) and the structure of our main argument
will fail. This problem is solved by introducing the
{\em{robust future}} of $\e_2$
(Definition \ref{Robust}), which in this case is an edge labelled
$E_1$ that `replaces' the preferred future of $\e_2$ when it
is cancelled.

\begin{definition} \label{HNPBite}
Suppose that $\chi_1$ and $\chi_2$ are beaded paths in $G$
and $\chi_1\chi_2$ is tight.
Suppose that there is a bead $\rho_1 \subset \chi_1$
and a bead $\rho_2 \subset \chi_2$ so that
\begin{enumerate}
\item\label{rho1} either $\rho_1$ is a displayed edge $\gamma$ in $\chi_1$ which is linear or else
$\rho_1$ is a displayed \pep \ in $\chi_1$ of the form
$E \bar{\tau}^k \bar{\nu} \gamma$, where $\gamma$ is a linear edge;
\item when tightening $f_{\#}(\chi_1)f_{\#}(\chi_2)$ to form
$f_{\#}(\chi_1\chi_2)$, $\rho_1$ bites $\rho_2$ and the edge $\e'$
in the exhibiting pair $( \e',\e)$ (see Definition \ref{d:bite}) is
in the future of $\gamma$;
\item moreover\footnote{The \rpep\ implicit in the symbol $\lhnp$ is not
the \pep\ in \eqref{rho1}.}, $\e' \rhnp \e$.
\end{enumerate}
Under these circumstances we say that $\rho_2$ is {\em left-HNP-bitten} by $\rho_1$
and we write $\rho_1 \lbite \rho_2$.
There is an entirely analogous definition of {\em right-HNP-biting} $\rho_1 \rbite \rho_2$,
and when we are unconcerned about the direction we will refer
simply\footnote{We swap orientation in Definition \ref{d:evHNPb}
so as to emphasize  this point immediately.}
 to {\em HNP-biting}.

We make the analogous definition for HNP-biting within diagrams.
\end{definition}

\begin{definition}\label{d:evHNPb}
Suppose that $\chi_1$ and $\chi_2$ are beaded paths
and that $\rho_1$ is a bead in $\chi_1$.  We say that $\rho_1$
is {\em eventually HNP-bitten} by $\chi_2$ if $\rho_1$ is eventually
bitten by $\chi_2$ (Definition \ref{d:eventbite}) and this biting is HNP-biting.

We make the analogous definition within diagrams.
\end{definition}

\begin{definition}
Suppose that $E$ and $E'$ are edges in $G$.  We say that 
$E$ and $E'$ are {\em indistinguishable} if there is a Nielsen path
$\tau$ and an integer $s > 0$ so that $f(E) = E \tau^s$ and
$f(E') = E' \tau^s$.
\end{definition}

The edges $E_1$ and $E_2$ in Example \ref{Ex:HNP} are
indistinguishable.

\subsection{Parabolic HNP-cancellation and robust futures}

The following is a simple (but key) observation, and has an obvious
application to HNP-cancellation of edges of parabolic weight.

\begin{lemma} \label{HNPreplace}
Suppose that $\tau$, $\nu$, $\nu'$ and $\sigma$ are Nielsen paths, 
with $\sigma$ irreducible and
$\tau = \nu' \overline{\sigma} \nu$.  Suppose further that $\gamma$ is the initial edge
of $\sigma$, and that $f(\gamma) = \gamma \odot \xi^l$ for some Nielsen
path $\xi$.  Then $\sigma$ has the form $\gamma \xi^r \overline{\gamma'}$
where $r$ is some integer and $\gamma'$ is an edge so that
$\gamma$ and $\gamma'$ are indistinguishable.

Moreover, suppose that $E$ is an edge so that $f(E) = E \odot \tau^m$,
and let $\rho = E \overline{\tau}^i \overline{\nu} \gamma$ be a \rpep \ with
$0 \le i < m$. Then $f_{\#}(\rho)$ has the form
$E \odot \tau^{m-i-1} \nu' \gamma' \overline\xi^j$ where
$\gamma$ and $\gamma'$ are indistinguishable.
\end{lemma}
\begin{proof}
The first assertion is an immediate consequence of the structure of indivisible Nielsen
paths of parabolic weight, and the second is then obvious (a detailed analysis
of the Nielsen paths of parabolic weight is undertaken in \cite[Section 1]{BG2}).
\end{proof}

\begin{definition}\label{ParabPrefFut} In general, non-displayed edges 
$\e$ in diagrams do not 
have preferred futures. But if $\check\e$ has parabolic weight, there is a
unique edge of the same weight in $f_\#(\check\e)$, and it is natural
to define the (immediate) preferred future of $\e$ to be the corresponding
edge in the immediate future of $\e$. (If $\e$ happens to be displayed,
this agrees with our earlier definition.)
\end{definition}

In Section \ref{PincerSection}, when proving the Pincer Lemma,
we will have to exclude HNP-biting.  This will also be the case
in the applications of the Pincer Lemma in Sections \ref{TeamsSection}
and \ref{BonusSection}.  Thus, in following the future of  a linear
edge $\gamma$ when HNP-cancellation occurs, we would like to ignore the 
preferred future (which disappears), and rather follow the future of 
the interchangeable edge $\gamma'$ from Lemma \ref{HNPreplace}
above.  Thus we make the following

\begin{definition} [Robust Futures for Parabolic
Edges] \label{Robust}
Suppose that $\e$ is a (not necessarily displayed) edge in a colour
$\mu(S)$, and that $\widecheck{\e}$ is contained in a parabolic
stratum.  If the preferred future of $\e$ is cancelled from the left 
[resp.~right] by
HNP-cancellation in $\top(S)$, then
Lemma \ref{HNPreplace} provides an edge $\gamma'$
that is indistinguishable from $\widecheck\e$ and survives in the
tightened path $f_{\#}(E \bar{\tau}^k\bar{\nu}\text{\em\o} \gamma)$
[resp.~its reverse] considered in Definition \ref{HNPDef}.

We define the {\em robust
future} of an edge $\e \subseteq \bot(S)$ as follows.  If the preferred
future of $\e$ survives in $\top(S)$, then the robust future of $\e$
is just the preferred future of $\e$.  If the preferred future is
cancelled by HNP-cancellation, then the robust future of $\e$ is the
above edge labelled $\gamma'$, provided this survives in
$\top(S)$.  Otherwise there is no robust future.
\end{definition}

\begin{definition} [Robust Pasts for Linear Edges]\label{RobustPast}
Let $\e'$ be an edge of $\top(S)$ and suppose
that both it and its immediate past are labelled by linear edges. If
$\e'$ is not the robust future of any edge then the robust past of
$\e'$ is the past of $\e'$.  But if $\e'$ {\em{is}} the (immediate) robust
future of $\e$ then the robust past of $\e'$ is $\e$.
\end{definition}

Just as for preferred futures, the notions of robust future and robust 
past can be extended arbitrarily many steps forwards or backwards
in time by iterating the definition.

\subsection{A setting where we require cancellation lemmas}
\label{ss:setting}

Consider the following situation.  Let $\chi_1\sigma\chi_2$
be a tight path in $G$ with $\chi_1$ and $\chi_2$
monochromatic and  $\sigma$  a path with a preferred
decomposition into monochromatic paths (each of which comes
equipped with a beaded decomposition).  
We will analyse the possible
interaction between $\chi_1$ and $\chi_2$ in iterates of
$\chi_1\sigma\chi_2$ under $f$ (where the tightening follows
the convention of Subsection \ref{AbstractFolding}).

As ever, the following lemma remains valid with left/right orientation
reversed.

 \begin{lemma} \label{OneHNP}
Suppose that $\chi_1$, $\chi_2$ and $\sigma$ are as above, and 
suppose that each non-vanishing bead in $\chi_2$ is 
eventually bitten by a bead from $\chi_1$ in some
iterated image  $f_{\#}^k(\chi_1\sigma\chi_2)$
of $\chi_1\sigma\chi_2$.

Suppose further that $\rho$ is a bead in $\chi_2$ so that
$f_{\#}(\rho)$ has parabolic weight,
and that $\rho$ is eventually left-HNP-bitten 
by a bead from $\chi_1$ in the evolution of $\chi_1\sigma\chi_2$. Then $\rho$
is the rightmost non-vanishing bead in $\chi_2$.
\end{lemma}
\begin{proof}
Pass to the iterate $f_{\#}^{k-1}(\chi_1\sigma \chi_2)$ so that
the preferred future of $\rho$ lies in a \rpep\  $\pi$, which exhibits the
(eventual) HNP-biting of $\rho$ in the tightening to form
$f_{\#}^{k}(\chi_1\sigma \chi_2)$.  Let $\rho_1$ be the preferred
future of $\rho$ in $f_{\#}^{k-1}(\chi_1\sigma \chi_2)$.  Since
$f_{\#}(\rho)$ has parabolic weight, $\rho_1$ has parabolic weight,
and is either a displayed edge or a displayed \pep\ or \gep.
 We must prove
that no bead to the right of $\rho_1$ is eventually bitten by
the future of $\chi_1$.

By Definition \ref{HNPBite} and Lemma \ref{l:oEmpty} the
\rpep \ $\pi$ has the form $\gamma \bar{\tau}^k\bar{\nu} \e$, where
\begin{enumerate}
\item $\gamma$ is an edge so that $f(\gamma) = \gamma \odot \tau^m$;
\item $\gamma$ is either a displayed edge in the future of $\chi_1$ in 
$f^{k-1}_{\#}(\chi_1\sigma\chi_2)$ or else if the rightmost edge in a 
displayed \pep; and
\item $\e$ is contained in $\rho_1$.
\end{enumerate}
Let $\alpha$ be the displayed edge or \pep\ containing $\gamma$.

Let $\rho_1'$ be the terminal part of $\rho_1$ from $\e$ to its right end,
and let $\chi_2'$ be the terminal part of the future of $\chi_2$ in
$f_{\#}^{k-1}(\chi_1\sigma\chi_2)$, from $\e$ to its right end.

Since $\rho_1$ is displayed, we have $\chi_2' = 
\rho_1' \odot \beta$
for some path $\beta$.  

By Lemma \ref{HNPreplace}, when tightening to form
$f_{\#}^{k}(\chi_1\sigma \chi_2)$, the edge $\e$ is replaced
by an indistinguishable edge $\e'$ which comes from the
future of $\alpha$.  Suppose that $\delta$ is that part
of $f_{\#}(\alpha \rho_1')$ from $\e'$ to the right end.
Since $\alpha$ is a (linear) edge or a \pep,
the edge $\e'$ survives in all iterates of $\alpha$ (under
any choices of cancellation.  Similarly, since $\e$ and
$\e'$ are indistinguishable, $\e'$ survives in all iterates
of $\delta$ (under any choices of tightening).  This implies that
we have a hard splitting 
$f_{\#}(\alpha \chi_2') =  f_{\#}(\alpha \rho_1') \odot f_{\#}(\beta)$,
and the fact that $\alpha$ is displayed implies that
no bead in $\beta$ can be eventually bitten by the future
of $\chi_1$, as required.
\end{proof}

In applications of Lemma \ref{OneHNP} (and of Lemmas \ref{DealWithExp} and \ref{WholeNielsen} below), we usually take
$\chi_1 = \widecheck{\mu_1(S)}$ and $\chi_2 = \widecheck{\mu_2(S)}$,
where $\mu_1$ and $\mu_2$ are colours and $S$ is some corridor,
and we will choose $\sigma$ to be the label of that part of $\bot(S)$
which lies strictly between $\mu_1(S)$ and $\mu_2(S)$.\footnote{
However, it will also be convenient sometimes to take $\chi_1$ to
be a subinterval of $\widecheck{\mu_1(S)}$ consisting of an interval of beads.}
Since the folding conventions of Subsections \ref{Folding} and 
\ref{AbstractFolding} are compatible, and because of the hardness of our splittings,
the interaction between $\mu_1$
and $\mu_2$ in the future of $S$ can be analysed by studying the
interaction between the futures of $\chi_1$ and $\chi_2$ in iterated images
of $\chi_1\sigma \chi_2$ under $f$.

\subsection{Reapers} \label{ss:reaper}

In \cite{BG1} proving the existence of reapers was straightforward
(see \cite[Section 9]{BG1}).
In the current context, however, we have to work harder to
prove that a suitable incarnation of
a reaper exists, because of the phenomena discussed
in the preceding subsection. At the heart of our difficulties
is the fact that Nielsen atoms need not be single edges.

\begin{definition}
A {\em beaded Nielsen path} in a corridor $S$ is a subinterval
$\sigma \subset \bot(S)$ so that $\widecheck{\sigma}$ is a beaded
path all of whose beads are Nielsen paths.
\end{definition}

Note that in the above definition we do not assume that $\sigma$
is a single colour, or even that each bead in $\widecheck{\sigma}$
is contained in a single colour.  
Examples of beaded Nielsen paths include
that part of a \gep\ between the extremal edges,
and the sub-paths $\overline{\tau}^i$ of a \rpep\ 
$E\overline{\tau^k}\bar{\nu}\text{\em\o}\gamma$.

Although the beads in a beaded Nielsen path might not be displayed in
a path $\widecheck{\mu(S)}$, it is still possible to define the future
of a bead in a beaded Nielsen path, and the notions of preferred
future and biting still make sense.  We will use this observation in 
the sequel.

The following notion is parallel to that of \cite[Definition 10.1]{BG1},
which was pivotal in the bonus scheme (cf. Section
\ref{BonusSection} below). Here, it plays a more central role.

\begin{definition}[Swollen present and swollen future]
\label{Swollen}
Suppose $S$ is a corridor and that $I \subseteq \bot(S)$ is a
beaded Nielsen path in $S$. 
The {\em swollen present} of $I$ is the\footnote{uniqueness is immediate from
the observation that if a terminal sub-path $\sigma$ of a Nielsen path $\tau$ is itself Nielsen
then   $\sigma$ is a concatenation of beads in $\tau$.}
maximal subinterval $I' \subseteq \bot(S)$ such that (i) $I \subseteq I'$;
(ii) $I'$ is a beaded Nielsen path in $S$; and (iii) the beads of $I$ 
are beads of $I'$.

The {\em left swollen present} of $I$ is that part of the swollen present from the
left end up to the right end of $I$, whilst the {\em right-swollen present} goes
from the left end of $I$ to the right end of the swollen present.

If the actual future of $I$ is a beaded Nielsen path the
(immediate) {\em swollen future} $sw_1(I)$ of $I$ is the swollen present of the
(actual) future of $I$. With a similar 
qualification, the {\em swollen future} $sw_k(I)$ at $\time(S) + k$ is
defined to be $sw_1(sw_{k-1}(I))$.

With the same qualifications, the left and right swollen futures are defined in
the obvious ways.
\end{definition}

The first qualification in the above definition is required because it is possible
that the immediate future of a beaded Nielsen path is not a beaded Nielsen
path.  Thus we must be careful only to apply
this concept in cases where we know the swollen future to exist.

\begin{definition}[Reapers] \label{ReaperDef}
Suppose that $S$ is a corridor and $I \subset \bot(S)$ is
a beaded Nielsen path in $S$ with nonempty swollen future $sw_1(I)$.
Suppose that $\alpha$ is an edge in $\bot(S)$
immediately adjacent to $I$ on the left.  We say that $\alpha$ is
a {\em left-reaper for $I$} if (i) $\widecheck{\alpha}$ is a linear edge;
(ii) $\widecheck\alpha$ bites some of the future of $\widecheck{I}$ in 
$f_{\#}(\widecheck{\alpha I})$; and
(iii) the robust future of $\alpha$ is immediately adjacent to $sw_1(I)$ in $\top(S)$.

There is an entirely analogous definition of {\em right-reapers}.
As usual, when we are unconcerned about the direction we will refer
to {\em reapers}.
\end{definition}

\begin{definition}[Left-edible] \label{d:left-edible}
Let $S$ be a corridor in a well-folded diagram, and $I \subset \bot(S)$
a beaded Nielsen path.  
We say that $I$ is {\em left-edible}
if each bead in $I$ is eventually bitten by a bead coloured $\mu$ in the
future of $S$, where $\mu(S)$ lies to the left of $I$.

{\em Right-edible} paths are defined with a reversal of the left-right orientation.
\end{definition}

In the remainder of this section we work towards proving
Propositions \ref{OneSideStable} and \ref{Reaper}.

\begin{proposition} \label{OneSideStable}
Let $S$ be a corridor in a well-folded diagram and $I \subset \bot(S)$ a
left-edible path so that $|I| \ge B+J$.  Then the immediate future of
$I$ in $\top(S)$ is left-edible.
\end{proposition}

The following lemma is straightforward, and allows us to focus our attention
on the time when cancellation between colours begins.

\begin{lemma} \label{l:gotofuture}
Let $S$ be a corridor in a well-folded diagram and let $I \subset \bot(S)$
be a left-edible colour, all of whose beads are eventually bitten by
beads coloured $\mu$.  Let $S^I$ be the corridor in the future of $S$ so that
the first  biting of a bead in the left swollen future of $I$ by something
coloured $\mu$ occurs in $S^I$.  Then the left swollen future of $I$ in $\bot(S^I)$
is left-edible.
\end{lemma}

 In the following statement $B$ is the Bounded Cancellation Constant
from Proposition \ref{BCL} and $J$ is the constant from the Beaded
Decomposition Theorem \ref{BDT}.  The corridor $S^I$ is as in Lemma 
\ref{l:gotofuture} above, and $I^\lambda$ is the left swollen future of $I$ in
$S^I$.

\begin{proposition} \label{Reaper}
Suppose that $S$ is a corridor in a well-folded diagram and $I \subset \bot(S)$
is a left-edible path, all of whose beads are eventually bitten by beads coloured
$\mu$.  Suppose also that $|I| \ge B+J$.  Then
\begin{enumerate}
\item the immediate future of $I^\lambda$ in $\top(S^I)$ has an associated
left reaper $\alpha$, which is coloured $\mu$; and
\item for each bead in the immediate future of $I^\lambda$, when it is eventually
bitten the biting is by the robust future of $\alpha$.
\end{enumerate}
\end{proposition}

\subsection{Two Cancellation Lemmas}

The following lemma is useful in the proof of Lemma
\ref{C1Lemma} below. We record it now because a variation
on it (Lemma \ref{WholeNielsen}) is needed in the proof of 
Proposition \ref{Reaper}.

We revert to the setting described in Subsection \ref{ss:setting}.

\begin{lemma} \label{DealWithExp}
Assume that in the iterates of $\chi_1\sigma\chi_2$
(i.e.~forward-images  under $f_{\#}$)
each bead in $\chi_2$ is eventually bitten by a bead in $\chi_1$.
Suppose that $\chi_2$ has weight $i$, where $H_i$ is an
exponential stratum, and that all beads of weight $i$ in $\chi_2$
are Nielsen beads.
Let $\rho$ be a bead of weight $i$ in $\chi_2$.
\begin{enumerate}
\item  If $\rho$ is not bitten in
$f_{\#}(\chi_1\sigma\chi_2)$ but is eventually bitten
in the image $f_{\#}^k(\chi_1\sigma\chi_2)$  
then $\rho$ is entirely consumed in $f_{\#}^k(\chi_1\sigma\chi_2)$.
\item If $\rho$ is bitten but not entirely consumed in
$f_{\#}(\chi_1\sigma\chi_2)$ then
$\rho$ is the rightmost bead in $\chi_2$.
\end{enumerate}
\end{lemma}
\begin{proof}
There is at most one indivisible Nielsen path of weight $i$
and the lemma is vacuous unless there
is exactly one.

Let $\beta$ be a bead in $\chi_2$ of weight $i$, and suppose that
an edge $\eta$ in the future of $\chi_1$ is the edge which cancels the
rightmost edge in the preferred future of $\beta$ to exhibit the
biting of $\beta$ by $\chi_1$.
Since $\beta$ is an indivisible Nielsen path, it has edges of weight $i$
on both ends, as does its preferred future, 
and so $\eta$ has weight $i$.  Suppose that the 
past of $\eta$ in $\chi_1\sigma\chi_2$ has weight $i$.  Then by 
\cite[Theorem 8.1]{BG2} and Assumption \ref{FinalPower}, $\eta$
is either a displayed edge in the future of $\chi_1$, or else is
contained in a Nielsen bead.  
Suppose first that $\eta$ is contained in a Nielsen bead $\tau$.
Since $\eta$ is to cancel with an edge in $\beta$, the path $\tau$
must have weight $i$.  Hence $\tau = \bar\beta$, and $\beta$ is 
entirely consumed when it is bitten.

Suppose then that $\eta$ is displayed in the future of $\chi_1$.  By
Assumption \ref{FinalPower}.(\ref{DisplayExp}) we may assume
that the edge $\eta$ is contained in a displayed path of the form
$f(\eta)$.  Since $f(\eta)$ is $i$-legal, and $\beta$ is not, it is not
possible for the illegal turn in $\beta$ (of weight $i$) to be cancelled
by any iterates of $\eta$.  However, $|f(\eta)| > |\beta|$, by Assumption
\ref{FinalPower}\eqref{ExpLong}, so it is not possible for the displayed
copy of $f(\eta)$ to be cancelled by the future of $\beta$.  Therefore,
in this case $\beta$ must be the rightmost bead in $\chi_2$.

Furthermore, suppose that $\beta$ and $\eta$ are as above, and
the past of $\eta$ in $\chi_1\sigma\chi_2$ has weight $i$, 
and suppose moreover that 
$\beta$ is not bitten in $f_{\#}(\chi_1\sigma\chi_2)$.  Then $\beta$ is 
bitten by $\eta$ in some  $f_{\#}^k(\chi_1\sigma\chi_2)$, and $k \ge 2$.  
Thus we may assume that the immediate past of $\eta$ is also displayed
and is $\eta$.  By applying Lemma 
\ref{EndStab1} and noting that the rightmost edge of $\beta$ must be
$\bar\eta$, we see that the sub-path between the immediate
past of $\beta$ and the immediate 
past of $\eta$ has the form $\cdots \bar\eta \omega \eta \cdots$
for some path $\omega$.  The path $\omega$ must start and finish at he
same vertex, and in
order for the written copy of $\bar\eta$ to cancel with the written
copy of $\eta$ it must be that $f_{\#}(\omega)$ is a point.  However,
$\omega$ is not a point, because otherwise the past of $\beta$ and
the past of $\eta$ would already cancel.  This contradicts the fact
that $f$ is a homotopy equivalence.  The same argument
shows that if $\eta$ is contained in a Nielsen bead and $\beta$
is not bitten in $f_{\#}(\chi_1\sigma\chi_2)$
then $\beta$ cannot be bitten by $\eta$.

Therefore, if $\beta$ is bitten by an edge $\eta$ whose past in 
$\chi_1\sigma\chi_2$ has weight $i$ then $\beta$ is close to
the left end of $\chi_2$, and is either entirely consumed when
bitten or is the rightmost bead in $\chi_2$.

We may now assume that the bead $\rho$ is cancelled by an edge
$\eta$ whose past in $\chi_2$ has weight greater than $i$.  The above
arguments show that we may assume that the immediate past of
$\eta$ also has weight greater than $i$, and by Lemma \ref{EndStab1}
we may assume that this past is contained in a displayed edge,
a \gep, or a \pep.  It is easy to see that the immediate past of 
$\eta$ cannot have exponential weight and cannot be a \gep.  
Thus we may assume that the immediate past of $\eta$ is either
the edge on the left end of a \pep\ of the form $\gamma\nu\tau^k\overline{E}$,
 (and that the edge $\gamma$
is parabolic) or else is displayed and parabolic.

Lemma \ref{EndStab1} and the above arguments imply that this
immediate past of $\eta$ must be a linear edge, and the above
arguments now imply that if $\rho$ is bitten in a corridor it must be
entirely consumed.
 \end{proof}

 The following variant of Lemma \ref{DealWithExp} is the one
we need in the  proof of Proposition \ref{Reaper}.
 We continue to study $\chi_1\sigma\chi_2$ as in 
Subsection \ref{ss:setting}.

\begin{lemma} \label{WholeNielsen}
Suppose that $\chi_2$ is a beaded
Nielsen path and each of its beads is eventually bitten
by a bead in $\chi_1$ in some iterated image of $\chi_1
\sigma\chi_2$ under $f$.

Let $\rho$ be a bead in $\chi_2$ which is not bitten in $f_{\#}(\chi_1\sigma\chi_2)$.
If $\rho$ is bitten but not consumed in some iterated image of $\chi_1\sigma\chi_2$
then $\rho$ is the rightmost bead in $\chi_2$.
\end{lemma}
\begin{proof}
We follow the proof of Lemma \ref{DealWithExp} above, with
the added wrinkle that there may be parabolic weight Nielsen
paths to consider in $\chi_2$.  In this case there needn't be 
a unique Nielsen path of weight $i$.

Suppose that $\rho$ is as in the statement of the Lemma.  If
$\rho$ has exponential weight, then the arguments of the proof
of Lemma \ref{DealWithExp} give the required properties.  
If $\rho$ has parabolic weight,  Lemma 
\ref{ParabolicOldEdges} implies that when $\rho$ is bitten by
an edge $\eta$ in the future of $\chi_1$,  the immediate
past of $\eta$ has weight greater than that of $\rho$.  Also,
this immediate past must be parabolic. Arguing as in the proof of 
Lemma \ref{DealWithExp}, one sees  that either $\rho$ is entirely 
consumed when bitten, or else $\rho$ is the rightmost bead
in $\chi_2$.
\end{proof}

\begin{corollary} \label{cor:reaper}
Suppose that $I$ is a beaded Nielsen path in $\bot(S)$ for some corridor
$S$ of a well-folded diagram, and suppose that all beads of $I$ are eventually 
bitten from the left  by beads in a single colour $\mu$.  Then, with the possible
exception of $B$ beads 
on the left end and one bead on the right (the final one bitten),
whenever $\mu$ bites a Nielsen bead in the future of $I$, 
it consumes it entirely.
\end{corollary}

\bigskip

\noindent{\bf{Proof of the Proposition \ref{OneSideStable}}}
\smallskip

\begin{proof}
If the immediate future of $I$ in $\top(S)$ were not left-edible, then Corollary
\ref{cor:reaper} would ensure that no bead in $I$ which is not bitten in $S$ is ever
bitten by $\mu$.  However, the assumption on the length of $I$ (and the Bounded 
Cancellation Lemma) ensure that there {\em are} beads in $I$ not bitten in $S$.
The fact that $I$ {\em is} left-edible therefore ensures that the future of $I$ in
$\top(S)$ is also left-edible.
\end{proof}

\bigskip

\noindent{\bf{Proof of the Proposition \ref{Reaper}}}
\smallskip

\begin{proof}
Let $S'$ be the corridor containing
the immediate past of $I^\lambda$.
Lemma
\ref{WholeNielsen} implies that in $\top(S')$ there is an edge $\rho$ in
$\mu$ which cancels a whole Nielsen path in the future of $I$.

Since $|I| \ge B +J$, there is a bead in $I$ not bitten in $\top(S)$.  The proof of Lemma \ref{WholeNielsen} now implies that there is a reaper as in the statement of the proposition.
\end{proof}

\section{Non-fast and Unbounded Beads} \label{C1Section}

With the technical exertions of the previous section behind us,
we are now able to return to the  main argument,
picking up the flow of \cite{BG1} at Subsection 6.6. Thus our
next purpose is to reduce the task of bounding
the bead norm of the intervals $C_{(\mu,\mu')}$ to that of bounding
the lengths of certain long blocks of Nielsen atoms.  These blocks
are the analogue of the intervals $C_{(\mu,\mu')}(2)$ from
\cite{BG1}, and will be the building blocks of the {\em{teams}} 
introduced in Section
\ref{TeamsSection} (in analogy with \cite[Section 9]{BG1}).

\begin{definition} \label{Slowpep}
Suppose that $\rho = \gamma \nu \tau^k \overline{E_i}$ is a
\rpep \ (with $k \ge 0$).  We say that $\rho$ is {\em left-slow} if $\gamma$ is empty
 or a concatenation of left-slow beads.

There is an entirely analogous definition of {\em right-slow} \rpep s of the form
$\rho = E_i \overline{\tau}^k \overline{\nu}  \overline{\gamma}$.
\end{definition}

Often, we will just speak of {\em slow} \rpep s, since a single \rpep \ 
can only be left-slow or right-slow, but not both.

\begin{definition}
Suppose that the bead $\rho$ is such that $f_{\#}(\rho)$ is not a Nielsen
bead.  Then the function $n \mapsto |f_{\#}^n(\rho)|$ grows at least linearly.
In this case, we call $\rho$ an {\em unbounded bead}.
\end{definition}

\begin{definition} \label{tame}
A beaded path is called {\em right-tame} if all of its beads are
\gep s, slow \pep s, Nielsen paths and \atom s which do not
have a right-fast displayed edge in their immediate future.
\end{definition}

The next lemma follows immediately from the definition.
\begin{lemma}
$A_4(S_0,\mu)$ is a right-tame path.
\end{lemma}

\begin{lemma}
Suppose that $\alpha$ is a non-vanishing atom which is not right-fast.  Then either
all of the beads in $f_\#(\alpha)$ are Nielsen paths and \gep s, or else
the preferred future of $\alpha$ is parabolic.
\end{lemma}
\begin{proof}
The only modification to Lemma \ref{fofatom} is the 
exclusion of exponential edges in the second case, which is 
valid because such an edge would obviously contradict the fact 
that $\alpha$ is not right-fast.
\end{proof}

\begin{definition} \label{untrapped}
Suppose that $\sigma$ is a right-tame path.  The {\em untrapped
weight} of $\sigma$ is the largest $j$ so that $f_{\#}(\sigma)$ contains a
bead of weight $j$ which is not Nielsen.
\end{definition}

\begin{definition} \label{DefEi}
Suppose that, for some pair $(\mu,\mu') \in \mathcal Z$, the
untrapped weight of $C_{(\mu,\mu')}$ is $j$.
For each $1 \le i \le j$, define $\rho_i$
to be the leftmost bead in $C_{(\mu,\mu')}$ so that $f_{\#}(\rho_i)$
has an unbounded bead of weight at least $i$ that is not
HNP-bitten in the future of $S_0$.\footnote{Note that it is
possible that $\rho_i = \rho_{i+1}$ for some $i$.}

Let $\mathcal E_i$ denote those beads in $C_{(\mu,\mu')}$
from the right end up to and including $\rho_i$, and let
$\mathcal D_i = \mathcal E_i \smallsetminus \mathcal E_{i+1}$.
\end{definition}

The following is the analogue of \cite[Lemma 6.7]{BG1}
\begin{lemma} \label{C1Lemma}
For all $1 \le i \le \omega$ there is a constant $C_1(i)$ so that
for each of the paths $C_{(\mu,\mu')}$ and decomposition into intervals
$\mathcal D_i$ as above, we have
\[	\| \mathcal D_i \|_\beta \le C_1(i) .	\]
\end{lemma}
\begin{proof} As far as possible, we try to
follow the proof of
\cite[Lemma 6.7]{BG1}.  However, due to the phenomena
described in Section \ref{TrappedHNP}, the proof here is
somewhat more complicated.

We go forward to the time, $t$ say, which is one step before the moment
when $\mu'$ first starts to bite the preferred futures.  By virtue
of Remark \ref{BeadsDontVanish}, and the definition of 
$\mathcal D_i$, there are at least as many beads in the future
of $\mathcal D_i$ at time $t$ as there are in $S_0$.  Therefore,
it is sufficient to bound the number of beads in the future of 
$\mathcal D_i$ at time $t$; to ease the
notation, we write $\mathcal D_i$ for this future,
i.e.~pretend that $t=\time(S_0)$.

It is possible that there exist beads $\rho \in \mathcal D_i$ so that
$f_{\#}(\rho)$ has weight greater than $i$.  In such a case,
all of the beads in $f_{\#}(\rho)$ of weight greater than $i$ are Nielsen beads. 

Consider the highest weight $k$ for which there is a 
bead $\rho$ in $\mathcal D_i$ with $f_{\#}(\rho)$ of weight $k$,
and suppose that $k > i$.  Suppose first that $\rho$ has exponential
weight.  Then by Lemma \ref{DealWithExp} either $\mathcal D_i$ has
bead norm at most $B$ (and length at most $\ell=JB(B+1)$), or else
$\rho$ is entirely consumed when it is bitten.  In the first case
$\rho$ is the leftmost bead in $\mathcal D_i$, and also in $C_{(\mu,\mu')}$.
A similar argument applies when $\rho$ has parabolic weight.

Thus, excluding cases where $|\mathcal D_i|<\ell$, we may
 treat  the Nielsen beads of weight higher than $i$ as indivisible
units, which are entirely consumed when bitten.  We are therefore
in the situation of the proof of \cite[Lemma 6.7]{BG1}, where the
unbounded beads in $\mathcal C_i$ grow apart at a linear rate,
and so must be cancelled quickly.  Otherwise, the proof is entirely
parallel to the one from \cite{BG1}.
\end{proof}

We are trying to reduce the task of bounding the bead norm to
that of bounding the size of intervals consisting entirely
of Nielsen beads, which are each consumed by a reaper.  In order
to make this reduction, we still have some HNP-biting to deal with.
In order to deal with this, we need an analogue of \cite[Lemma 9.4]{BG1}.

Recall that $L$ is the maximal length of $f(E)$ where $E$ is an edge
in $G$.

\begin{proposition} [cf. Lemma 9.4, \cite{BG1}] \label{C4Lemma}
There is a constant $C_4$ depending only on $f$ which satisfies the
following properties.
If $I$ is an interval on $\top(S)$ labelled by a beaded path all of whose beads
are Nielsen \atom s, then the path labelling the past of $I$ in $\bot(S)$
is of the form $u\alpha v$ where $\alpha$ is a beaded path all of whose beads
are Nielsen \atom s and $|u|$ and $|v|$ are less than $C_4$. 

If the past of $I$ begins (respectively ends) with a point fixed by
$f$, then $u$ (respectively $v$) is empty.

In particular, $| I | \le |\alpha| + 2LC_4$.
\end{proposition}
\begin{proof} 
The interval $I \subset \top(S)$ is a beaded path, all of whose beads are  
Nielsen paths
of length at most $J$.  Therefore, along $I$ there are points where
$I$ admits a hard splitting and these points occur with a frequency of
at least one every $J$ edges.  Since these points are vertices, the
set of labels of points at which the splitting occurs is finite.
Consider the path from $\top(S)$ to $\bot(S)$ starting from one of
these vertices.  The label of this path is $w\bar{t_i}$ where
$w$ is a (possibly empty) path in $G$ of length at most $L$, and
$t_i$ is one of the edges from the mapping torus $M(f)$.  (We are about
to use a finiteness argument and it will be important that the
repetition we infer includes the labels of the points on
$\bot(S)$.  Thus it is important which of the $t$-edges this path
includes.)

Since the data we record --- the label of the vertex on
$\top(S)$, the path $w\bar{t_i}$ and the label of the end of this path
on $\bot(S)$) --- run over a finite set, there is a constant $C'$ such
that in the interval within $C'$ vertices of the left end of
$I$ there will be repetition of these data.  Since the vertices
occur at least every $J$ edges, this repetition occurs within $C'J$ of
the left end of $I$.

Once we have found this repetition, we have an interval $\lambda
\subset \bot(S)$, an interval $\eta \subset \top(S)$ and a path $w_0$
of length at most $L$ such that $f_{\#}(\lambda) = w_0 \eta
\bar{w_0}$. Therefore, the free homotopy class of $f_{\#}(\lambda)$ is
the same as that of $\eta = f_{\#}(\eta)$, since $\eta$ is a
beaded path all of whose beads are Nielsen paths. Since $f$ is a homotopy equivalence,
the free homotopy class of $\lambda$ must be the same as that of
$\eta$.

Suppose that $\eta = p_1 \dots p_m$ where each $p_i$ is an indivisible
Nielsen path.  Now, $\lambda$ is tight, so $\lambda = \sigma
p_ip_{i+1} \dots p_m p_1 \dots p_{i-1} \bar{\sigma}$, for some path
$\sigma$. Thus, if `$\sim$' denotes free homotopy,
\[ f(\lambda) \sim f_{\#}(\sigma) p_i \dots p_{i-1}
f_{\#}(\bar{\sigma}), \] 
which tightens to
\[ w_0 p_1 \dots p_m \bar{w_0}.  \]
By the Bounded Cancellation Lemma, tightening the path $f(\lambda)$ as
written above reduces the length of $f_{\#}(\sigma)$ by less than $B$,
and the result has length at most $2L + |\eta|$.  This implies that
$|f_{\#}(\sigma)| < L+B$.  Therefore, $\| \sigma
\|$ is bounded, and by a small increase we may also assume that
$i=1$.  By considering only one vertex out of every $B(L+B)$,
we can find such a path $\eta$ where there is some $p_j$ in the
middle of $\lambda$ such that the path from the copy of $p_j \subset \top(S)$ to the copy of
$p_j \subset \bot(S)$ is a single edge labelled $t$, for some $j$.

We have argued that, for some path $\eta$ of bounded length which lies
on the left end of $I$, the past of $\eta$ is of the form $u \eta u'$
where $|u|$ and $|u'|$ are bounded, and the paths from the splitting
points in $\eta \subset I$ to $\bot(S)$ consist of single edges
labelled $t$.

Consider the analogous situation on the right end of $I$.  We can
find a path $\eta' \subset I$ lies at the right end of $I$
such that the past of $\eta'$ is of the form $v' \eta' v$
where $|v|$ and $|v'|$ are bounded and the paths from the
vertices of $\eta' \subset I$ to $\bot(S)$ consist of single
edges labelled $t$.  

Consider the paths along $\bot(S)$ and $\top(S)$ from the left
end of $\eta$ to the right end of $\eta'$.  We have a path $\rho
\subset \bot(S)$ with fixed points of $f$ on either end which maps to
a Nielsen path $f_{\#}(\rho) \subset I \subset \top(S)$.  The same
argument as in the proof of \cite[Lemma 1.14]{BG2} then shows that
$\rho = f_{\#}(\rho)$.  Hence $\rho$ is a beaded path, all of whose
beads are Nielsen
paths, and the paths $u$ and $v$ on either side of $\rho$ are
of bounded length as required.  This proves the first assertion in the
statement of the lemma.

The second assertion follows similarly, and the final assertion
follows immediately from the first. 
\end{proof}

Consider a pair $(\mu,\mu') \in \mathcal Z$, and recall the
definition of the subintervals $\mathcal E_i$ from Definition \ref{DefEi}.

\begin{proposition} \label{GetC2}
There is a constant $C_5$, depending only on $f$ so that the
following holds.  For each $(\mu,\mu') \in \mathcal Z$, the
interval $C_{(\mu,\mu')} \smallsetminus \mathcal{E}_1$ 
in $A_4(S_0,\mu)$
has the form $u N v$ where $u$ and $v$ are such that
$\| u \|_{\beta}, \| v \|_\beta \le C_5$ and $N$ is a beaded path all
of whose beads are Nielsen beads.
\end{proposition}
\begin{proof}
By Lemma \ref{OneHNP}, for each adjacency of colours $(\mu,\mu')$
there can only be one bead in $\mu(S)$ which is eventually HNP-bitten
by $\mu'$.

The result now follows from Proposition \ref{C4Lemma} and the 
definition of $\mathcal E_1$.
\end{proof}

\begin{definition}
For $(\mu,\mu') \in \mathcal Z$, define $C_{(\mu,\mu')}(2):=
N$, the beaded Nielsen path from Proposition \ref{GetC2}.
\end{definition}

The sum of our arguments to this point has reduced the
task of bounding the sum of the bead norms of the 
intervals $\mu(S_0)$ in $S_0$ to that of bounding the
sum of the lengths of the intervals $C_{(\mu,\mu')}(2)$ for
pairs $(\mu,\mu') \in \mathcal Z$.

We summarise
the results from this section as follows.

\begin{proposition} \label{PointofC1}
There is a constant $C_1$, depending only on $f$, so that
\[      \| C_{(\mu,\mu')}\|_\beta \le \|C_{(\mu,\mu')}(2)\|_{\beta} + C_1.      \]
\end{proposition}

\begin{remark} \label{BeadtoLength}
Since the intervals $C_{(\mu,\mu')}(2)$ consist entirely of Nielsen
beads, we have the following obvious relationship between
length and bead norm:
\[ |C_{(\mu,\mu')}(2)| \le \|C_{(\mu,\mu')}(2) \|_\beta
\le J|C_{(\mu,\mu')}(2)|.       \]
Therefore, in order to finish the bound on bead norm, it is sufficient
to bound the total lengths of the intervals $C_{(\mu,\mu')}(2)$.
\end{remark}

It is important for the remainder of the paper that the path
$C_{(\mu,\mu')}(2)$ is a beaded path that consists entirely of
Nielsen atoms.  This is a stronger statement than just asserting
it is a Nielsen path, since we require a decomposition into beads
of uniformly bounded size, each of which is a Nielsen path.
This makes the path $C_{(\mu,\mu')}(2)$ very similar to the long
blocks of constant letters which played such a prominent role in
\cite{BG1}

At this point the reader may benefit from consulting 
\cite[Section 7]{BG1}, which outlines the strategy for the remainder 
of the proof of Theorem \ref{MainThm} (the strategy from the positive
case still holds here).  For the remainder of this paper, we will
mostly continue without reminding the reader of this strategy.

\section{The Pleasingly Rapid Disappearance of Colours} 
\label{PincerSection}

We are now at the point in our arguments where we need to formulate
and prove the Pincer Lemma, as in \cite[Section 8]{BG1}.
In \cite{BG1} the Pincer Lemma was proved by counting colours
which {\em essentially vanished}, which is to say they came to consist
entirely of constant letters.  For positive automorphisms, this is a
well-defined event and can only occur once for each colour.  For
general automorphisms, the analogues of constant letters are
indivisible Nielsen paths.  However, since Nielsen paths can contain
non-constant edges, indivisible Nielsen paths are not
indivisible in an absolute sense (the terminology refers to the fact that an indivisible
Nielsen path cannot be split into two Nielsen paths).  Thus, it is
possible that a colour can be labelled by a Nielsen path at some time
$t$ but not at some later time $t+k$. There are two ways to circumvent
this problem.  The first is to concentrate on the times when a colour
decreases in weight, whilst the second is to focus on the times when a
colour becomes Nielsen and seek compensation when a colour 
subsequently ceases to
be Nielsen. We mostly pursue the second idea but there are aspects of the first also.

The version of the Pincer Lemma which we need in this paper is
Theorem \ref{PincerLemma}.

The ideas in the proof of the Pincer Lemma here are very similar to
those in \cite{BG1} but the execution is somewhat different.

\begin{definition} \label{StablyNielsen}
Suppose that $I$ is a non-empty beaded Nielsen path
 and that $U$ and $V$ are
beaded 
paths.  We say that $I$ is {\em stably Nielsen} in the
path $UIV$ if the future\footnote{as defined in (\ref{stackDiags})}
of $I$ in $f_{\#}(UIV)$ is also a non-empty Nielsen beaded path.

Suppose that $\mu_1, \mu_2$ and $\mu_3$ are colours in a 
well-folded diagram and that the intervals
$\mu_1(S), \mu_2(S)$ and  $\mu_3(S)$ are non-empty and
adjacent in $\bot(S)$. If
$\widecheck{\mu_2(S)}$ is a non-empty Nielsen path, then we say 
that $\mu_2(S)$ is {\em stably Nielsen} if, in the above sense,
 $\widecheck{\mu_2(S)}$ is stably
Nielsen in $\widecheck{\mu_1(S)}
\widecheck{\mu_2(S)} \widecheck{\mu_3(S)})$.
\end{definition}

\begin{lemma}[Relative Buffer Lemma] \label{RelBuffer}
Let $i \in \{ 1 , \dots , \omega-1 \}$ and let $I \subset \bot(S)$ be an
edge-path labelled by edges in $G_i$. Suppose that the colours
$\mu_1(S)$ and $\mu_2(S)$ lie either side of $I$, adjacent to it.
Provided that the whole of $I$ does not die in $S$, no edge in the future 
of $\mu_1(S)$  with label in
$G \smallsetminus G_i$ will ever cancel
with an edge in the future of $\mu_2(S)$  with label in
$G \smallsetminus G_i$.
\end{lemma}
\begin{proof}
Given Lemmas \ref{EndStab1} and \ref{EndStab2},  the proof of \cite[Lemma 8.1]{BG1} applies modulo changes of terminology.
\end{proof}

We now need the following `two-sided' version of Proposition \ref{OneSideStable}.

\begin{lemma} \label{StableNielsenLemma}
Let $\mu_1$, $\mu_2$, $\mu_3$ and $S$ be as in Definition
\ref{StablyNielsen}, and suppose that $\mu_2(S)$ is stably Nielsen.  Then
for all corridors $S'$ in the future of $S$, if $\mu_1(S')$ and
$\mu_3(S')$ are nonempty then $\mu_2(S')$ is a (possibly empty)
Nielsen path.
\end{lemma}
\begin{proof} 
Whilst $\mu_1(S')$ and $\mu_3(S')$ are non-empty, any bead in $\mu_2$
which is bitten must be bitten by a bead coloured either $\mu_1$ or $\mu_3$.  Let
$I_1$ be the set of (Nielsen) beads in $\mu_2(S)$ which are eventually bitten
by a bead coloured $\mu_1$ (and are bitten whilst $\mu_1(S')$ and $\mu_3(S')$ are non-empty).
Define $I_2$ to be those beads in $\mu_2(S)$ which are bitten by a bead coloured $\mu_3$
(with the same proviso).

Suppose that $I_1$ and $I_2$ are non-empty.  They form intervals, and $I_1$
is to the left of $I_2$.

Proposition \ref{Reaper}, and the fact that $\mu_2(S)$ is stably Nielsen, implies
that unless $I_1$ is immediately consumed there is a left reaper coloured
$\mu_1$ associated to $I_1$, and similarly there is a right reaper coloured
$\mu_3$ associated to $I_2$.  The properties of reapers in Definition \ref{ReaperDef}
imply the result.  

In case one or both of $I_1$ and $I_2$ are empty (or immediately consumed), 
there is at most one reaper to consider, but the result follows in the same way.
\end{proof}

\begin{lemma}[Buffer Lemma] \label{BufferLemma}
Suppose, for some corridor $S$ in a well-folded diagram, that
$I \subset \bot(S)$ is a beaded Nielsen path and
that $\mu_1(S)$ and $\mu_2(S)$ lie either side of $I$,
 immediately adjacent to it.
 Suppose further that $\check{I}$ is stably Nielsen in $\check{\mu_1(S)}\check{I}
\check{\mu_2(S)}$.  Provided that the whole of $I$ does not die in $S$, no
bead in $\mu_1(S)$ can be eventually bitten by a bead coloured $\mu_2$ (and vice versa), 
unless it is (eventually) HNP-bitten.
\end{lemma}
\begin{proof}
Given Lemmas \ref{EndStab1}, \ref{EndStab2}
and \ref{StableNielsenLemma}, and the exclusion of HNP-biting, the proof 
of \cite[Lemma 8.1]{BG1} applies.
\end{proof}

The proof of the following lemma follows that
of \cite[Lemma 8.1]{BG1}.

\begin{lemma}[Weighted Buffer Lemma] \label{WtBuffer}
Suppose, for some corridor $S$ in a well-folded diagram, 
that $I \subset \bot(S)$ is a beaded path consisting
of Nielsen beads and beads of weight at most $i$, and 
that $\mu_1(S)$ and $\mu_2(S)$ lie either side of $I$,
 immediately adjacent to it.
Suppose further that the only  beads of $f_{\#}(\check{\mu_1(S)}\check{I}\check{\mu_2(S)})$ 
that are in the future of  $I$ and have weight greater than $i$ are
Nielsen beads.  

Then, provided that the whole of $I$ does not die in $S$, no bead
in $\mu_1(S)$ can be eventually bitten by a bead coloured $\mu_2$ (and vice versa), 
unless it is (eventually) HNP-bitten.
\end{lemma}

\subsection{The Two Colour Lemma}

Example \ref{Ex:HNP} can be used to construct examples where the
above two results are false if HNP-biting is not excluded. The same
is true of the results in this section. This accounts for the caution
that the reader will note in Sections
\ref{TeamsSection}, \ref{BonusSection} and 
\ref{LongGepsandPepsSection}, where we are careful to ensure
that the Pincer Lemma is applied only to  
pincers that  involve no HNP-biting.

\begin{definition} [Stable $f$-neutering]  
Suppose that $U$ and $V$ are beaded paths, that for some $k$ the futures of $V$ in
$f_{\#}^k(UV)$ and $f_{\#}^{k+1}(UV)$ are Nielsen, but that the future of $V$ in
$f_{\#}^{k-1}(UV)$ contains a non-Nielsen bead.  

Denote the futures of $U$ and $V$ in $f_{\#}^{k-1}(UV)$ by $U^{k-1}$ and $V^{k-1}$,
respectively.  Let $\beta$ be the  rightmost non-Nielsen bead in $f_{\#}(V^{k-1})$.
If the biting of $\beta$ in the tightening of $f_{\#}(U^{k-1})f_{\#}(V^{k-1})$ to form $f_{\#}^k(UV)$
is not HNP-biting then we say that $U$ {\em stably left $f$-neuters} $V$
in $k$ steps.

The definition of {\em stable right $f$-neutering} is identical with the roles of $U$ and $V$
reversed, and when we are unconcerned about the direction we will refer simply
to {\em stable $f$-neutering}.
\end{definition}
In the light of Proposition \ref{OneSideStable}, once stably $f$-neutered, the subsequent futures
of $V$ remain beaded Nielsen paths.

\begin{proposition} [Two Colour Lemma, cf. Proposition 8.4 \cite{BG1}]
\label{TwoColour} 
There exists a constant $T_0$, depending only on $f$, so that if $U$
and $V$ are beaded
paths and $U$ stably $f$-neuters $V$ then it does so in at most $T_0$ steps.
\end{proposition}

\begin{proof}
Denote the future of $U$ in $f_{\#}^i(UV)$ by $U^i$ and the future
of $V$ by $V^i$.

As in the proof of \cite[Proposition 8.4]{BG1}, we will decompose each
of the paths $V^i$ into an {\em{unbounded part}} and a 
{\em{bounded part}}.  The bounded
part will be an interval on the right end of $V^i$ whose immediate (abstract)
future is a beaded Nielsen path.  The unbounded interval lies on the left end
of $V^i$, and we will bound its length.  

This would
be a straightforward adaptation of the proof from \cite{BG1} if Proposition
\ref{PointofC1} provided a bound of the length of that part of $C_{(\mu,\mu')}$
not contained in $C_{(\mu,\mu')}(2)$.  However, the bound in Proposition
\ref{PointofC1} is just a bound on bead norm.
Thus, we need to deal with the possibility of long \gep s and \pep s.

The following enumerated claims will together yield an upper bound
on the length of the unbounded part of $V^i$, which in the course of 
the proof will be decomposed into $V^i_{\text{fast}}$ and $V^i_{\text{nc}}$ 

Three of the claims concern the existence of a constant $k_j$ that
depends only on $f$; we use the abbreviation $\exists k_j=k_j(f)$.

\smallskip

{\bf Claim 1:}  $\exists k_1=k_1(f)$ such that any \gep \ in $V^i$
has length less than $k_1$.
\smallskip

This follows in a straightforward way from the Buffer Lemma \ref{BufferLemma}
and the fact that the obvious preferred future of the rightmost edge in any \gep \ in
$V^i$ must eventually cancel with an edge from the future of $U^i$.

\smallskip

Next we consider long \pep s in $V^i$.  Suppose that $\rho$ is a \pep \ in
$V^i$.  Then the label on ${\rho}$ or ${\bar{\rho}}$ 
has the form $E\bar{\tau}^k\bar{\nu}\gamma$,
where $\tau$ is Nielsen path, $f(E) = E \odot \tau^m$ and $\bar{\gamma}\nu$ is
a terminal segment of $\tau$.  We consider a number of different cases.
First we dismiss a case that 
follows immediately from Lemma \ref{C0Lemma} and from the fact that
exponential edges are left-fast:

\smallskip

{\bf Claim 2:}  If $\check{{\rho}} = E\bar{\tau}^k\bar{\nu}\gamma$ and $\gamma$ is
an exponential edge then the right end of $\rho$ lies within $C_0$
of the left end of $V^i$.

\smallskip

Next we consider $V^i_{\text{fast}}$, which is defined 
to consist of those beads from the left end of $V^i$
up to and including the rightmost bead in $V^i$ whose immediate (abstract) future
contains a left-fast bead.  
\smallskip

{\bf Claim 3:} $\exists k_2=k_2(f)$ such that 
$|V^i_{\text{fast}}|\le k_2$.
\smallskip

This follows immediately from Lemma \ref{C0Lemma} unless the rightmost
bead in $V^i_{\text{fast}}$ is a \pep .  (Note that this rightmost bead is not a \gep ,
since a \gep \ does not have a left-fast bead in its immediate abstract future.)

Suppose, then, that the rightmost bead in $V^i_{\text{fast}}$ is a \pep , say $\rho$.
If $\check{{\rho}} = E\bar{\tau}^k\bar{\nu}\gamma$, then 
we are done by Claim 2. So
suppose that $\check{{\rho}} = \bar{\gamma} \nu \tau^k \bar{E}$.  Let
$\e$ be the edge in $\rho$ whose label is $\bar{E}$.  The preferred
future of $\e$ is to be cancelled by an edge in the future of $U^i$.  
By an obvious finiteness argument (as in the proof of \cite[Proposition 8.4]{BG1}), 
there is a constant $p$ so that the path $V^p$ contains
no left-fast beads.  This gives a bound on the amount of time before the future
of $\rho$ is bitten, and hence a bound on the amount that the future of $\rho$ can
shrink before then.  Suppose that $V^j$ is the first future of $V^i$ in which the future 
of $\rho$ has been bitten.  Because the preferred future of $\e$ is to be cancelled, 
Proposition \ref{Reaper} and the Buffer Lemma \ref{BufferLemma} imply that
the length of the future in $V^j$ of $\rho$ is bounded above by a 
constant depending only on $f$.

The required bound on $|V^i_{\text{fast}}|$
is now at hand:  Lemma \ref{C0Lemma} bounds the length of 
$V^i_{\text{fast}} \smallsetminus \rho$, and the combination of the bound on $j$
and the bound on the length of the future of $\rho$ in $V^j$ gives a bound on
the length of $\rho$.  This completes the proof of Claim 3.  We remark that the
above argument also gives a bound on the amount of time it takes for $V^1_{\text{fast}}$
to be entirely consumed.

\smallskip

We now define a set $V^i_{\text{nc}}$ as follows: Let $\rho_{\text{nc}}$ be the rightmost
bead in $V^i$ whose immediate abstract future is not Nielsen.  We define
$V^i_{\text{nc}}$ as follows:
\begin{enumerate}
\item if $\rho_{\text{nc}} \in V^i_{\text{fast}}$ then $V^i_{\text{nc}} = \emptyset$;
\item if $\rho_{\text{nc}}$ is not a \pep , then $V^i_{\text{nc}}$ consists of those
beads from (but not including) the rightmost bead in $V^i_{\text{fast}}$ up to
and including $\rho_{\text{nc}}$;
\item if $\rho_{\text{nc}}$ is a \pep \ with label of the form $\bar{\gamma} \nu \tau^k \bar{E}$
or $\rho_{\text{nc}}$ is a \pep \ with label of the form $E\bar{\tau}^k\bar{\nu}\gamma$
and $\gamma$ is not a Nielsen path, then $V^i_{\text{nc}}$ consists of those
beads in $V^i$ from (but not including) the rightmost bead in $V^i_{\text{fast}}$
up to and including $\rho_{\text{nc}}$;
\item \label{LooseEndRight} finally, if $\rho_{\text{nc}}$ is a \pep \ with label of the form
$E \bar{\tau}^k\bar{nu} \gamma$
and $\gamma$ is either empty or a Nielsen path, then $V^i_{\text{nc}}$ consists of
that interval from (but not including) the rightmost bead in $V^i_{\text{fast}}$ up
to and including the leftmost edge in $\rho_{\text{nc}}$ (the label of this
leftmost edge is $E$).
\end{enumerate}
Note that in Case \ref{LooseEndRight} the bead $\rho_{\text{nc}}$ is certainly
not contained in $V^i_{\text{fast}}$.

\smallskip

{\bf Claim 4:} $\exists k_3=k_3(f)$  such that  $|V^i_{\text{nc}}|\le k_3$.
\smallskip

The proof of Claim 3 above established an upper bound on the time before
all of $V^i_{\text{fast}}$ is entirely consumed, and hence also on the time
before the future of $V^i_{\text{nc}}$ begins to be consumed.
We now follow the proof of Lemma \ref{C1Lemma}, which
establishes an upper bound 
 on the time that can elapse before the final non-constant bead in $V^i$ is bitten.
We will be done if we can bound this time from below by a positive
constant times $|V^i_{\text{nc}}|$. 

In the current setting, we have non-constant beads in $V^i_{\text{nc}}$ 
that may not be growing apart like those in the proof of Lemma \ref{C1Lemma}.\footnote{This is because we are now measuring length rather than bead-norm.}
But  there {\em{is}} a lower
bound on the rate at which the surviving futures of these
beads can come together.
Hence the length of $V^i_{\text{nc}}$ provides a lower bound on 
the amount of time that must elapse before $V^j$ becomes  stably Nielsen,
since the future of $V^i_{\text{nc}}$ must be entirely consumed before
this time. (Note that in Case \ref{LooseEndRight},
the preferred future of the edge $\bar{E}$ in $\rho_{\text{nc}}$ must be 
eventually consumed by the future of $U^i$.)  This proves Claim 4.

\smallskip

The {\em unbounded part} of $V^i$ is the union of $V^i_{\text{fast}}$ and
$V^i_{\text{nc}}$, whilst the {\em bounded part} is the remainder of $V^i$.
The sum of the previous four claims bound the length of
the unbounded part of $V^i$ by a constant that depends only on $f$.

There is a similar bound on the number of edges in $U^i$ that have
an edge in their future that cancels with an edge in the future of $V^i$.
(Here we need the hypothesis that the path $V^k$ becoming stably Nielsen
does not arise from HNP-biting.)

At this stage, we can follow the proof of \cite[Proposition 8.4]{BG1} directly.
After an amount of time bounded by a constant
that depends only on $f$, either the future of $V$ 
becomes stably Nielsen or empty, or else there is a repetition 
of the following data: (i) the unbounded part of
$V^i$ plus the leftmost $B+J$ edges of the bounded part; (ii) a
terminal segment of $U^i$ containing all of the edges that can ever
interact with the future of $V$.  Once we have such a repetition, if
the future of $V$ has not become stably Nielsen or vanished then it never will,
contrary to hypothesis.
\end{proof}

We need a weighted version of neutering and the two-colour lemma.

\begin{definition} [$(f,i)$-neutering]
Fix $i \in \{ 1, \ldots , \omega \}$ and let $U$ and $V$ be beaded paths.
Suppose  that for some $k$ the future of $V$ in $f_{\#}^k(UV)$ has weight
less than $i$, but that the future of $V$ in $f_{\#}^{k-1}(UV)$ has weight at 
least $i$.

Denote the futures of $U$ and $V$ in $f_{\#}^{k-1}(UV)$ by $U_{k-1}$ and
$V_{k-1}$, respectively.  Let $\beta$ be the rightmost bead in $f_{\#}(V_{k-1})$
of weight at least $i$.  If the biting of $\beta$ in the tightening of
$f_{\#}(U_{k-1})f_{\#}(V_{k-1})$ to form $f_{\#}^k(UV)$ is not HNP-biting
then we say that $U$ {\em $(f,i)$-neuters} $V$ in at most $k$ steps.
\end{definition}

\begin{proposition}[Weighted Two Colour Lemma] \label{WeightTwo}
There exists a constant $T_0'$, depending only on $f$, so that for any
$i \in \{ 1 , \ldots , \omega \}$, if $U$ and $V$ are beaded
paths and $U$ $(f,i)$-neuters $V$ then it does so in at most $T_0'$ steps.
\end{proposition}
\begin{proof}
We decompose the futures of $U$ and $V$ in $f_{\#}^k(UV)$
as in Lemma \ref{TwoColour}.

The proof is similar to that of Lemma \ref{TwoColour}, 
except that when we appeal to the proof of
Proposition \ref{C1Lemma} we assume that we have a path
$\mathcal E_j$ with $j \ge i$.
Otherwise, the proof of Lemma \ref{TwoColour} above and that
of \cite[Proposition 8.4]{BG1} can now be followed {\em mutatis
mutandis}.
\end{proof}

By replacing $T_0$ by $T_0'$ if necessary, we may assume that $T_0 \ge
T_0'$. We henceforth make this assumption. 

\subsection{The disappearance of colours: Pincers and implosions}

\begin{definition} \label{PincerDef}
Consider a pair of non-constant edges $\e_1$ and $\e_2$ which cancel
in a corridor $S_t$ of $\Delta$, and suppose that, for $i = 1,2$,  the
immediate past of $\e_i$ lies in a bead of some $\mu_i(S_t)$ that is 
either a unbounded atom, a
\gep \ or a \pep .  Suppose further that the cancellation of
$\e_1$ and $\e_2$ is not HNP-cancellation, and that
$\mu_1 \neq \mu_2$.  Consider the paths $p_1,
p_2$ in $\F \subset \Delta$ tracing the histories of $\e_1$ and
$\e_2$.  Suppose that at time $\tau_0$ the paths $p_1$ and $p_2$ lie in a
common corridor $S_b$.  Under these circumstances, we define the {\em
  pincer} $\Pin = \Pin(p_1,p_2,\tau_0)$ to be the sub-diagram of
$\Delta$ enclosed by the chains of $2$-cells along $p_1$ and $p_2$,
and the chain of $2$-cells connecting them in $S_b$.

We define $S_{\Pin}$ to be the earliest corridor of the pincer in
which $\mu_1(S_{\Pin})$ and $\mu_2(S_{\Pin})$ are adjacent.  Define
$\wt{\chi}(\Pin)$ to be the set of colours $\mu \not\in \{ \mu_1 ,
\mu_2 \}$ such that there is a $2$-cell in $\Pin$ coloured
$\mu$. Finally, define
\[      \life(\Pin) = \time(S_{\Pin}) - \time(S_b) . \]
\end{definition}
See \cite[Section 8]{BG1} for illustrative pictures.

\begin{proposition}[Unnested Pincer Lemma, cf. Proposition 8.7 \cite{BG1}] \label{prePincer}

\ 

\noindent There exists a constant $T_1$, depending only on $f$, 
such that for any pincer $\Pin$
\[      \life(\Pin) \le T_1(1 + |\wt{\chi}(\Pin)|)      .       \]
\end{proposition}

In the proof of Proposition 8.7 (Regular Implosions) in
\cite{BG1}, the strategy was to identify a constant $T_1$ such
that over each period of time of length $T_1$ within a pincer, at
least one colour became constant.  There are a number of impediments
to implementing this strategy in the current situation.   The first is
that Nielsen paths can consist of edges which are not constant edges,
so if a colour {\em becomes Nielsen} then it may cease to be
Nielsen at some stage in the future.  In order to overcome this
impediment, we make the following

\begin{definition}
Suppose that for some colour $\mu$ and some corridor $S$, the path
$\check{\mu(S)}$ is stably Nielsen, and let $\nu_1$ and $\nu_2$ be the colours
immediately on either side of $\mu$ in $S$.  If there is some
corridor $S'$ in the future of $S$ in which $\check{\mu(S')}$ is not
Nielsen and $S'$ is the earliest such corridor, then we say that
$\mu$ is {\em resuscitated} in $S'$.  By Lemma
\ref{StableNielsenLemma}, at least one of $\nu_1$ and $\nu_2$ is not adjacent
to $\mu$ in $S'$,  so either $\nu_1(S')$ or $\nu_2(S')$ is empty.
If  $\nu_i(S')$ is empty, we say that $\nu_i$ {\em sacrifices
itself} for $\mu$.  
\end{definition}

\begin{remark}
A colour can sacrifice itself for at most one colour.

A colour may become stably Nielsen and be resuscitated a number of
times, but a different colour must sacrifice itself for each resuscitation.

The concept of `becoming stably Nielsen' is analogous to that of a colour
`essentially vanishing' in \cite[Section 8]{BG1}.  However, the concept
of `resuscitation' does not have an analogue in \cite{BG1}.
\end{remark}

Fix a pincer $\Pin$ and assume that $\life(\Pin) > 1$.  The strategy
to prove Proposition \ref{prePincer} is to identify a constant $T_1$
so that  during the life of $\Pin$, in each
$T_1/2$ steps of time there is a colour that becomes
stably Nielsen (perhaps vanishing) 
In order to obtain the bound in the statement of
Proposition \ref{prePincer},
we then count the colours which
become stably Nielsen or vanish, and the colours which sacrifice themselves 
for those that are
resuscitated.  A colour can therefore be counted twice
-- once for disappearing (or for the last time it becomes stably
Nielsen), and once as a sacrifice -- but no colour is counted more than
twice.  Thus Proposition \ref{prePincer} is an immediate consequence of
the following result whose proof will occupy the remainder of this
subsection.

\begin{proposition} \label{pSmall}
There is a constant $T_1$, depending only on $f$, so that for any
pincer $\Pin$ in a minimal area van Kampen diagram over $M(f)$, in any
interval of time of length $T_1 / 2$, at least one colour in
$\tilde{\chi}(\Pin)$ becomes stably Nielsen or vanishes.
\end{proposition}

\begin{definition} [$p$-implosive arrays] \label{implosive}
Let $p$ be a positive integer and $S$ a corridor.  A {\em
  $p$-implosive array} of colours in  $S$ is an ordered tuple $A(S) = [
  \nu_0(S), \ldots ,   \nu_r(S) ]$, with $r > 1$, such that
\begin{enumerate}
\item each pair of colours $\{ \nu_j, \nu_{j+1} \}$ is separated in
  $S$ only by a stably Nielsen (or empty) path;
\item in each of the corridors $S = S^1, S^2, \ldots , S^p$ in the
  future of $S$, no $\nu_j(S^i)$ is empty or a stably Nielsen path,  $j = 1
  , \ldots , r-1$;
\item in $S^p$, {\em either} an edge coloured $\nu_0$ from a
   unbounded \atom , a \gep \ or a
   \pep \ cancels with an edge coloured $\nu_r$ from a
   unbounded \atom , a \gep \ or a
  \pep \ (and hence the colours $\nu_j$ with $j = 1, \dots , r-1$ are
  consumed entirely), or {\em else} each of the colours $\nu_j$ ($j =
  1, \dots , r-1$) become stably Nielsen or vanish, while $\nu_0$ and $\nu_r$
  are not Nielsen in $f_{\#}(\widecheck{\nu_0(S^p)} \cdots
  \widecheck{\nu_r(S^p)})$ (although they may nevertheless become stably Nielsen
  or even disappear in $S^p$ because of colours external to the
  array).
\end{enumerate}
\end{definition}
Arrays satisfying the first of the conditions in (3) are said to be of
{\em Type I}, and those satisfying the second condition are said to be of
{\em Type II}.  (These types are not mutually exclusive).

The {\em residual block} of an array of Type II is the stably Nielsen
path which lies between $\nu_0(S^p)$ and $\nu_r(S^p)$ (if either
$\nu_0(S^p)$ begins or $\nu_r(S^p)$ ends with an interval of
Nielsen \atom s include these  in the residual block). Note that the
residual block may be empty.  The {\em
  enduring block} of the array is the set of stably Nielsen paths in
$\bot(S)$ that have a future in the residual block.

Note that there may exist some {\em unnamed colours} between $\nu_j(S)$ and
$\nu_{j+1}(S)$; if they exist, these form a stably Nielsen path.

\begin{remark}
Let $[ \nu_0(S), \dots , \nu_r(S)]$ be a $p$-implosive array.
\begin{enumerate}
\item Any $q$-implosive sub-array of $[ \nu_0(S), \dots , \nu_r(S)]$
  has $q = p$.
\item If an edge of $\nu_i$ cancels with an edge of $\nu_j$ and $j-i >
  1$, then this cancellation can only take place in $S^p$.  If the
  edges cancelling come from displayed unbounded \atom s, \gep s or
  \pep s, then the sub-array $[ \nu_i(S), \dots , \nu_j(S)]$ is
  $p$-implosive of Type I.
\item If $u, v$ and $w$ are beaded
  edge-paths such that $u$, $v$ and $f_{\#}(uwv)$ are Nielsen paths
  then $w$ is a Nielsen path.  It follows that the residual block of
  any array of Type II contains edges from at most two of the colours
  $\nu_j$, and if there are two colours then they 
are consecutive,  $\nu_j,\,\nu_{j+1}$.
\item Likewise, the enduring block of an implosive array of
  Type II is an interval involving at most two of the $\nu_j$ and if
  there are two such colours they must be consecutive.
\end{enumerate}
\end{remark}

\begin{lemma} \label{ExistsImpArray}
Let $\Pin$ be a pincer.  The ordered list of colours along each
corridor before $\time(S_{\Pin})$ in a pincer $\Pin$ must contain
a $p$-implosive array for some $p$.
\end{lemma}

\begin{proof}
The definition of $p$-implosive array is designed so that when a
colour becomes stably Nielsen (or disappears) in a pincer there is a
$p$-implosive array.  See the proof of \cite[Lemma 8.10]{BG1} for
more details.
\end{proof}

\begin{definition} \label{HNP-ImpArray}
Suppose that $A(S) = [\nu_0(S), \ldots , \nu_r(S)]$ is a $p$-implosive
array. We say that $A(S)$ is an {\em HNP-implosive array} if either
\begin{enumerate}
\item $A(S)$ is of Type I and in $S^p$ the cancellation between
  $\nu_0$ and $\nu_r$ is HNP-biting, or
\item $A(S)$ is of Type II and in $S^p$, for some $0 < i < r$,
  $\nu_0$ and $\nu_i$ are involved in HNP-biting or
  for some $0 < j < r$, $\nu_j$ and $\nu_r$ are involved in
  HNP-biting.
\end{enumerate}
\end{definition}

In order to follow the
arguments from \cite{BG1}, we need to sharpen
Lemma \ref{ExistsImpArray}:
HNP-cancellation can beget $p$-implosive arrays with $p$
arbitrarily large, and therefore we must
argue for the frequent occurrence of 
$p$-implosive arrays that are not HNP-implosive. A first step in this
direction is given by the following

\begin{lemma} \label{NoHNP}
Let $\Pin$ be a pincer, and let $\mu_1$ and $\mu_2$ be the
colours associated to the bounding-paths $p_1$ and $p_2$ 
of $\Pin$.  Then there is no HNP-biting between beads in 
$\mu_1$ and $\mu_2$ within $\Pin$.
\end{lemma}

\begin{proof} Follows from Lemmas \ref{OneHNP} and \ref{DealWithExp}.
\end{proof}

When we are unconcerned about $p$ in a $p$-implosive array, we refer
merely to an {\em implosive array}.  The first restriction to note
concerning implosive arrays is this:

\begin{lemma} \label{NotManyColours}
If $[\nu_0(S), \dots , \nu_r(S)]$ is implosive of Type I, then $r \le
B$.  If it is implosive of Type II, then $r \le 2B$.
\end{lemma}

\begin{proof}
In Type I arrays, the interval $\nu_1(S^p)\cdots \nu_{r-1}(S^p)
\subset \bot(S^p)$ is to die in $S^p$, so the bound is an immediate
consequence of the Bounded Cancellation Lemma.  For Type II arrays,
one applies the same argument to the intervals joining $\nu_0(S^p)$
and $\nu_r(S^p)$ to the residual block.
\end{proof}

\begin{proof}[Proof of Proposition \ref{pSmall}]
We give a suitable formulation of `short' so that in
any corridor $S$ within $\Pin$, $S$ contains a short
$p$-implosive array.  Proposition \ref{pSmall} then follows from an
obvious finiteness argument.

Let $A(S) = [\nu_0(S) , \ldots , \nu_r(S)]$ be the implosive array
guaranteed to exist by Lemma \ref{ExistsImpArray}, and suppose that $p
\ge 2T_0$ (if not then a colour becomes stably Nielsen or vanishes within $2T_0$
of $\time (S)$). 

We can decompose each of the colours $\nu_j(S)$ in 
analogy with \cite{BG1}, using the decomposition in
Section \ref{Chromatic} above.

We fix a constant $\Lambda_1$ so that if $\| A(S) \| >
\Lambda_1$ then one of the following must occur in $S^{T_0}$:
\begin{enumerate}
\item there is a block of displayed Nielsen \atom s in some
  $\nu_j(S^{T_0})$ of length at least $J+4B$,
\item there is a displayed \gep \ in some $\nu_j(S^{T_0})$ of length
  at least $J + 4B + 2$,
\item there is a displayed \pep \ in some $\nu_j(S^{T_0})$ of length
  at least $J + 4B + L + 1$, or
\item there is an interval of unnamed colours in $A(S)$ (which form a
  stably Nielsen block) of length at least $J + 4B$ between
  $\nu_0(S^{T_0})$ and $\nu_r(S^{T_0})$.
\end{enumerate}

In the remainder of the proof, we shall use the term {\em{block}}
to refer generically to the identified interval in whichever
of the above cases we find ourselves. Increasing $\Lambda_1$
if necessary, we may assume 
that the past of the block in $S$ satisfies the
relevant condition from (1) -- (4)  with the bound increased
by $2BT_0$.

For such a block $I$ in $S^{T_0}$, consider the first edge on either side
of this block which is not contained in a Nielsen path.  These edges
may be on one end of a  \gep \ or a \pep \ (including the \gep
\ or \pep\ from condition (2) or (3)), or may be contained
in unbounded \atom s.  Call these edges $\e_1$ and $\e_2$.

The Buffer Lemma \ref{BufferLemma} implies that either (i) one of
$\e_1$ and $\e_2$ must be `stabbed in the back' -- we do not exclude
the possibility that this stabbing happens by HNP-biting, or
(ii) there is HNP-cancellation across the above block.

We first dispose of case (ii).  Suppose, for ease of notation,
that the edge $\e_1$ HNP-bites the edge $\e_2$ across the above
block $I$.  Let $\e_1$ have weight $k$.  Then all edges in $I$ and
$\e_2$ must have weight less than $k$.  Let $\e_2'$ be the first edge
to the right of $I$ that has weight at least $k$.  Then the Relative
Buffer Lemma \ref{RelBuffer} implies that either $\e_1$ or $\e_2'$
must be stabbed\footnote{Note that if there is no such edge $\e_2'$
in $\Pin$ then $\e_1$ must be stabbed in the back, by Lemmas
  \ref{RelBuffer} and \ref{NoHNP}.}
in the back (again, this could be by
HNP-biting). 

We have argued that some edge must be stabbed in the back.
Suppose that this stabbing is of an edge $\e$ in $S^{T_0}$ and that $\e$ has
weight $k_1$.  Consider first the possibility that $\e$ is stabbed in
the back via HNP-biting.  Then this occurs by an edge $\e'$ of
weight at least $k_1 + 1$.  Now, either this stabbing in the back
occurs within $T_0$ of $S^{T_0}$, or by the Weighted Two Colour Lemma
(\ref{WeightTwo}) there is another block as in (1) -- (4) above.  This
block has higher weight than the previous block, and as above leads to
another stabbing in the back.  If this stabbing is HNP-biting,
pass to a yet higher weight stabbing, and so on.

Eventually (after less than $\omega$ iterations of this argument), we 
get an edge $\e$ stabbed in the back with the stabbing not 
HNP-biting.  Suppose that $\e$ has weight $k_2$.  Suppose for
ease of notation that $\e$ is to the left of the long block, and
suppose that $\e$ is coloured $\nu_i$.  Because of the block of Nielsen
\atom s to the non-stabbing side of $\e$, the Two Colour Lemma
(Proposition \ref{TwoColour}) implies that if the edge $\e'$ which
stabs $\e$ in the back is coloured by $\nu_j$ then $i - j > 1$;
we then write $\nu_j \searrow \nu_i$.

Passing to an innermost pair $\nu_{l_1} \searrow \nu_{l_2}$ between
$\nu_i$ and $\nu_j$ we can see that there are no blocks in $S^{T_0}$
satisfying any of (1) -- (4) above, for otherwise there would be a
further stabbing, leading to a related pair of colours between our
innermost pair, contradicting the innermost nature of this pair.

Once there are no such blocks, we have a bound on the length of the
$p$-implosive array implicit in the
relation $\nu_{l_1} \searrow \nu_{l_2}$. An obvious
finiteness argument now finishes the proof.
\end{proof}

We have already seen how Proposition \ref{pSmall} implies 
Proposition \ref{prePincer}.  Just as in \cite[Section 8]{BG1},
we must now deal with the 
possibility of `nested pincers'.  

\subsection{Super-buffers}

\begin{definition} We consider sequences of $5$-tuples of tight edge-paths in $G$.
\[ U_k := \Big( u_{k,1}, u_{k,2}, u_{k,3}, u_{k,4}, u_{k,5} \Big),
\ \ k= 1,2,...  \]
with $|u_{k,1}|$ and $|u_{k,2}|$ at most $C_0 + C_1 + 2B(B+1) + 1$, 
while $|u_{k,2}|$ and
$|u_{k,4}|$ are at most $C_0+C_1+J$ and $|u_{k,3}| \le 4B(B+1)+1$.\footnote{The purpose of these constants is just as in \cite[Definition 8.19]{BG1}, with appropriate changes due to Lemmas \ref{A1Short}
and \ref{VanishLemma} and Proposition \ref{PointofC1}.}   
We fix an integer $T_1'$ sufficiently large to ensure that 
for any sequence of length $T_1'$ there will be a repetition,
i.e.~some
$t_1 < t_2 \le T_1'$ with
\[  \Big( u_{t_1,1}, u_{t_1,2}, u_{t_1,3}, u_{t_1,4}, u_{t_1,5} \Big) =
\Big( u_{t_2,1}, u_{t_2,2}, u_{t_2,3}, u_{t_2,4}, u_{t_2,5} \Big) . \]
We also choose
$T_1' \ge T_1$.
\end{definition}

With appropriate changes of terminology and the results
of the previous subsection in hand, the
proof of 
\cite[Proposition 8.21]{BG1} yields:

\begin{lemma} \label{SuperBuffer}
Let $V = V_1 V_2 V_3$ be a tight concatenation of three beaded paths in $G$.
If the future of $V_2$ is not stably Nielsen in $f_{\#}^{T_1'}(V)$ then
the future of $V_2$ is not stably Nielsen in $f_{\#}^k(V)$ for any
$k \ge 0$.
\end{lemma}

\subsection{Nesting and the Pincer Lemma}

Let $\lambda_0 = J+2B(T_0+1)+1$, which is the obvious analogue
of the constant of the same name in \cite[Section 8]{BG1}.  As in
\cite[Remark 9.5]{BG1}, it is convenient to assume that
$LC_4 < \lambda_0$, and we increase $\lambda_0$ 
to make this so. (This makes certain statements in Section \ref{TeamsSection}
easier, but has no serious affect.)

\begin{definition} \label{NestDef}
Consider one pincer $\Pi_1$ contained in another $\Pi_0$.
Suppose that in the corridor $S \subseteq \Pi_0$ at the top of
$\Pi_1$ (where its boundary paths $p_1(\Pi_1)$ and $p_2(\Pi_1)$
come together) the future in $\top(S)$ of at least one of the edges
containing $p_1(\Pi_1) \cap \top(S)$ or $p_2(\Pi_1)\cap \top(S)$
is not contained in any stably Nielsen path and this future
\footnote{We allow this future to be empty, in which case ``contained
in" means that the immediate past of the long stably Nielsen path
is not separated from $\Pi_1$ by any edge that has a future in
$\top(S)$.} lies in a beaded path consisting of Nielsen beads and
beads of weight strictly less than the weight of the edges containing
$p_1(\Pi_1) \cap \top(S)$ and $p_2(\Pi_1)\cap \top(S)$, and that this
beaded path has at least $\lambda_0$ non-vanishing beads.
Then we say that $\Pi_1$ is {\em nested} in
$\Pi_0$. 
\end{definition}

\begin{remark} \label{remark:nesteddef}
Besides the obvious translations, the above differs from
\cite[Definition 8.22]{BG1} in that the path at the top
of the pincer may now consist of Nielsen beads and lower 
weight beads, whereas in \cite{BG1} it consisted entirely
of constant letters.  This more general setting does not make
any of the proofs in this section harder (because of the
Weighted Two Colour Lemma), but is needed because
of the more complicated definition of the `cascade of pincers'
below (Definition \ref{CascadeDef}).
\end{remark}

\begin{definition}
For a pincer $\Pi_0$, let $\{ \Pi_i \}_{i \in I}$ be the set of all
pincers nested in $\Pi_0$.  Then define
\[      \chi(\Pi_0) = \tilde{\chi}(\Pi_0) \smallsetminus \bigcup_{i \in I}
\tilde{\chi}(\Pi_i)     .       \]
\end{definition}

The corridor $S_t$ was defined in Definition \ref{PincerDef}.

\begin{lemma} \cite[Lemma 8.25]{BG1}
If the pincer $\Pi_1$ is nested in $\Pi_0$ then
$\time(S_t(\Pi_1)) < \time(S_{\Pi_0})$.
\end{lemma}
\begin{proof}
The existence of the beaded path at the top of the pincer
$\Pi_1$ makes this an immediate consequence of the 
Weighted Buffer Lemma \ref{WtBuffer}.
\end{proof}

Define $T_1 = T_1' +2T_0$.  The following theorem is the main result
of this section, and is the strict analogue of \cite[Theorem 8.26]{BG1}.
The proof in the current context follows the proof from \cite{BG1} {\em mutatis mutandis}.

\begin{theorem} [Pincer Lemma] \label{PincerLemma}
For any pincer $\Pi$
\[ \life(\Pi) \le T_1 (1 + |\chi(\Pi)|) .       \]
\end{theorem}

\section{Teams} \label{TeamsSection}

By virtue of Lemma \ref{PointofC1}, Remark \ref{BeadtoLength} 
and the results of Section \ref{FastBeads},
we have reduced the task of bounding the bead norm 
of $S_0$ to that of 
bounding the lengths of certain blocks $C_{(\mu,\mu')}(2)$ which
consist of Nielsen beads coloured $\mu$ all of which are to
be eventually bitten by beads coloured $\mu'$ in the future of $S_0$.  
By Proposition \ref{Reaper}, if
such a block has length at least $B+J$, then there is an 
associated reaper, which consumes Nielsen beads in
$C_{(\mu,\mu')}(2)$ at a constant rate (and entirely consumes
any bead it bites, up to the final bead).  Note that to each
pair $(\mu,\mu')$ there is at most one associated reaper.

This puts us in the situation where we can develop the technology
of {\em teams} as in \cite[Section 9]{BG1}. However, there are
a number of key differences to \cite{BG1}: we
already had to work hard in 
Section \ref{TrappedHNP} to establish  the existence of 
a {\em reaper} for
$C_{(\mu,\mu')}(2)$, and now we have to work harder 
to identify the times
$\hat{t}_1(\mu,\mu')$ and $t_1(\mathcal{T})$ attached to a pair
$(\mu,\mu') \in \mathcal Z$ and a team $\mathcal T$, using
the {\em robust} past of the reaper instead of the actual past;
this is required  in order that the 
 Pincer Lemma apply to teams of genesis
(G3).  It is worth remarking that once we have identified the pincer
$\Pi_{\mathcal{T}}$ associated to a team $\mathcal T$ of genesis (G3),
we revert to an analysis of actual pasts (as in the definition of pincer).

Note that the colour of the edges in the robust future of an edge may 
not always be the same, contrary to the actual future. In fact, whenever
the robust past is not the actual past, the colour changes.
 This explains a 
slight difference between Definition \ref{preTeam} below
and \cite[Definition 9.1]{BG1}.

Consider an interval $C_{(\mu,\mu')}(2)$ so that
$|C_{(\mu,\mu')}(2)| > B+J$, and let $\epsilon^\mu$
be the reaper associated to $C_{(\mu,\mu')}(2)$ in 
Proposition \ref{Reaper} above.  Let $t_0$ be the 
time at which $\epsilon^\mu$ first bites a Nielsen bead
in $C_{(\mu,\mu')}(2)$, and let $\beta_\mu$ be the
rightmost bead in the future of $C_{(\mu,\mu')}(2)$ at 
this time.  Note that $\beta_\mu$ is a Nielsen bead.
Let $\epsilon_\mu$ be the rightmost edge in $\beta_\mu$.

\begin{remark}
Since $|C_{(\mu,\mu')}(2)| > B+J$, and each bead of 
$C_{(\mu,\mu')}(2)$ is to be bitten by $\mu'$, the colour
 of $\epsilon^\mu$ is $\mu'$.
\end{remark}

\begin{lemma}
Suppose that the immediate past of $\epsilon_\mu$
exists (i.e.~that $\epsilon_\mu$ does not lie on
$\partial \Delta$).  Then the immediate past of $\epsilon_\mu$ 
lies in some bead $\sigma$, and $\sigma$ contains
the immediate past of each edge in $\beta_\mu$.
\end{lemma}

The above lemma, applied at each stage in the past,
implies that we can follow the past of the edge $\epsilon_\mu$
and deduce consequences about the past of all edges
in $\beta_\mu$.

We now define a time $\hat{t}_1(\mu,\mu')$ as follows:
We go back to the last point in time when
(i) the past
of $\epsilon_\mu$ and the robust past of $\epsilon^\mu$
lay in a common corridor; and
(ii) $\epsilon_\mu$ is contained in a beaded Nielsen path
whose swollen present is immediately adjacent to the
robust past of $\epsilon^\mu$.

{\em{We denote this corridor $S_{\uparrow}$.}}

\begin{definition} \label{preTeam}
The robust past of $\epsilon^\mu$ at time $\hat{t}_1(\mu,\mu')$
is called the {\em reaper}, and is denoted $\hat\rho(\mu,\mu')$.
The interval $\hat{\mathfrak T}(\mu,\mu')$ is the maximal
beaded Nielsen path in $\bot(S_{\uparrow})$ all of whose
beads are eventually bitten by $\hat{\rho}(\mu,\mu')$.
The {\em pre-team} $\hat{\mathcal T}(\mu,\mu')$ is defined to be the 
set of pairs $(\mu_1,\mu_2) \in \mathcal Z$ so that (i)
the robust past of $\epsilon^\mu$ is coloured $\mu_2$ at some time
between $\hat{t}_1(\mu,\mu')$ and $t_0$; and (ii) 
$\hat{\mathfrak T}(\mu,\mu')$ contains some edges coloured
$\mu_1$.
The number of beads in $\hat{\mathfrak T}(\mu,\mu')$ is denoted
$\| \hat{\mathcal T} \|$.
\end{definition}

As in \cite[Section 9]{BG1}, we will define {\em teams} to be pre-teams
satisfying a certain maximality condition (see Definition \ref{TeamDef}
below).

\begin{remark}
Just as in \cite[Remark 9.2]{BG1}, if $\hat{t}_1(\mu,\mu') < \time(S_0)$ then near
the right-hand end of $\hat{\mathfrak{T}}(\mu,\mu')$ one may have an
interval of colours $\nu$ for which $\nu(S_0)$ is empty.
\end{remark}

\begin{lemma} [cf. Lemma 9.3, \cite{BG1}] \label{TeamBound}
If $\hat{t}_1(\mu,\mu') \ge \time(S_0)$ then
\[      \sum_{(\mu_1,\mu_2) \in \hat{\mathcal{T}}(\mu,\mu')} |C_{(\mu,\mu')}(2)|
\le \| \hat{\mathcal{T}}(\mu,\mu') \| + B(B+1). \]
\end{lemma}
\begin{proof}
The extra $B(B+1)$ is to account for the beads consumed
before the reaper comes into play. Otherwise
the proof is just as in \cite{BG1}.
\end{proof}

\subsection{The Genesis of pre-teams} [cf. Subsection 9.2, \cite{BG1}]

We consider the various events that may occur at $\hat{t}_1(\mu,\mu')$ which prevent
us pushing the pre-team back one step in time.  Recall that $S_{\uparrow}$ is the corridor
at time $\hat{t}_1(\mu,\mu')$ which contains $\hat{\mathcal T}(\mu,\mu')$.  Suppose that
$\mu_2$ is the colour of $\hat{\rho}(\mu,\mu')$.

There are four types of events:
\begin{enumerate}
\item[(G1)]  The immediate past of $C_{(\mu,\mu_2)}(S_\uparrow)$ is separated from the 
robust past
of $\hat{\rho}(\mu,\mu')$ by an intrusion of $\partial \Delta$.
\item[(G2)]  We are not in Case (G1), but the immediate past of $C_{(\mu,\mu_2)}(S_\uparrow)$
is separated from the robust past of $\hat{\rho}(\mu,\mu')$ because of a singularity.
\item[(G3)]  The immediate past of $C_{(\mu,\mu_2)}(S_\uparrow)$ is still in the same corridor
as the robust past of $\hat{\rho}(\mu,\mu')$, but the swollen present of the immediate
past of $C_{(\mu,\mu_2)}(S_\uparrow)$ is not immediately adjacent to the robust past
of $\hat{\rho}(\mu,\mu')$.
\item[(G4)]  We are not in any of the above cases, but the immediate past of the rightmost
edge in $C_{(\mu,\mu_2)}(S_\uparrow)$ is not contained in a beaded Nielsen path.
\end{enumerate}

We now make the definition of a team.

\begin{definition}[cf. Definition 9.6, \cite{BG1}]\label{TeamDef}
All pre-teams $\hat{\mathcal T}(\mu,\mu')$ with $\hat{t}_1(\mu,\mu') \ge \time(S_0)$ are
defined to be teams, but the qualification criteria for pre-teams with $\hat{t}_1(\mu,\mu') < \time(S_0)$ are more selective.

If the genesis of $\hat{\mathcal T}(\mu,\mu')$ is of type (G1) or (G2), then the rightmost
component of the pre-team may form a pre-team at times before $\hat{t}_1(\mu,\mu')$.  In
particular, it may happen that $(\mu_1,\mu_2) \in \hat{\mathcal T}(\mu,\mu')$ but
$\hat{t}_1(\mu,\mu') > \hat{t}_1(\mu_1,\mu_2)$ and hence $(\mu,\mu') \not\in
\hat{\mathcal T}(\mu_1,\mu_2)$.  To avoid double-counting in our estimates on $\| \mathcal T \|$
we disqualify the (intuitively smaller) pre-team $\hat{\mathcal T}(\mu_1,\mu_2)$ in these
settings.

If the genesis of $\hat{\mathcal T}(\mu,\mu')$  is of type (G4), then again it may happen that
what remains to the right of $\hat{\mathcal T}(\mu,\mu')$ at some time before $\hat{t}_1(\mu,\mu')$
is a pre-team.  In this case, we disqualify the (intuitively larger) pre-team $\hat{\mathcal T}(\mu,\mu')$.

The pre-teams that remain after these disqualifications are now defined to be {\em teams}.
\end{definition}

A typical team will be denoted $\mathcal T$ and all hats will be dropped from the notation
for their associated objects (just as in \cite[Section 9]{BG1}).

A team is said to be {\em short} if $\| \mathcal T \| \le \lambda_0$ or
$\sum_{(\mu_1,\mu_2) \in \mathcal T} |C_{(\mu_1,\mu_2)}(2)| \le \lambda_0$.
Let $\Sigma$ denote the set of short teams.

\begin{lemma} [cf. Lemma 9.7, \cite{BG1}]
Teams of genesis (G4) are short.
\end{lemma}

We wish our ultimate definition of a team to be such that every pair $(\mu,\mu')$
with $C_{(\mu,\mu')}(2)$ non-empty is assigned to a team.  The above definition
fails to achieve this because of two phenomena:  first, a pre-team $\mathcal T(\mu,\mu')$
with genesis of type (G4) may have been disqualified, leaving $(\mu,\mu')$
teamless; second, in our initial discussion of pre-teams we excluded pairs $(\mu,\mu')$ with $| C_{(\mu,\mu')}(2)| \le B + J$.  The following
definitions remove these difficulties.

\begin{definition} [Virtual team members] \label{Virtual}
If a pre-team $\pTmm$ of type (G4) is disqualified under the terms of Definition \ref{TeamDef} 
and the smaller team necessitating disqualification is $\pT(\mu_1,\mu_2)$,  
then we define $(\mu,\mu')\vin\pT(\mu_1,\mu_2)$ and $\pTmm\subset_v\pT(\mu_1,\mu_2)$. 
We extend the relation $\subset_v$ to be transitive and extend $\vin$ correspondingly. 
If $(\mu,\mu')\vin\T$ then $(\mu_2,\mu')$ is said to be a {\em virtual member} of  
the team $\T$. 
\end{definition}

\begin{definition} 
If  $(\mu,\mu')$ is such that  $1\le |C_{(\mu,\mu')}(2)|\leq B+J$ and 
$(\mu,\mu')$ is neither a member nor a virtual member of any previously 
defined team, then we define $\T_{(\mu,\mu')}:=\{(\mu,\mu')\}$ to be a
(short) team with $\|\T_{(\mu,\mu')}\|=|C_{(\mu,\mu')}(2)|$.  
\end{definition}
 
\begin{lemma} [cf. Lemma 9.10, \cite{BG1}]\label{allIn} 
Every   $(\mu,\mu')\in\vecZ$ with $C_{(\mu,\mu')}(2)$ non-empty is a member 
or a virtual member of exactly one team, and there are less than $2\n$ teams. 
\end{lemma}  
\begin{proof} 
The first assertion is an immediate consequence of the preceding 
three definitions, and the second  follows 
from the fact that $|\vecZ| < 2\n$. 
\end{proof} 

\subsection{Pincers associated to teams of genesis (G3)} [cf. Subsection 9.3,\cite{BG1}]

In this subsection we describe a pincer $\Pi_\T$ canonically associated to each team
of genesis (G3), as in \cite[Subsection 9.3]{BG1}.  The only real difference between
the definitions here and those in \cite{BG1} is the use of robust past and 
beaded Nielsen paths.  Sadly, this variation leads to complications
in the cascade of pincers; see
Definition \ref{CascadeDef} and Remark \ref{remark:nesteddef}.

\begin{definition} [cf.~Definition 9.11, \cite{BG1}]
The {\em narrow past} of a team $\T$ at time $t$
consists of those beaded Nielsen paths whose beads are displayed in their
colour and whose future is contained in $\mathfrak T$.
The narrow past may have several components at each time, the set of which
are ordered left to right according to the ordering in $\mathfrak T$ of their
futures.  We call these components {\em sections}.
\end{definition}

{\em For the remainder of this subsection we consider only long teams of
genesis (G3).}

The following lemma follows from the definition of teams
of genesis (G3) in a straightforward manner.

\begin{lemma} \label{BeadsForPincer}
Let $\T$ be a team of genesis (G3).
There exist beads $y(\T)$ and $y_1(\T)$ of different colours,
both lying strictly between
the immediate past of the swollen present of $\T$
and the robust past of $\hat\rho(\mu,\mu')$, so that
$y(\T)$ is bitten by $y_1(\T)$ and this is not HNP-biting.
\end{lemma}

\begin{definition} [The Pincer $\tilde{\Pi}_\T$] \label{prePincerDef}
Choose a leftmost pair of beads $y(\T), y_1(\T)$ satisfying Lemma
\ref{BeadsForPincer}, and let $x(\T)$ be the leftmost
edge in $y(\T)$.  Let $x_1(\T)$ be the edge in $y_1(\T)$ which
is the past of the edge which cancels with the leftmost edge 
in the immediate future of $x(\T)$.

Define $\tilde{p}_l(\T)$ to be the path in the family forest $\mathcal F$
that traces the history of $x(\T)$ to $\partial\Delta$, and let $\tilde{p}_r(\T)$
be the path that traces the history of $x_1(\T)$.

Define $\tilde{t}_2(\T)$ to be the earliest time at which the paths
$\tilde{p}_l(\T)$ and $\tilde{p}_r(\T)$ lie in the same corridor.
\end{definition}

\begin{remark}
Since the pair $y(\T), y_1(\T)$ in Definition \ref{prePincerDef}
are the leftmost pair satisfying Lemma 
\ref{BeadsForPincer}, any non-vanishing beads which lie
between $\mathfrak T$ and this pair are involved
in HNP-biting and are of lower weight than $y_1(\T)$, by the
Weighted Buffer Lemma \ref{WtBuffer}.
\end{remark}

\begin{lemma} \label{Lemma:GetPincer}
The segments of the paths $\tilde{p}_l(\T)$ and $\tilde{p}_r(\T)$, 
together with the  path joining them along the bottom of the corridor 
at time  $\tilde{t}_2(\T)$ form a pincer.
\end{lemma}
\begin{proof}
Note that when choosing the beads $y(\T)$
and $y_1(\T)$ we excluded HNP-cancellation.  That the
paths in the statement of the lemma form a pincer then follows
immediately from the definition of pincers.
\end{proof}
We denote the pincer described in Lemma 
\ref{Lemma:GetPincer} above by $\tilde{\Pi}_\T$.

\subsection{The cascade of pincers}
The Pincer Lemma argues for the regular disappearance of colours
within a pincer during those times when more than two colours continue
to survive along its corridors.  However, when there are only two colours,
the situation is more complicated.

Recall that the constant $T_0$ is as in Proposition \ref{TwoColour}, subject
to the requirement that $T_0 \ge T_0'$ as in the assumption immediately
after Proposition \ref{WeightTwo}.  The pincer $S_{\Pin}$ associated to a pincer
$\Pin$ is defined in Definition \ref{PincerDef}.

\begin{lemma} \label{CascadeLemma}
One of the following must occur: 
\begin{enumerate}
\item \label{GetPincer} $\time(S_{\tilde{\Pi}_\T}) > t_1(\T) - T_0$;
\item the path $\tilde{p}_l(\T)$ and the entire narrow past
of $\T$ are not in the same corridor at time $t_1(\T) - T_0$; or
\item \label{GetCascade} at time $t_1(\T) - T_0$ the path 
$\tilde{p}_l(\T)$ and the
narrow past of $\T$ are separated by a path which does not split
as a beaded path whose beads are either
 Nielsen paths or of weight less than $\tilde{p}_l(\T)$.
\end{enumerate}
\end{lemma}
\begin{proof}
If not, the Weighted Two Colour Lemma (Lemma \ref{WeightTwo}) would give a contradiction, since there {\em is} to be interaction between
the beads $y(\T)$ and $y_1(\T)$ at time $t_1(\T)$, and this interaction
is not HNP-biting.
\end{proof}

We now consider each of the three cases in turn, seeking
a definition of times $t_2(\T)$ and $t_3(\T)$ and (possibly)
a pincer $\Pi_\T$.  The following definition is entirely
analogous to \cite[Definition 9.13]{BG1}, with the appropriate
translations.

\begin{definition} [cf. Definition 9.13, \cite{BG1}] \label{CascadeDef} 
\ 
\begin{enumerate}
\item \label{Case1} Suppose some section of the narrow past of $\T$ 
is not in the same corridor as $\tilde{p}_l(\T)$ at time $t_1(\T)
- T_0$:  In this case\footnote{this includes the possibility
that $\tilde{p}_l(\T)$ does not exist at time $t_1(\T) - T_0$}
we define $t_2(\T) = t_3(\T)$ to be the earliest time
at which the entire narrow past of $\T$ lies in the same
corridor as $\tilde{p}_l(\T)$ and has length at least $\lambda_0$.
\item \label{Case2} Suppose that Case (\ref{Case1}) does not
occur and $\time(S_{\tilde{\Pi}_\T}) > t_1(\T) - T_0$.  We define
$\Pi_\T = \tilde{\Pi}_\T$ and $t_3(\T) = \time(S_{\Pi_\T})$.  If the
narrow past of $\T$ at time $t_1(\T) - T_0$ has length less than
$\lambda_0$, we define $t_2(\T) = t_3(\T)$, and otherwise
$t_2(\T) = \tilde{t}_2(\T)$.
\item \label{Case3} Suppose that neither Case (\ref{Case1}) or
Case (\ref{Case2}) occurs:  In this case,
Lemma \ref{CascadeLemma}(\ref{GetCascade}) pertains.  We pass to the
latest time at which there is a path between $\tilde{p}_l(\T)$ and
the narrow past of $\T$ which has an edge of at least the same
weight as $\tilde{p}_l(\T)$ at this time and is not contained
in a Nielsen path.  Choose a pair of beads $y'(\T)$, $y_1'(\T)$
as in Lemma \ref{BeadsForPincer}, as
well as edges $x'(\T)$, $x_1'(\T)$.  Let $\tilde{p}'_l(\T)$ 
be the path 
tracing the history of $x'(\T)$.  Let
$\tilde{p}'_r(\T)$ trace the history of the edge $x_1'(\T)$ that
cancels $x'(\T)$.  Let $\tilde{t}_2'(\T)$ be the earliest time at which the
paths $\tilde{p}'_l(\T)$ and $\tilde{p}'_r(\T)$ lie in the same
corridor and consider the pincer formed by these paths after time 
$\tilde{t}_2'(\T)$ and the path joining them along the bottom
of the corridor at time $\tilde{t}_2'(\T)$.

We now repeat our previous analysis with the primed objects 
$\tilde{p}_l'(\T), \tilde{t}_2'(\T)$, {\em etc.} in place of
$\tilde{p}_l(\T), \tilde{t}_2(\T)$, {\em etc.}, checking whether
we now fall into Case (\ref{Case1}) or (\ref{Case2}); if we do not
then we pass to $\tilde{p}''_l(\T)$, {\em etc.}. We iterate this
analysis until we  fall into Case (\ref{Case1}) or 
(\ref{Case2}), at which point we acquire the desired definitions
of $\Pi_\T, t_2(\T)$ and $t_3(\T)$.
\end{enumerate}
Define $p_l(\T)$ (resp. $p_r(\T)$) to be the left (resp. right)
boundary path of the pincer $\Pi_\T$ extended backwards in
time through $\mathcal F$ to $\partial \Delta$.  Define 
$p_l^+(\T)$ to be the sequence of edges (one at each time) lying on
the leftmost of the primed $\tilde{p}_l(\T)$ from the top of
$\pi_{\T}$ to time $t_1(\T)$.
\end{definition}

\begin{definition} [cf. Definition 9.14, \cite{BG1}]
Let $\T$ be a long team of genesis (G3).  
We define $\chi_P(\T)$ to
be the set of colours containing the paths $\tilde{p}_l(\T),
\tilde{p}'_l(\T), \tilde{p}''_l(\T), \ldots$ that arise in 
Case (\ref{Case3}) of Definition 
\ref{CascadeDef} but do not become $p_l(\T)$.
\end{definition}

\begin{lemma} [cf. Lemma 9.15, \cite{BG1}]

\ 

\begin{enumerate}
\item If $\T$ is a long team of genesis (G3),
\[ t_1(\T) - t_3(\T) \le T_0(|\chi_P(\T)| + 1 ).  \]
\item If $\T_1$ and $\T_2$ are distinct teams then
$\chi_P(\T_1) \cap \chi_P(\T_2) = \emptyset$.
\end{enumerate}
\end{lemma}

\subsection{The length of teams}

This subsection follows \cite[Subsection 9.4]{BG1}.  We
consider the lengths of arbitrary teams.

\begin{definition} [cf. Definition 9.16, \cite{BG1}]
\label{Down1Def}
Let $\T$ be a team.  Define $\mbox{\rm{down}}_1(\T) \subset 
\partial \Delta$ to consist of those edges $e$ that are labelled
by some $t_i$ and satisfy one of the following conditions:
\begin{enumerate} 
\item[1.] $e$ is at the left end of a corridor containing a section of the narrow 
past of $\T$ that is not leftmost at that time; 
\item[2.]  $e$ is at the right end of a corridor containing a section of the narrow 
past of $\T$ that is not rightmost at that time;  
\item[3.]  $e$ is at the right end of a corridor which 
contains the rightmost section of the narrow past of $\T$ at that time but which does 
not intersect $p_l(\T)$.\\ 
\end{enumerate} 
\end{definition}

\begin{definition} [cf. Definition 9.17, \cite{BG1}]
Define $\partial^\T \subset \partial \Delta$ to be the intersection
of the narrow past of $\T$ with $\partial \Delta$.
\end{definition}

\begin{lemma} [cf. Lemma 9.18, \cite{BG1}]

\ 

\begin{enumerate}
\item For distinct teams $\T_1$ and $\T_2$, the sets
$\partial^{\T_1}$ and $\partial^{\T_2}$ are disjoint.
\item For distinct teams $\T_1$ and $\T_2$, the 
sets $\text{\rm{down}}_1(\T_1)$ and $\text{\rm{down}}_1(\T_2)$
are disjoint.
\end{enumerate}
\end{lemma}

\begin{definition} [cf. Definition 9.19, \cite{BG1}]
Suppose that $\T$ is a team of genesis (G3).  We define  $\QT$ be the
set of edges $\e$ with the following properties:
$p_l(\T)$ passes through $\e$  before time $t_3(\T)$,  the corridor $S$
with $\e\in\bot(S)$ contains the entire narrow past of $\T$, and
this narrow past has length at least $\lambda_0$.
\end{definition}

The following lemma reduces the task of bounding the total length
of teams to that of bounding the size of the sets $\QT$.  Its proof
follows that of \cite[Lemma 9.20]{BG1}.

\begin{lemma} [cf. Lemma 9.20, \cite{BG1}]

\

\begin{enumerate}
\item If the genesis of $\T$ is of type (G1) or (G2), then
\[ \| \T \| \le 2LC_4 |\text{\rm{down}}_1(\T)| + |\partial^\T|.      \]
\item If the genesis of $\T$ is of type (G3), then
\[ \| \T \| \le 2C_4 |\text{\rm{down}}_1(\T)| + |\partial^\T|
+ 2LC_4 |\QT| + 2LC_4T_0 (|\chi_P(\T)| + 1) + \lambda_0.        \]
\end{enumerate}
\end{lemma}

\subsection{Bounding the size of $\QT$}

Let $\mathcal{G}_3$ be the set of long teams of genesis
(G3) for which $\QT$ is nonempty. Our goal for the remainder
of this section is to find a bound for $\sum_{\T \in \mathcal{G}_3}
|\QT|$.

\begin{lemma} [cf. Lemma 9.22, \cite{BG1}]
For all $\T \in \mathcal{G}_3$
\[ t_3(\T) - t_2(\T) = \life(\Pi_\T) \le T_1(|\chi(\Pi_\T)| + 1).       \]
\end{lemma}

\begin{lemma} [cf. Lemma 9.23, \cite{BG1}]
If $\T_1, \T_2 \in \mathcal{G}_3$ are distinct teams then
$\chi(\Pi_{\T_1}) \cap \chi(\Pi_{\T_2}) = \emptyset$.
\end{lemma}
\begin{proof}
The pincers $\Pi_{\T_i}$ are disjoint or else one is 
contained in the other.  In the latter case, say $\Pi_{\T_1}
\subset \Pi_{\T_2}$, the definition of nesting (Definition
\ref{NestDef}), and of the pincer associated to a team
(Definition \ref{CascadeDef}) ensure that $\Pi_{\T_1}$
is actually nested in $\Pi_{\T_2}$ (cf. Remark 
\ref{remark:nesteddef}).
\end{proof}

\begin{corollary} [cf. Corollary 9.24, \cite{BG1}]
\[\sum_{\T \in \mathcal{G}_3} t_3(\T) - t_2(\T)
\le 3T_1 |\partial \Delta|.     \]
\end{corollary}

We have now reduced our task for this section to bounding
the number of edges in the $\QT$ which occur before
$t_2(\T)$; this is the cardinality of the following set.

\begin{definition} [cf. Definition 9.25, \cite{BG1}]
\label{Down2Def}
For a team $\T \in \mathcal{G}_3$ we define $\text{\rm{down}}_2(\T)$
to be the set of edges in $\partial \Delta$ that lie at the right-hand
end of a corridor containing an edge in $\QT$ before time $t_2(\T)$.
\end{definition}

Just as in \cite{BG1}, it is not necessarily the case that 
the sets $\text{\rm{down}}_2(\T)$ are disjoint for distinct
teams, and we must deal with the possibility of `double-counting'.

The left-to-right ordering defined on paths in $\mathcal F$
in \cite[$\S$9]{BG1} is defined in the current context exactly
as in \cite{BG1}.

\smallskip

\noindent{\bf Notation:} Let $\Gthree'$ be the set of teams $\T
\in \Gthree$ with $\down_2(\T) \neq \emptyset$.

\smallskip

\begin{lemma} [cf. Lemma 9.26, \cite{BG1}]
Consider $\T \in \Gthree'$.  If a path $p$ in $\mathcal F$ is to
the left of $p_l(\T)$ and a path $q$ is to the right of $p_r(\T)$,
then there is no corridor connecting $p$ to $q$ at any time
$t < t_2(\T)$.
\end{lemma}

\begin{definition} [cf.~Definition 9.27, \cite{BG1}]
$\T_1 \in \Gthree'$ is said to be {\em below} $\T_2 \in \Gthree'$
if $p_l(\T_1)$ and $p_r(\T_1)$ both lie between $p_l(\T_2)$ and
$p_r(\T_2)$ in the left-to-right ordering.

$\T_1$ is {\em to the left} of $\T_2$ if both $p_l(\T_1)$ and 
$p_r(\T_2)$ lie to the right of $p_r(\T_1)$.

We say that $\T$ is at {\em depth} $0$ if there are no teams above
it.  Then, inductively, we say that a team $\T$ is at depth $d+1$ if $d$ 
is the maximum depth of those teams above $\T$.

A {\em final depth} team is one with no teams below it.
\end{definition}

Note that there is a complete left-to-right ordering of those teams
in $\Gthree'$ at any given depth.

\begin{lemma} [cf. Lemma 9.28, \cite{BG1}]
If there is a team from $\Gthree'$ below a team $\T \in \Gthree'$,
then $t_1(\T) \ge \time(S_0) \ge t_2(\T)$.
\end{lemma}
\begin{proof}
The proof from \cite{BG1} works almost verbatim.  In particular,
the same proof shows that $\time(S_0) \ge t_2(\T)$.

To see that $t_1(\T) \ge \time(S_0)$, suppose that $\T'$ is a team
below $\T$.  Associated to
the team $\T'$ we have the beaded Nielsen path $\mathfrak T'$,
which is to be consumed by some reaper.  
The definitions of nesting and of the pincer $\Pi_{\T'}$ 
ensure that this
consumption of $\mathfrak T'$ must occur before time $t_1(\T)$.
On the other hand, $\mathfrak T$ has 
a non-empty future or past in $S_0$.
\end{proof}

With the preceding results in hand, a direct translation
of the proof of Lemma 9.29, \cite{BG1} finishes the
work of this section:

\begin{lemma} [cf. Lemma 9.29, \cite{BG1}]
There exist sets of colours $\chi_c(\T)$ and $\chi_\delta(\T)$
associated to each team $\T \in \Gthree'$ such that the sets
associated to distinct teams are disjoint and the following 
inequalities hold.

For each fixed team $\T_0 \in \Gthree'$ (of depth $d$ say), the
teams of depth $d+1$ that lie below $\T_0$ may be described
as follows:
\begin{enumerate}
\item[$\bullet$] There is at most one distinguished team
$\T_1$, and
\[ \| \T_1 \| \le 2B \Big( T_1(1+ |\chi(\Pi_{\T_0})|) +
T_0(|\chi_P(\T_0)| + 1) \Big).  \]
\item[$\bullet$] There are some number of final-depth teams.
\item[$\bullet$] For each of the remaining teams $\T$ we
have
\[ |\down_2(\T_0) \cap \down_2(\T)| \le T_1 \Big(
1+|\chi_c(\T)| \Big) + T_0 \Big( |\chi_\delta(\T)|+2 \Big). \]
\end{enumerate}
\end{lemma}

\begin{corollary} [cf.~Corollary 9.30, \cite{BG1}]
Summing over the set of teams $\T \in \Gthree'$ that are
not distinguished, we get
\[ \sum_\T \Big| \down_2(\T) \Big| \le
2 \Big| \bigcup_\T \down_2(\T)\Big| + \sum_\T T_1\Big(
1 + |\chi_c(\T)| \Big) + \sum_\T T_0 \Big( |\chi_\delta (\T)|
+ 2 \Big).      \]
\end{corollary}

Summing over the same set of teams again, we finally obtain:

\begin{corollary}
\[      \sum_\T |\down_2(\T)| \le |\partial \Delta| (2 + 3T_1 + 5T_0).  \]
\end{corollary}

\section{The Bonus Scheme} \label{BonusSection}

This section closely follows \cite[Section 10]{BG1}.  
We have at last reached a stage where the proofs from
\cite{BG1} can be translated without significant
modification.

In the previous section we defined teams and obtained a 
global bound on $\sum \| \T \|$.  If $C_{(\mu,\mu')}(2)$ is
non-empty then $(\mu,\mu')$ is a member or virtual member
of a unique team.  If the team is such that $t_1(\T) \ge \time(S_0)$,
then no member of the team is virtual and we have the inequality
\[ \| \T \| \ge \sum_{(\mu_1,\mu_2) \in \T} |C_{(\mu_1,\mu_2)}| - B(B+1),
\]
established in Lemma \ref{TeamBound}.  
This inequality might fail in case $t_1(\T) < \time(S_0)$.  The {\em
bonus scheme} assigns additional edges to teams in order to 
compensate for this failure.

By definition, at time $t_1(\T)$ the reaper $\rho = \rho_\T$ lies
immediately to the right of $\mathfrak T$.  The beads of $\mathfrak T$
not consumed from the right by $\rho$ by $\time(S_0)$ have a 
preferred future in $S_0$.  This preferred future, if contained in a single
colour, lies in $C_{(\mu_1,\mu_2)}(2)$ for some member 
$(\mu_1,\mu_2) \in \T$.  It could also intersect more than one colour
\footnote{Since Nielsen beads have bounded length, and there is a bound on the number of adjacencies of colours, there are relatively few such beads.}.
However, not all beads in the
 $C_{(\mu_1,\mu_2)}(2)$ need arise in this way:  some may not
 have a Nielsen bead as an ancestor at time $t_1(\T)$.  And if 
 $(\mu_1,\mu_2)$ is only a virtual member of $\T$, then no bead of
 $C_{(\mu_1,\mu_2)}(2)$ lies in the future of $\mathfrak T$.  The
 {\em bonus} beads in $C_{(\mu_1,\mu_2)}(2)$ are a certain subset of 
 those that do not have a Nielsen bead as an ancestor at time $t_1(\T)$.
 They are defined as follows.

\begin{definition}
Let  $\T$ be a team with  $t_1(\T) < \time(S_0)$ and consider a time
$t$ with $t_1(\T) < t < \time(S_0)$.

The {\em swollen future} of $\T$ at time $t$ is defined as in
Definition \ref{Swollen} with respect to the interval $\mathfrak{T}$,
which lies at time $t_1(\T)$.

Let $\epsilon$ be a non-Nielsen bead that lies immediately to the
left of the swollen future of $\T$, but whose immediate
ancestor is not a right
linear edge in this position.  If the path from $\epsilon$ to the reaper
$\rho_\T$ of $\T$ is a \gep, then we say that $\epsilon$ is a
{\em rascal}.  Otherwise, if $\epsilon$ provides
more Nielsen beads than the reaper consumes, then $\epsilon$ 
is a {\em terror}.

In both cases, the {\em bonus provided by $\epsilon$} is the 
set of beads in the swollen future of $\T$ in $S_0$ that have
$\epsilon$ as their most recent ancestor which is not a 
Nielsen bead,
and which are eventually consumed by $\rho_\T$.

The set $\bonus(\T)$ is the union of the bonuses provided
to $\T$ by all rascals and terrors.
\end{definition}

\begin{lemma} [cf. Lemma 10.2, \cite{BG1}]
For any team $\T$,
\[      \sum_{(\mu_1,\mu_2) \in \T\ or\ (\mu_1,\mu_2) \vin \T}
| C_{(\mu_1,\mu_2)}(2)| \le \| \T \| + |\bonus(\T)| +B +J.      \]
\end{lemma}
 Note that the \gep \ which contains
a rascal in the above definition is not displayed.
We now proceed to bound the total bonus provided 
to teams by all
rascals and terrors.  Terrors
are straightforward to deal with.

\begin{lemma} [cf. Lemma 10.3, \cite{BG1}]
The sum of the lengths of the bonuses provided to all
teams by terrors is less than $2L\n$.
\end{lemma}
\begin{proof}
Let $\epsilon$ by a terror, associated to a team $\T$.
Since the region from $\epsilon$ to the reaper of $\T$ is not
a \gep, $\epsilon$ must be right-fast.  Therefore, it will be 
separated from the team to which it is associated after one unit
of time.  Hence the bonus that $\epsilon$ provides is at most $L$.

That there can be at most one terror per 
adjacency of colours follows in a straightforward manner from
Lemma \ref{EndStab2} and the definition of terror.

Thus the total contribution of all terrors is less than $2L \n$.
\end{proof}

  In
parallel with \cite[Definition 10.4]{BG1}, we make the following

\begin{definition}
Fix a team $\T$ with $t_1(\T) < \time(S_0)$ and consider the
interval of time $[\tau_0(\epsilon),\tau_1(\epsilon)]$, where 
$\tau_0(\epsilon)$ is the
time at which a rascal $\epsilon$ appears at the left end of the 
swollen future of $\T$, and $\tau_1(\epsilon)$ is the time at 
which the robust future of $\epsilon$ is no longer to the immediate
left of the future of the swollen future of $\T$.

In the case where the robust future $\hat{\epsilon}$ of $\epsilon$ at 
time $\tau_1(\epsilon)$ is cancelled from the left by an edge
$e'$, we define $\tau_2(\epsilon)$ to be the earliest time when the
pasts of $\hat{\epsilon}$ and $e'$ are in the same corridor.  The path
in $\mathcal F$ that traces the past of $\hat{\epsilon}$ is denoted
$p_\epsilon$ and the past following the ancestors of $e'$ from 
$\tau_2(\epsilon)$ to $\tau_1(\epsilon)$ is denoted $p'_\epsilon$.  The 
pincer\footnote{we include the degenerate case here where the
``pincer" has no colours other than those of $\epsilon$ and $e'$.} 
formed by $p_\epsilon$, $p'_\epsilon$ and the corridor joining
them at time $\tau_2(\epsilon)$ is denoted $\Pi_\epsilon$.
\end{definition}

The only essential difference between the above definition and 
\cite[Definition 10.4]{BG1} is the use of the robust future of 
$\epsilon$ rather than the pp-future. 

With this definition in hand, the remaining results from
\cite[Section 10]{BG1} may be translated directly,
yielding in particular:

\begin{proposition} [cf. Lemma 10.13, \cite{BG1}]
Summing over all teams that are not short, we have
\[ \sum_{\T}|\bonus(\T)| \le \Big( (B+3)(3T_1 + 2T_0)L +
6BT_1 + 4BT_0 + 2\lambda_0 + 2B + 5L + 1 \Big) \n .\]
\end{proposition}

\section{From Bead Norm to Length} \label{LongGepsandPepsSection}

The output of the results up to now is a
bound for the bead norm of our corridor $S_0$.  In order to
complete the proof of Theorem \ref{BoundS} in the case of the
specified IRTT $f$ (which implies
Theorem \ref{MainThm}) we need to turn this into a bound
on the length of $S_0$.  For this we need to bound the total
length of the \gep s and \pep s in $S_0$ which have length more
than $J$ (or indeed any other fixed length).  In this section we 
explain how the techniques of the bonus scheme can be used
to establish such a bound. 

If a bead $\rho$ in $\mu(S_0)$  has length greater 
than $J$, it is either a \gep\ or a \pep.  
If it is a \pep then we may trace its past: at each time,
this past is either of length at most $J$ or 
else is a \pep\ or a \gep.  Whilst
this past remains a \pep, the number of Nielsen paths will decrease with
each 
backwards step in time, so at some point in the past of $\rho$, 
it must become
a \gep.  

Suppose now that $\rho$ is a \gep.  The past of a \gep\ is either a 
\gep\ or else has length at most $J$.  Thus, the length of the
\gep\ decreases as we go into the past until eventually it is
of length at most $J$.

There is a strong analogy between teams of genesis (G4) and long 
\gep s and \pep s.  On
one end of a long bead is a linear edge which consumes the
Nielsen beads in the middle.  This linear edge can be considered
as a reaper.  On the other end of a \gep\ is a linear edge which can
be considered as a rascal.   The moment when the past of a \pep\ becomes a
\gep\ is analogous to $\tau_1(\epsilon)$ from the bonus scheme,
and so a \pep\ in $S_0$ can be thought of as a team with a rascal 
$\epsilon$ with $\tau_1(\epsilon) \le \time(S_0)$.  Similarly, a long 
\gep\ in $S_0$ can be thought of as a team with a rascal $\epsilon$
so that $\tau_1(\epsilon) > \time(S_0)$.  

We can define the bonus associated to such a rascal exactly
as we did in the previous section.  Since we
are in the setting of  genesis  type (G4), all of the Nielsen beads in a long
\gep\ or \pep\ are in the bonus.  Thus it is enough to bound the 
total of the bonuses associated to long \gep s and \pep s.

The only thing we need to be able to follow the bonus
scheme directly is a bound on the number of long \gep s
and \pep s in $S_0$.

\begin{lemma}
The number of beads of length greater than $J$ in
$S_0$ is less then $4 \n$.
\end{lemma}
\begin{proof}
Let $\rho$ be a bead in $S_0$ of length greater than $J$, and
assign a time $\tau_1(\rho)$ to $\rho$ as described above.  If
$\rho$ is a \gep\ then $\tau_1(\rho) > \time(S_0)$, whilst if 
$\rho$ is a \pep\ then $\tau_1(\rho) \le \time(S_0)$.

Let $\rho'$ be the past or future of $\rho$ at time $\tau_1(\rho) -1$.
Consider the `event' at time $\tau_1(\rho)$ which stops the robust
future of $\rho'$ being a \gep.  

This `event' is either an intrusion of the boundary, a singularity, or else
there is an associated pincer caused by a cancellation from another
colour.  There are less than $\n$ events of each of the first two
types.

The Buffer Lemma ensures that there is at most one event of the
third type for each adjacency of colours.  An application of 
Lemma \ref{NoOfAdj} completes the proof.
\end{proof}

A bound on the total length of long beads in $S_0$ now follows
exactly as in the bonus scheme from Section \ref{BonusSection}
(the detailed arguments being in \cite[Section 10]{BG1}).

\subsection{The end of the main road}

In Section \ref{Linear} we discussed how Theorem \ref{MainThm}
follows from Theorem \ref{noName} and Proposition \ref{longGEPS}.
The bound that we just established
on the total length of long beads in $S_0$ proves
Proposition \ref{longGEPS}. The output of our estimates in 
the previous sections bounded the bead norm of $S_0$ by a linear
function of $|\partial\Delta|$, and Theorem \ref{noName} follows
from this because 
\[	[S]_{\beta} \le B \| S \|_\beta	,	\]
(see Lemma \ref{BLengthNorm}).

Thus the proof of
Theorem \ref{MainThm} is finally at an end, and the reader can
join us in wondering why a statement as simple and engaging
as Theorem \ref{BoundS} should require such a complicated proof.

\section{Corridor Length Functions and Bracketing} \label{BracketingSection}

In this section we prove Theorem \ref{BoundS} in full generality and
deduce the Bracketing Theorem from it. Our proof of Theorem \ref{BoundS}
proceeds via a discussion of {\em{corridor length functions}} for
more general semidirect products and mapping tori. 
Such functions should be regarded as measuring the complexity of
van Kampen diagrams in the spirit of isoperimetric and isodiametric
functions. We prove the following results (see Subsection \ref{ss:defL}
for precise definitions of the terms involved).

\begin{proposition}\label{CLchange} Let $G_1$ and $G_2$ be 
compact combinatorial complexes with
fundamental group $\Pi$,
and for $i=1,2$ let $f_i:G_i^{(1)}\to G_i^{(1)}$ be
an edge-path map of 1-skeleta inducing $\phi\in\rm{Out}(\Pi)$. Then the
$t$-corridor length function for the mapping torus $M(f_1)$
is $\simeq$ equivalent to that of $M(f_2)$.
\end{proposition}

\begin{proposition}\label{CLfi} If $\Pi$ is finitely generated and 
$\G=\Pi\rtimes_\phi \Z$ is finitely presented, then for every positive
integer $p$, the corridor length function of $\Pi$ is $\simeq$ equivalent
to that of
$\G_p=\Pi\rtimes_{\phi^p} \Z$
\end{proposition}

In the previous section we completed the proof of 
Theorem \ref{BoundS} in the case of one particular IRTT representative $f$
of a certain power of an arbitrary free-group automorphism $\phi$.  
The above results complete the proof in the general case.
Before turning to the proof of these results, we explain
how the Bracketing Theorem stated in the introduction 
is obtained by applying Theorem \ref{BoundS} to the most
naive topological representation of a free group automorphism $\phi$.

\subsection{The Bracketing Theorem}

The terms in   the following theorem were defined in the
introduction.

\medskip
\noindent{\bf Theorem \ref{Bracket}.}
{\em
There exists a constant $K = K(\phi,\mathcal B)$
such that any word $w \equiv e_1 \dots e_n$ that represents
the identity in $F\rtimes_\phi\Z$ 
admits a $t$-complete bracketing $\beta_1,\dots , \beta_m$
such that the content $c_i$ of each $\beta_i$
satisfies $d_F(1,c_i)\le Kn$.
}

\smallskip

\begin{proof}
We work with the mapping torus $M$
of the obvious realisation of $\phi$ on
the graph with one vertex whose edges are indexed by $\mathcal B$.
Given a word $w$, we consider
a minimal-area van Kampen diagram over $M$ with boundary
label $w$. 
 We insert a bracket 
$w_1(w_2)w_3$ if and only if there is a $t$-corridor whose
ends are labelled by the initial and terminal letters of $w_2$.
(One must allow $t$-corridors of zero length in this description;
one would exclude them by making the easy reduction to words that
have no proper sub-words that are null-homotopic.)

These brackets are pairwise compatible because distinct $t$-corridors
cannot cross. And
because every $t$-edge in the boundary of a van Kampen
diagram is the end of a (perhaps zero-length) corridor, the
bracketing is complete. The content of the bracket is the freely
reduced form of the label  along the top or bottom of the corridor
(according to the orientation of the sentinels). In the former case, the
length of the corridor bounds the length of this label, and in the
latter case one has to multiply the length by at most $L=\max\{|\phi(b)|:b\in
\mathcal B\}$.
\end{proof}

\subsection{Corridor length functions} \label{ss:defL}
If $\Pi$ is a group with finite generating
set $\A$ and $\phi\in\rm{Aut}(\Pi)$ is such that $\G=\Pi\rtimes_\phi\Z$
is finitely presented, then $\G$ has a finite presentation of the
form 
\[	\big\langle \A,\,t \mid \mathcal R,\, t^{-1}at=\hat\phi(a)\ (a\in\A)\big\rangle ,	\]
where $t$ is the generator of the visible $\Z$, the relations $\mathcal R$
involve only the letters $\A$, and $\hat\phi(a)\in F(\A)$ is equal to
$\phi(a)$ in $\Pi$. 

We are concerned with the geometry of $t$-corridors
in van Kampen diagrams over such presentations. Thus we associate to the
presentation the {\em{$t$-corridor length function}} $\Lambda :\mathbb N\to
\mathbb N$, which is defined as follows. For each $w\in F(\A\cup\{t\})$
with $w=1$ in $\G$, we choose a van Kampen diagram for $w$ in which the
length of the longest $t$-corridor is as small as possible, and we
define $\lambda_t(w)$ to be this length. We then define
$$
\Lambda(n) := \max \{\lambda_t(w) \mid w=_\G 1,\ |w|\le n\}.
$$

More generally, since we have a well-defined notion of van Kampen
diagram and $t$-corridor in the setting of mapping tori of
{\em{edge-path maps}}\footnote{an {\em{edge-path map}} is a cellular
map that sends edges to edge-paths} of
combinatorial complexes, we can define the {\em{$t$-corridor length function}}
for such a complex.

\subsection{Invariance under change of topological representative}

The scheme of the following proof follows the standard method of
 showing that
features of the geometry of van Kampen diagrams are preserved under
quasi-isometry. However, one has to be careful to deal only with fibre-preserving
maps in order to retain control over the $t$-corridor structure.

\smallskip

\begin{proof}[Proof of Proposition \ref{CLchange}.]
We have a cocompact action of $\G =\Pi\rtimes_\phi\Z$
on the universal cover $X_i=\tilde M(f_i)$ for $i=1,2$, where the
action of $\Pi$ leaves invariant the connected components $C_{i,m}$
of the preimage
of $G_i\subset M(f_i)$ and the generator $t$
of $\Z$ acts so that $t^r.C_{i,m} = C_{i,m+r}$.

The cocompactness of the actions means that there exist constants $\delta_1, \delta_2$
so that every vertex in $C_{i,m}$ is within a distance $\delta_i$ of 
any $\Pi$-orbit of vertices in $C_{i,m}$, where distance is measured
in the combinatorial metric on the 1-skeleton (unit edge lengths).

We define $\G$-equivariant quasi-isometries between the 1-skeleta
of the $X_i$  as follows. First we
pick  base vertices $x_i\in C_{i,0}$ and define
$g_1:\gamma.x_1\mapsto \gamma.x_2$ and $g_2:\gamma.x_2\mapsto \gamma.x_1$.
Then, for each  vertex  $v\in C_{i,m}\smallsetminus \Gamma.x_i$ we choose
a closest element  $v'\in \Gamma.x_i\cap C_{i,m}$
and define $g_i(v) := g(v')$. Next, we extend to the edges in 
$C_{i,m}$ by sending each to a shortest edge path connecting the
images of its vertices. Finally, we extend  $g_i$ to $t$-edges
in $X_i$ so that it sends each such homeomorphically onto
the $t$-edge joining the images of its endpoints.

With the maps $g_1, g_2$ in hand, we can now
push van Kampen diagrams back and forth between $X_1$ and
$X_2$ as in the standard proof of the qi-invariance of Dehn functions
(cf.~\cite{BH}, page 143).
Thus, 
given a loop $\ell$ in the 1-skeleton of $X_1$, labelled $u_1t^{\e_1}u_2
\dots u_lt^{\e_l}$
we consider the loop $ g_1\circ \ell$ in $X_2^{(1)}$
and fill it with a
van Kampen diagram $\Delta$ so as minimize the length
of the longest $t$-corridor. We will be done if we can bound
$\lambda_t(\ell)$ by a linear function of this length.

Viewing $\Delta$ as a
map from a cellulated 2-disc to $X_2$, we compose it with $g_2$
to obtain a map to $X_1$. This new map is obtained from $\Delta$
by simply changing the labels on the edges: the $t$-edges are unchanged
while the edges labelled by 1-cells in $G_2$ are now labelled by 
edge-paths in the 1-skeleton of $G_1$ whose length is
bounded by the constants of the quasi-isometry $g_2$; the boundary
label of the diagram will be $\ell'=v_1t^{\e_1}v_2
\dots v_lt^{\e_l}$, where the $v_j$ are edge-paths of uniformly bounded
length and each $v_j$ is contained in the same component $C_{1,m_j}$
as $u_j$. (This is the point at which we use the fact that we chose
our quasi-isometries to respect fibres.) The faces of this diagram
can be filled with van Kampen diagrams in $X_1$; in the case of 2-cells
with no $t$-labels, we use only lifts of 2-cells from $G_1$; in the
case of 2-cells labelled $t^{-1}\rho t \sigma$ we divide them into 
(short) $t$-corridors in the obvious manner. 
The result\footnote{A familiar problem in this type of argument
arises from degeneracies that threaten the planarity of the
diagram; such problems are removed by surgery \cite{LS}. In
the current setting these surgeries take place only in the regions
between the $t$-corridors and therefore do not affect our discussion.} is a van
Kampen diagram for $\ell'$ in $X_1$ whose $t$-corridors are in bijection
with those of $\Delta$ and whose length is bounded by $k$ times the
length of those in $\Delta$, where $k$ is a constant that depends only
on our quasi-isometries.

To complete the desired diagram filling our original loop $\ell$, 
we need an annular diagram between $\ell$ and $\ell'$ that does not
 disrupt the structure of $t$-corridors in $\Delta'$. To this end,
 we join the vertices of $u_j$ to those of $v_j$ by paths in $C_{i,m_j}$
 of minimal length and fill the resulting loop with a diagram mapping to
 $C_{i,m_j}$; this gives a diagram
 $\Delta''$ with holes corresponding to the
 occurrences of $t^{\pm 1}$ in $\ell$.
 Next, if the
 arc joining the termini of $u_j$ and $v_j$ is labelled $\rho_i$,
 then we  insert a $t$-corridor into the
 hole associated to $\dots u_jt u_{j+1}\dots$, where the bottom
 of the $t$-corridor is labelled $\rho_j$. (If $t$ is replaced by $t^{-1}$,
 the bottom of the corridor is the arc $\sigma_{j+1}$ joining the
   initial vertex of $u_{j+1}$
 to that of $v_{j+1}$.) To complete the construction of $\Delta$, one
 uses 2-cells in $C_{i,m_{j+1}}$ to fill the loop formed by the top of
 the $t$-corridor and $\sigma_{j+1}$.
 \end{proof}

\begin{corollary} If $\Pi$ is finitely generated and $\G=\Pi\rtimes_\phi\Z$
is finitely presented then,
up to $\simeq$ equivalence, the $t$-corridor length function of 
$\Pi\rtimes_\phi\Z$ depends only on the semidirect product (i.e.~
although it
depends on the form of the finite presentation, it does
 not depend on the
choice of $\A$ and $\hat\phi$).
\end{corollary}

\subsection{Passing to Powers} 

The purpose of this subsection is to prove Proposition \ref{CLfi}. 

Let $\At$ be as above. Identifying
 $\G_p=\Pi\rtimes_{\phi^p}\Z$ with the subgroup 
$\Pi\rtimes p\Z$ of $\G$, we take generators $\A\cup\{\tau\}$
where $\tau=t^p$ in $\G$.
To each word $w\in\FAt$ that
equals $1\in\G$ we associate a word $w_p$ in the free group
on $\A\cup\{\tau\}$ according to the following
scheme. First
 we draw a path on the integer lattice in $\mathbb R^2$ that begins at the
origin and proceeds up one space as we read $t$, down one as we read $t^{-1}$
and moves one space to the right as we read a letter from $\Apm$.
We shall modify $w$ by replacing certain open 
segments of this path that lie in the vertical
 intervals $[mp,(m+1)p]$; these segments
are of two types, called {\em{bumps}} and {\em{steps}}. 

If both endpoints of the subpath are at height $mp$
and none of its edge are at height $(m+1)p$, then the segment
is called an {\em{up-bump}}. If the initial
endpoint is at height $mp$, the terminus at
height $(m+1)p$ and all other vertices are at 
heights in $(mp, (m+1)p)$,
then the segment is called an {\em{up-step}}. A {\em{down-bump}}
and {\em{down-step}} are defined similarly.

When we have replaced all steps and bumps from the path defined by $w$,
  the horizontal segments of the
resulting path will all run
at heights divisible by $p$.

To this end, we write $w=u_1v_1u_2v_2\dots$ where
$u_1$ is the first non-trivial prefix of $w$
whose exponent sum in $t$ is $0\mod p$ and $v_1$ is the (possibly empty)
subword  before the next $t^{\pm 1}$, then $u_2$
 is the first non-trivial prefix of $w$
whose exponent sum in $t$ is $0\mod p$, and so on.
Each $u_i$ labels either a bump or a step.

If $u_i$ labels a bump
then  we replace it by the reduced word $
U_i\in \FA$ that is equal in $\G$ to $u_i$.
 If $u_i=t^\e u_i',\, \e=\pm 1,$ is a step, then
we replace it by the unique reduced word $t^{\e p}U_i$ with
$U_i\in\FA$ and $t^\e U_i=u_i$ in $\G$.

Let $\tilde w_p\in\FAt$ be the word obtained from $w$ by the above process and
let $w_p\in\FAtp$ be the word obtained from $\tilde w_p$ by (starting from
the left) replacing sub-words labelled $t^{\pm p}$ by $\tau^{\pm p}$
and then freely reducing.

As usual, in the following lemma $L=\max\{|\phi(a)| : a\in\mathcal A\}$.

\begin{lemma}\label{lengthwp} $w=\tilde w_p=w_p$ in $\G$
and $ |w_p| \le |\tilde w_p| \le L^{p-1} |w|$.
\end{lemma}

\begin{proof} The bound on $|\tilde w_p|$ comes from the following observation.
For  a bump labelled $u_i$,
one can pass from $u_i$ to $U_i$ by deleting all letters $t^{\pm 1}$
from $u_i$ and replacing each occurrence of $a\in \A$ in $u_i$, 
say $u_i=\alpha a\beta$, by the freely
reduced word in $\FA$ representing $\phi^r(a)$, where $-r$ is the exponent
sum of $t$ in $\alpha$. Similarly, if a step is labelled $u_i=t^\e u_i'$,
then $U_i$ is obtained by deleting all $t$ from $u_i'$ and replacing
 each occurrence of $a\in \A$ in $u_i$,
say $u_i'=\alpha a\beta$, by the freely
reduced word in $\FA$ representing $\phi^{\e(p-r)}(a)$, 
where $\e r$ is the exponent
sum of $t$ in $\alpha$.
\end{proof}

The replacement scheme described in the preceding proof corresponds
to the construction of a singular-disc diagram  $A(w)$ exhibiting
the equality $w = \tilde w_p$ in $\G$. 
 Specifically, for each bump or
step, one draws the vertical line joining each vertex to the height where
it will be pushed, one labels it by the appropriate power of $t$, and then
one
fills-in the resulting line of rectangles  with 2-cells whose boundary
labels have the form $t^{-1}at\phi^{-1}(a)$.
(Starting from this
specific planar embedding one will in general
have to flip some of the components of the interior in order to get
an embedded diagram $A(w)$ with boundary cycle $\tilde w_p w_p^{-1}$.)

\begin{lemma}\label{corrinAw} $A(w)$ is a union of $t$-corridors;
each has at most one of its ends on the boundary arc labelled $\tilde w_p$,
and the length of a $t$-corridor in $A(w)$ is at most
  $L^{p-1}\max |u_i|$, where the $u_i$ are the sub-words
of $w$ labelling  bumps and steps.
\end{lemma}

\begin{proof} The diagram $A(w)$ consists of a string of disc diagrams, one
for each bump or step. A  $t$-corridor in a disc corresponding to
a bump labelled $u_i$ has both of its ends on the arc labelled $u_i$,
while a $t$-corridor in a disc corresponding to a step
labelled $tu_i'$ may have one end on the corresponding arc labelled $t^p$
in $\tilde w_p$ and one on the arc labelled $u_i'$ or (if the change in height
along $u_i'$ is not monotone) both ends on the arc labelled $u_i'$. In all
cases, the label on the bottom side of the corridor is a concatenation of
 less than $|u_i|$ words of the form $\phi^r(a)$ with $a\in\A$ and $|r|\le
p-1$.
\end{proof}

\noindent{\bf{Proof of Proposition \ref{CLfi}.}}
As we discussed immediately before subsection \ref{Iterate},
 the set of diagrams for $\G_p$ is,
after $p$-refinement, a subset of the diagrams over $\G$,
and hence the corridor length function of the latter
$\preceq$-dominates that of the former. (There are some
constants to take account of here, such as a factor
of $p$ in length
coming from the $p$-refinement, and an $L^{p-1}$ needed
to estimate the area of a $t$-corridor in terms of the corresponding
$\tau$-corridor, but these are trivial matters.)
Thus the true content of the proposition is that  
the corridor length function of $\G$ is
$\preceq$-bounded above by that of the $\G_p$.

For each freely-reduced
 word $W\in \FAtp$ that is null-homotopic in $\G_p$
we fix a van Kampen diagram $\Delta(W)$ whose $\tau$-corridors have length
at most $\Lambda(|W|)$. Then, 
for each freely-reduced $w\in \FAt$ that is 
null-homotopic in $\G$ we define a van Kampen diagram $\Delta_p(w)$ 
as follows. First, we replace
 $\Delta(w_p)$ by its $p$-refinement (which
has boundary label $\tilde w_p$). We then attach
to this the singular-disc diagram $A(w)$ along the portion of
its boundary labelled $\tilde w_p$.

We claim that the length of each
$t$-corridor in $\Delta_p(w)$ is at most
$$
L^{p-1}\, (2+ \Lambda(L^{p-1}|w|)).
$$

It follows from Lemma \ref{corrinAw} that 
each of the 
$t$-corridors in $\Delta_p(w)$ is either contained in 
the annular diagram $A(w)$, or else is a layer in the
$p$-refinement of a $\tau$-corridor from $\Delta(w_p)$,
possibly augmented on each end by a $t$-corridor in $A(w)$.
(The fact that there are no $t$-corridors in $A(w)$ with
both ends on the boundary arc labelled $\tilde w_p$ is
crucial here.)

The length of a $t$-corridor in $A(w)$ is at most
$L^{p-1}|v|$. The length of a $\tau$-corridor from
$\Delta(w_p)$
is at most $\Lambda(|w_p|) \le \Lambda(L^{p-1}|w|)$, and the
length of each layer in its refinement is therefore at most
$L^{p-1}\, \Lambda(L^{p-1}|w|)$.
\hfill$\square$
  
\appendix\section{On a Result of Brinkmann} \label{BrinkSection}

The following theorem is the main result in Peter Brinkmann's
paper \cite{BrinkDyn}. It plays a vital role in the first proof
that the conjugacy problem is solvable for free-by-cyclic groups
\cite{BMMV}  (our Corollary \ref{Conj}).

\begin{theorem} \label{Brink} \cite[Theorem 0.1]{BrinkDyn}
Let $\phi : F \to F$ be an automorphism of a finitely generated
free group.  Then there exists a constant $K \ge 1$ such that for
any pair of exponents $N, i$ satisfying $0 \le i \le N$, the following
two statements hold:
\begin{enumerate}
\item\label{Brink1} If $w$ is a cyclic word in $F$, then
\[ \| \phi^i(w) \| \le K \Big( \| w \| + \| \phi^N(w) \| \Big) ,        \]
where $\| w \|$ is the length of the cyclic reduction of $w$ with
respect to some word metric on $F$.
\item \label{Brink2} If $w$ is a word in $F$, then
\[ | \phi^i(w)| \le K \Big( |w| + |\phi^N(w)| \Big),    \]
where $|w|$ is the word length of $w$.
\end{enumerate}
\end{theorem}

The purpose of this appendix is to explain how to extract
Theorem \ref{Brink} from our proof of Theorem \ref{MainThm}.
We regard words and cyclic words in $F_n$
 as, respectively, based and unbased loops
in the graph $R$ with one vertex and $n$ edges; the assertions
of Theorem \ref{Brink} are then statements about how the 
lengths of the tightened
images of such loops grow when one applies the obvious
topological realisation $\overline\phi$ of $\phi$.
As in the previous subsection, these assertions will follow
if we can establish the corresponding bounds with $\overline\phi:R\to R$
replaced by a topological (IRTT) representative $f:G\to G$ of
a power of $\phi$ satisfying Assumption
\ref{FinalPower}.

\begin{remark} The proof given below shows that the constant $K$
of Theorem \ref{BoundS} suffices for Theorem \ref{Brink}.
Brinkmann \cite{BrinkDyn}
states that (his constant) $K$ can be computed
effectively, but we do not see how to prove this. Indeed,
given his approach (and ours), this assertion would
seem to require an effective construction of an improved
relative train track representative for $\phi$, and a
proof that such a construction exists does not
seem to be available at the moment.
\end{remark}

 The following lemma allows a proof of the assertions in
\eqref{Brink1} and \eqref{Brink2} to be undertaken simultaneously.

\begin{lemma}
If $\sigma$ is a nontrivial loop in $G$, then
for some $j\ge 1$, the loop $f_{\#}^j(\sigma)$ admits a 
splitting at a vertex.
\end{lemma}
\begin{proof}
According to \cite[Lemma 4.1.2, p.554]{BFH}, $\sigma$ admits a splitting
$\sigma = \sigma_1$, where $\sigma_1$ is a path, but we argue further
to arrange for this splitting to be at a vertex.

We divide the argument into a number of cases, 
depending on the largest $i$ so that the stratum $H_i$
contains an edge of
$\sigma_1$.  If this $H_i$
is a zero stratum,  $f_{\#}(\sigma_1) \subset G_{i-1}$
and an obvious induction applies.     If $H_i$
parabolic, then we apply
\cite[Lemma 4.1.4]{BFH} to the circuit $\sigma$ to obtain a splitting into paths,
at least one of which is a basic path, and so has a vertex at one end. 
If $H_i$ is an exponential stratum, then there is a positive integer
$K$ so that the number of $i$-illegal turns in $f^k_{\#}(\sigma_1)$ is the 
same for all $k \ge K$.  In this case, since all Nielsen paths of exponential
weight are edge-paths and all periodic paths are Nielsen, \cite[Lemma 4.2.6]{BFH} 
implies that $f^K_{\#}(\sigma_1)$ admits a splitting into sub-paths which are either 
$r$-legal or pre-Nielsen paths. 
If all sub-paths of $f^K_{\#}(\sigma_1)$ are pre-Nielsen paths, then 
$f^{K+1}_{\#}(\sigma_1)$ is a Nielsen path, and we ensured in 
\cite[Section 1]{BG2} that all Nielsen paths are edge-paths.

Suppose, then, that $f^K_{\#}(\sigma_1)$ contains an $r$-legal path $\rho$
of weight $r$ in its splitting.  Then an iterate $f^i_{\#}(\rho)$ of $\rho$ contains a 
displayed edge $\e$ of weight $r$, and the path $f^{K+i}_{\#}(\sigma_1)$
splits immediately on either side of $\e$.  
Since $\sigma$ has weight $i$, the splitting of $f^{K+i}_{\#}(\sigma_1)$ 
induces a splitting of $f^{K+i}_{\#}(\sigma)$ at a vertex, as required.
\end{proof}

In order to prove the statements \eqref{Brink1} and \eqref{Brink2}, 
we analyze the van Kampen diagram $\Delta$ over the
mapping torus of $f:G\to G$ that has boundary label
$t^{-k}\sigma t^k f_{\#}^k(\sigma)^{-1}$. This is a
simple stack of corridors as consider in Subsection \ref{stackDiags}.

In the restricted setting of stack diagrams, 
many of the difficulties that had to be overcome in the proof of 
Theorem \ref{MainThm} do not arise
(there are no singularities, for example), but
there remain difficulties that one does not encounter
in the context of positive automorphisms.

The number of edges in  $\partial\Delta$ not
labelled
$t$ is the quantity that determines the upper bound we seek,
 $n := |\sigma| + |f^N(\sigma)|$).
We must bound the length of each corridor in $\Delta$
linearly in terms of $n$.  Theorem \ref{BoundS} provides a
bound in terms of  $\n$, so we must argue is that in the context
of stack diagrams, one can dispose of the contribution of the
$t$-edges to this bound.
In order to do so, we make an exhaustive list of those
places in the  proof of Theorem
\ref{BoundS} where $t$-edges were accounted for, 
and we explain why, in each case, 
they are not required in the setting of simple stack diagrams.

(1) The $t$-edges contributed to the bound on the size of 
$S_0(2)$ and $S_0(3a)$ in Section \ref{FastBeads}, but
these sets do not arise in stack diagrams.

(2) The $t$-edges were required in determining the
sets $\down_1(\T)$ used to bound the
lengths of teams (see Definition \ref{Down1Def}).
But $\down_1(\T)$ was
used only to bound the lengths of those
teams whose narrow past had several components at 
some time in the past, and this cannot happen in a stack
diagram.

(3) The $t$-edges entered the definition of $\down_2(\T)$,
which was
used to bound the
number of edges in $\QT$ before time $t_2(\T)$ (see Definition 
\ref{Down2Def}). But
there are no such edges in a stack of corridors, so
we do not have to worry about double-counting, and an
improved bound
on the lengths of teams can be derived directly from the
Pincer Lemma, noting that there are less than
$2|\partial \Delta|$ adjacencies of colours.

(4) In the bonus scheme, the set $\partial^e$ is used to 
bound the size of the interval of time $[\tau_0(e),\tau_2(e)]$,
but in a stack of corridors it is clear that $\tau_0(e)
=\tau_2(e)$, so the edges $\partial^e$ are not required.

(5) Likewise, when bounding the size of the bonuses provided by
rascals, we do not need to use
the edges $\down_2(e)$ if our diagram is simply a stack of corridors

(6) A final use of $t$-edges is hidden in our references to \cite{BG1}
in the implementation of the Bonus scheme, specifically
the bound on the sum of the lengths of blocks satisfying condition
(iv) of the `tautologous tetrad'. This is unnecessary in stack
diagrams because there are no singularities and
no edges that are cancelled by edges from outside the future of
$S_0$, so the paths $\pi_l$ and $\pi_r$ travel forwards
in time  until they hit the boundary and
$\sum |\text{bdy}(\mathfrak{B})|<n$ bounds the size of the sum of
all such blocks.
$\square$


\begin{thebibliography}{99}
\bibitem{Alonso} J. Alonso, In\'egalit\'es isop\'erim\'etriques et quasi-isom\'etries, \textit{C. R. Acad. Sci. Paris}, {\bf 311} (1990), 761-764.
\bibitem{BFH} M. Bestvina, M. Feighn and M. Handel, The Tits
alternative for $Out(F_n)$ I: Dynamics of exponentially growing
automorphisms, \textit{Ann. of Math. (2)}, {\bf 151} (2000), 517--623.
\bibitem{BFH2} M. Bestvina, M. Feighn and M. Handel, The Tits alternative for $\text{\rm Out}(F_n)$ II: A Kolchin type theorem, \textit{Ann. of Math. (2)}, {\bf 161} (2005), 1--59.
\bibitem{BH2} M. Bestvina and M. Handel, Train tracks and
automorphisms of free groups, \textit{Ann. of Math. (2)}, {\bf 135}
(1992), 1--51.
\bibitem{BMMV} O. Bogopolski, A. Martino, O. Maslakova and
E. Ventura, Free-by-cyclic groups have solvable conjugacy
problem, preprint.
\bibitem{steer} M.R. Bridson, The geometry of the word problem, in \textit{Invitations to geometry and topology} (M.R. Bridson and S.M. Salamon, eds), Oxford University Press, 2002.
\bibitem{BG1} M.R. Bridson and D.P. Groves, The quadratic
isoperimetric inequality for mapping tori of free group automorphisms
I: Positive automorphisms, preprint at  http://arxiv.org/math.GR/0211459.
\bibitem{BG2} M. R. Bridson and D. Groves, Free-group automorphisms,
train tracks, and the beaded decomposition, preprint at
http://arxiv.org/math.GR/0507589.
\bibitem{BGgrowth} M.R. Bridson and D. Groves, The growth of conjugacy
classes under free group automorphisms, in preparation.
\bibitem{BH} M.R. Bridson and A. Haefliger, \textit{Metric spaces of 
non-positive curvature}, Springer-Verlag, Berlin, 1999. 
\bibitem{BrinkDyn} P. Brinkmann, Dynamics of free group
automorphisms, preprint.
\bibitem{Cooper} D. Cooper, Automorphisms of free groups have finitely generated fixed point sets, \textit{J. Algebra}, {\bf 111} (1987), 453--456.
\bibitem{LS} R.C. Lyndon and P.E. Schupp, Combinatorial group theory, 
Springer-Verlag, Berlin, 1977. 
\bibitem{OS} A.Yu. Ol'shanskii and M.V. Sapir, Groups with small
Dehn functions and bipartite chord diagrams, \textit{GAFA}, to appear.
\bibitem{schleimer} S. Schleimer, Polynomial time word problems, preprint.
\end{thebibliography}
\end{document}